\documentclass[12pt]{article}
\usepackage{graphicx,psfrag}
\usepackage{amsmath,amsthm,amssymb}
\usepackage{amsfonts}
\newtheorem{intro}{Theorem}[section]
\newtheorem{lintro}[intro]{Lemma}

\newtheorem{rintro}[intro]{Remark}

\newtheorem{pintro}[intro]{Proposition}

\newtheorem{dintro}[intro]{Definition}
\newtheorem{theorem}{Theorem}[subsection]
\newtheorem{definition}[theorem]{Definition}

\newtheorem{remarks}[theorem]{Remarks}
\newtheorem{remark}[theorem]{Remark}
\newtheorem{lemma}[theorem]{Lemma}
\newtheorem{cor}[theorem]{Corollary}
\newtheorem{proposition}[theorem]{Proposition}

\newcommand {\N}{\mathbb{N}} 
\newcommand {\Z}{\mathbb{Z}}            
\newcommand {\R}{\mathbb{R}} 
\newcommand {\E}{\mathbb{E}} 
\newcommand {\sph}{\mathbb{S}} 
\newcommand {\disc}{\mathbb{D}} 
\newcommand {\Q}{\mathbb{Q}} 
\newcommand {\C}{\mathbb{C}} 
\newcommand {\hip}{\mathbb{H}} 
\begin{document}
\makeatletter
\title{ Filling in solvable groups and in lattices in semisimple groups }
\author{
Cornelia DRU\c{T}U\thanks{drutu@agat.univ-lille1.fr}\\ \\
  UFR de Math\'ematiques et UMR 8524 au CNRS,\\
  Universit\'e de Lille I,
  F--59655 Villeneuve d'Ascq, France }
\date{}
\maketitle
\makeatother

\begin{abstract}
\noindent  We prove that the filling order is quadratic for a large class 
of 
solvable groups and asymptotically quadratic for all 
$\Q$-rank one 
lattices in semisimple groups of $\R$-rank at least $3$. As a byproduct of auxiliary results we give a shorter proof of the theorem on the nondistorsion of horospheres providing also an estimate of a nondistorsion constant. 
\end{abstract}

\section{Introduction}

In this paper we give an estimate of the filling order in some 
particular cases 
of infinite finitely generated groups and of Lie groups. We may talk about 
filling area of a loop in Riemannian manifolds, in 
finitely 
presented groups and, more generally, in metric spaces. In a metric space we 
begin by fixing a small $\delta $. By $\delta $-filling area of a loop one means the 
minimal number of small loops of length at most $\delta $ 
(``bricks'') one has to put one next to the other in order to obtain a net 
bounded by the given loop. Usually we choose $\delta =1$. By means of the filling area one can define the filling function and the filling order in a metric space. The filling 
order is the order in $\ell $ of the maximal 
area needed to fill a loop of length $\ell$ (see 
Section 
\ref{fil} 
for definitions and details). With the terminology introduced in Section \ref{fil}, if in a metric space $X$ a function in the same equivalence class as the filling function is smaller that $\ell ,\; \ell^2 $ or $ e^\ell $, it is sometimes said that the space $X$ {\it satisfies a linear, quadratic or exponential isoperimetric inequality}.

It is 
interesting to study these notions for two 
reasons at least : because a finitely presented group $\Gamma $ has a solvable Word problem if and only if its filling function is (bounded by) a recursive function and because the filling 
order is a 
quasi-isometry 
invariant \cite{Al*}.

We recall that a {\it quasi-isometry} between two metric spaces 
$(X_1,d_1)$ and 
$(X_2,d_2)$ is a map $q:X_1\to X_2$ such that 
$$
\frac{1}{L}d_1(x,y)-c\leq d_2(q(x),q(y))\leq Ld_1(x,y)+c,\; \forall x,y \in 
X_1,
$$ 
for some fixed positive constants $c$ and $L$, and $X_2$ is at a finite Hausdorff distance from the image of $q$. If such a map exists between $X_1$ and 
$X_2$, the 
two metric spaces are called {\it quasi-isometric}. A property invariant up to quasi-isometry is called {\it a geometric property.}

We shall study two classes of metric spaces : Lie solvable groups endowed with a left invariant Riemannian structure and lattices in 
semisimple 
groups, endowed with a word metric. The solvable groups are interesting as far as 
few 
things 
are known on their behaviour up to quasi-isometry. In the case of nilpotent groups there is more information. First, a consequence of Gromov's Theorem on polynomial growth \cite{Gr1*} is that virtual nilpotency is a geometric property in the class of groups (we recall that a group is called {\it virtually nilpotent} if it contains a nilpotent subgroup of finite index). In free nilpotent groups, the filling order is polynomial of degree $c+1$, where $c$ is the class of the group (\cite{Gr2*}, \cite{Pi*}). In the Heisenberg group ${\mathcal H}^3$ it was shown by Thurston in \cite{ECH*} that the filling order is cubic (which implies that  ${\mathcal H}^3$ is not automatic). The other Heisenberg groups ${\mathcal H}^{2n+1},\; n\geq 2$, have quadratic filling. This was conjectured by Thurston in \cite{ECH*}, Gromov gave an outline of proof in \cite{Gr2*} and D. Allcock gave the complete proof \cite{All*} by means of symplectic geometry. Y. Olshanskii and M. Sapir later gave a combinatorial proof \cite{OS*}.

The behaviour is much more diversified in the case 
of solvable groups. To begin with, the property of being virtually solvable is not  a geometric property in the class of groups anymore \cite{Di*}. On the other hand, certain solvable groups are very rigid with respect 
to quasi-isometry (\cite{FM1*}, \cite{FM2*}, \cite{FM3*}). Therefore the 
estimate of a quasi-isometry invariant for a class of solvable groups, as the 
one given in this paper, should be interesting. The filling order is already known for some solvable groups. W. Thurston has shown that the group $Sol$ has exponential filling (so it is not automatic) \cite{ECH*}. S. Gersten showed that the Baumslag-Solitar groups $BS(1,p)$ have exponential filling order for $p\neq 1$, which in particular implies that they are not automatic, while they are known to be asynchronously automatic. M. Gromov showed that $\R^n \rtimes \R^{n-1},\; n\geq 3,$ has polynomial filling order, without specifying the degree \cite{Gr2*}.
Also, G. Arzhantseva and D. V. Osin \cite{AO*} constructed a sequence of discrete non-polycyclic solvable groups with filling orders that are at most cubic. We should note here that the discrete solvable groups we deal with in this paper are all polycyclic, as lattices in Lie groups.

In the case of lattices, the filling order is already known for almost all ambient 
semisimple groups 
of $\R$-rank one (see the comments following Remark \ref{fact}). If the semisimple group has $\R$-rank 
$2$ then 
the filling order 
is exponential \cite{LP*}. We note that this has already been proven by W. Thurston in \cite{ECH*} in the particular case of $SL_3(\Z) $, from which result he deduced that $SL_3(\Z )$ is not combable. Also in \cite{ECH*} W. Thurston stated (without 
giving a 
proof) that the filling order of $SL_n(\Z),\; n\geq 4$, is quadratic. 
In the 
cases where the $\R$-rank is at least $3$ it is known that the filling order 
is at most 
exponential \cite{Gr2*}. In a 
previous paper we have proved that for some $\Q$-rank 
one 
lattices (among which the Hilbert modular groups) and for some 
solvable groups 
the filling order is at most ``asymptotically cubic'' \cite{Dr2*}. In this paper we 
prove the 
following more general and stronger results.

\begin{intro}\label{T1}
Let $X$ be a product of symmetric spaces and Euclidean buildings, $X$ of 
rank at 
least $3$. 
Let $\rho $ be a geodesic ray which is not contained in a rank one or in a rank two 
factor 
of $X$.
\begin{itemize}
\item[(1)] The 
Riemannian filling order in the horosphere $H(\rho )$ of $X$ is quadratic.
\item[(2)] Let $S$ be a Lie group acting by isometries, transitively with compact 
stabilizers on the horosphere $H(\rho )$ of $X$. In $S$ 
endowed with 
any 
left invariant metric the filling order is quadratic. The same is true for every discrete group $\Gamma$ acting properly discontinuously cocompactly on $H(\rho )$, $\Gamma$ endowed with a word metric.
\end{itemize} 
\end{intro}

The previous result has been obtained independently by E. Leuzinger and Chr. Pittet \cite{LP2*} in the case $X=SL_n(\R )/SO(n)$.

For a definition of Euclidean buildings see \cite{KlL*}, for a definition of 
horospheres see Section 
\ref{h}. We note that the quadratic estimate on the filling function is sharp.  This is because in a product of symmetric spaces and Euclidean buildings, which is of 
rank $r\geq 3$, every horosphere contains isometric copies of the Euclidean space $\E^{r-1}$. 
      
\begin{intro}\label{T2}
The filling order in any irreducible $\Q-$rank one lattice of a 
semisimple group 
of $\R$-rank at least $3$ is at most asymptotically quadratic. That is, for 
every 
$\varepsilon >0$ there exists $\ell_\varepsilon $ such that 
$$
A(\ell)\leq \ell^{2+\varepsilon },\; \forall \ell \geq \ell_\varepsilon \; .
$$
\end{intro}

\medskip

We emphasize that the filling order in all $\Q$-rank one lattices in 
semisimple groups of rank at least $2$ is at least quadratic. This is due to the fact that 
there are maximal flats in the symmetric space associated to the semisimple group on which the lattices act cocompactly. Thus Theorem \ref{T2} gives an ``asymptotically sharp'' estimate.

\medskip

{\it Examples where the two theorems apply} :

\medskip

(1) All solvable groups $Sol_{2n-1}(\alpha )=\R^n \rtimes_\alpha \R^{n-1},$ 
where $n\geq 3$ and $\alpha 
:\R^{n-1} \to Gl(n,\R)$ is an injective homomorphism with diagonalizable image, verify the 
hypothesis of 
Theorem \ref{T1}. The same is true for the groups of upper (lower) 
triangular 
matrices of order at least $4$ with the set of diagonals 
of the form $\{ (e^{\lambda_1}, e^{\lambda_2},\dots , e^{\lambda_n})\mid 
\sum_{i=1}^n \lambda_i =0 , \sum_{i=1}^n \alpha_i\lambda_i =0\}$ for a fixed 
vector $(\alpha_1 ,\alpha_2 ,\dots ,\alpha_n)$ with $\sum_{i=1}^n \alpha_i =0$.

\smallskip

(2) Any irreducible lattice in a semisimple group having a $\R$-rank 
one factor 
is a $\Q$-rank one lattice (\cite{Pr*}, Lemma 1.1). In particular the Hilbert modular groups $PSL(2,{\mathcal O}_K )$, where ${\mathcal O}_K$ is the ring of integers of a totally real field $K$ with $[K: \Q]\geq 3$.

 In \cite{Ti2*} one 
can also 
find examples of $\Q$-rank one lattices in simple 
groups.

\medskip

The paper is organized as follows. In Section \ref{prel} we recall some basic facts about filling area, asymptotic cones, buildings and ($\Q$-rank one) lattices in semisimple groups. In particular we recall that for every lattice there exists a space $X_0$ obtained 
from the 
ambient symmetric space $X=G/K$ by deleting a countable family of disjoint open horoballs on which 
the lattice acts with compact quotient (see 
Section \ref{r}). It follows that the filling order in the 
lattice is the same as the filling order in $X_0$.

In the case of the solvable groups we consider, one just has to look at horospheres. 
Since the 
projection of the exterior of the corresponding open horoball on the 
horosphere 
diminishes distances, one can study the whole exterior of the open 
horoball, 
which can also be denoted by $X_0$, instead of the horosphere.

Thus, both in the case of lattices and in the one of solvable groups it suffices to estimate the filling order in a metric space $X_0$ 
obtained 
from a symmetric space by deleting a family of disjoint open 
horoballs. The main 
tool we use is the asymptotic cone. This notion has been introduced by 
M. Gromov 
in \cite{Gr1*} and consists, philosophically speaking, of giving ``an image seen 
from infinitely far away''of a 
metric space (see 
Section \ref{fil} 
for definition and properties). To every metric space one associates a whole 
class of asymptotic cones (possibly isometric). There are similarities of 
arbitrary 
factor ``acting'' on this class, that is, sending a cone into another (Remark 
\ref{si}).

We dispose of a result, due to M. Gromov and P. Papasoglu, allowing to 
deduce from 
an 
uniform estimate of the filling order in all asymptotic cones an 
estimate of the filling order for 
the 
initial metric space (Section \ref{fil}, Theorem \ref{crit}). Thus, instead of 
considering 
the space $X_0$ one can consider its asymptotic cones. Let ${\mathbf K}_0$ be such an asymptotic cone. 
It is not 
difficult to prove that each ${\mathbf K}_0$ is obtained from an Euclidean 
building 
${\mathbf K}$ by 
deleting a family of disjoint open horoballs.

In Section \ref{horo} we first place ourselves in an Euclidean building. Essentially, an 
Euclidean 
building is a bunch of ``flats'', that is, of isometric copies of an 
Euclidean 
space. By deleting disjoint horoballs one makes polytopic 
holes into these 
flats. These polytopic holes can take, up to similarity, only a finite 
number of 
shapes. Thus we may hope to reduce the problem of filling a loop in a space like ${\mathbf K}_0$ to that of filling a loop in 
an 
Euclidean 
space with polytopic holes. We first prove some global and local properties of such a polytopic hole, that is, of the trace of a horoball in a maximal flat. We also provide a way to join two points on a connected horosphere by polygonal lines with length comparable to the distance, in two different cases (see Lemma \ref{fetze} and Lemma \ref{fetze2}). As a byproduct we give new proofs of the theorems on the nondistorsion of horospheres in Euclidean buildings and symmetric spaces. (Theorems \ref{nondist1} and \ref{nondist2} in this paper).

In Section \ref{k} we prove the main theorem, Theorem \ref{fillgen}, on the filling order in the exterior of a disjoint union of open horoballs. We give here an outline of the proof. First it is shown that a loop contained in the exterior of one open horoball, $X\setminus Hbo(\rho )$, where $\rho $ is not parallel to any rank 2 factor, has a quadratic filling area, as soon as the loop is contained in only one apartment (Proposition \ref{polit3}). Then we show that under the same hypothesis the same conclusion is true for certain loops contained in the union of two apartments (Proposition \ref{camera}). Then we prove Theorem \ref{fillK0}, which is a weaker version of Theorem \ref{fillgen}. This is done as follows. Up to spending a linear filling area, a generic loop can always be assumed to be contained in the union of a finite set of maximal flats, the number of maximal flats being of the same order as the length $\ell $ of the loop. In filling the loop, there is a problem when one passes 
from one flat 
to another. More precisely, by means of Proposition \ref{camera} one can show that when passing 
from one flat 
to another, in the process of filling the loop, one might have to spend an area of order $\ell^2$. This explains why instead of a quadratic filling area, one 
obtains a cubic 
filling area, in this first approach. For certain curves contained in a union of maximal flats of uniformly bounded cardinal a quadratic filling area is obtained.

In Section \ref{qsol} it is shown that part (b) of Theorem \ref{fillK0} implies part (b) of Theorem \ref{fillgen} and in particular Theorem \ref{T1}. The proof of the previous implication is done in two steps. Firstly it is shown that loops composed of a uniformly bounded number of minimizing almost polygonal curves have quadratic filling area (see the end of Section \ref{horo} for a definition of minimizing almost polygonal curves). Secondly an induction procedure is applied.

Section \ref{aq} contains the proof that Theorem \ref{fillK0}, (a), implies Theorem \ref{fillgen}, (a), and Theorem \ref{T2}. According to Theorem \ref{crit}, in order to prove Theorem \ref{fillgen}, (a), it is enough to prove that the filling order is quadratic in all the asymptotic cones of the space $X_0$ on which a $\Q$-rank one lattice $\Gamma $ acts properly discontinuously cocompactly. By Theorem \ref{fillK0}, (a), it is already known that in every asymptotic cone ${\mathbf K}_0$ of $X_0$ the filling order is at most cubic. It is shown that in reality the filling order is not cubic but quadratic.

To understand the proof of the previous statement we should see first why the filling order in the Euclidean plane is quadratic. In the Euclidean plane, any loop $\frak C$ of length $\ell$ 
can 
be filled with at most $\left( \frac{\ell }{\lambda_1}\right)^2$ bricks of 
length $\lambda_1$ (one may think of the bricks as being small squares). To fill it with bricks of 
length  $\lambda_2 << \lambda_1$, it is enough if we fill the  
$\lambda_1$-bricks with  $\lambda_2$-bricks, which can be done with $\left( 
\frac{\lambda_1 }{\lambda_2}\right)^2$ bricks for each $\lambda_1$-brick. In this way the initial loop $\frak C$ is filled with at most $\left( \frac{\ell }{\lambda_2}\right)^2$ bricks of 
length $\lambda_2$. Thus, the fact that for every 
$\lambda$, however small, we may fill $\frak C$ with at most $\left( \frac{\ell 
}{\lambda}\right)^2$ $\lambda $-bricks is due to the fact that the quadratic 
filling is preserved in the small. This should also happen, under replacement 
of 
the exponent $2$ by the exponent $3$, in a space with a cubical filling order. But in the space ${\mathbf K}_0$ the first 
important 
remark is that, when one fills a loop of length $\ell$ with bricks of length 
$1$, one puts $k_2\ell^2$ bricks with uncontrolled shapes and $k_1\ell^3$ bricks 
which bound small Euclidean squares entirely contained in ${\mathbf K}_0$ (Remark \ref{rc}). This 
means 
that, for $\lambda < 1$, the $\lambda $-filling area of the loop tends to become more and more quadratic as $\lambda $ becomes smaller and smaller. Since we may choose bricks as small as we want, the quadratic factor 
will end by dominating the cubic factor. And then, by applying similarities 
(Remark \ref{si}), since this is a 
reasoning 
which is done simultaneously in all asymptotic cones, one can come back to 
bricks of length one and obtain a quadratic filling order.

In the Appendix we provide an isoperimetric inequality for a hypersurface in an Euclidean space composed of points at a fixed distance from a certain polytope ${\mathcal P}$ (Proposition \ref{lpolit}). From this we derive an isoperimetric inequality for every polytopic hypersurface whose points are at a distance between $R>0$ and $a R,\; a>1$, from the polytope ${\mathcal P}$. The constants appearing in the isoperimetric inequality in the first case depend only on ${\mathcal P}$ while in the second case they also depend on $a$. The second result is useful in the proof of our main theorem.

{\it Notations} : Throughout the whole paper, in a metric space $X$, $B(x,r)$ denotes the open ball of center $x\in X$ and radius $r>0$, $S(x,r)$ its boundary sphere, ${\mathcal{N}}_r(A):=\{ x\in X \mid d(x,A)\leq r \}$, $\partial {\mathcal{N}}_r(A):=\{ x\in X \mid d(x,A) = r \}$ and $\breve{{\mathcal{N}}}_r(A):=\{ x\in X \mid d(x,A)<r \}$, where $A\subset X$.
 
\smallskip 

\noindent {\bf Acknowledgments.} The main part of this paper was written during the author's stay at Max-Planck-Institut f\"ur Mathematik in Bonn. I would like to thank this institution for its hospitality. I would also like to thank Mark Sapir for having shown me the induction trick, which I used in the proof of the implication Theorem \ref{fillK0}, (a) $\Longrightarrow $ Theorem \ref{fillgen}, (a).

\section{Preliminaries}\label{prel}

\subsection{Filling area, filling order, asymptotic cone}\label{fil}

The notion of filling area of a loop is well defined in the setting of Riemannian manifolds as well as in finitely presented groups (see for instance \cite[Chapter I, $\S 8A.4$]{BH*}). In the sequel we recall the meaning of this notion in geodesic metric spaces. Let $X$ be such a space and $\delta >0$ a fixed constant. 
We call ``loops'' lipschitz maps $\frak C$ from $\sph^1$ to $X$. We 
call {\it 
filling partition of} $\frak C$ a pair consisting of a triangulation of the 
planar unit 
disk $\disc^2$ and of an injective map from the set of 
vertices of 
the triangulation to $X,\; \pi :{\mathcal V}\to X$, where $\pi$ coincides 
with $\frak C$ 
on ${\mathcal V} \cap \sph^1$. The image of the map $\pi $ is called {\it 
filling disk of} $\frak C$. We can join the images of the vertices of each triangle by geodesics (for the 
vertices which are ends 
of arcs of 
$\sph^1$, we replace the geodesic by the arc of $\frak C$ contained 
between their images). We call the geodesic triangles thus obtained {\it bricks}. The length of a brick is the sum of the distances between 
vertices (for the vertices which are ends 
of arcs of 
$\sph^1$, we replace the distance by the length of the arc of $\frak C$ contained 
between the two 
images). The maximum of the lengths of bricks 
in a 
partition is called {\it the mesh of the partition}. The partition is called 
$\delta$-{\it 
filling partition of} $\frak C$ if its mesh is at most $\delta$. The corresponding filling disk is called $\delta$-{\it 
filling disk of} $\frak C$. We call 
$\delta$-{\it filling area of} $\frak C$ the minimal number of triangles in a triangulation associated to a 
$\delta$-filling 
partition of $\frak C$. We denote it with the double notation $A_\delta 
(\frak C)=P(\frak C,\delta 
)$. This generalisation of the notion of area is due to M. Gromov 
\cite{Gr2*}.

In each of the three cases (Riemannian manifolds, finitely presented groups, geodesic metric spaces) when we have defined a notion of 
``filling area'' 
for loops, we can now define the {\it filling function} $A:\R_+^* \to \R_+^*$, $A(\ell ):=$ the maximal area 
needed to 
fill a loop of length at most $\ell$. In metric spaces, 
for $\delta $ fixed, we use the notation $A_\delta (\ell )$ and we call this function the $\delta $-{\it filling function}.

Two filling functions corresponding to different presentations 
in a group, or to different $\delta $ in a metric 
space, or, more generally, to quasi-isometric metric spaces, satisfy an equivalence relation. We define this equivalence relation below.

Let $f_1$ and $f_2$ be two functions of real variable, taking real values. We say that {\it the order of the function $f_1$ is at most the order of the function $f_2$}, and we denote it by $f_1\prec f_2$, if $f_1(x)\leq af_2(bx+c)+dx+e ,\; \forall x$, where $a,\, b,\, c,\, d,\, e$ are fixed positive constants. We say that {\it $f_1$ and $f_2$ have the same order}, and we denote it by $f_1\circeq f_2$, if $f_1\prec f_2$ and $f_2\prec f_1$. The relation $\circeq $ is an equivalence relation. The equivalence class of a numerical function with respect to this relation is called {\it the order of the function}. If a function $f$ has (at most) the same order as the function $x,\, x^2,\, x^3,\, x^d$ or $\exp x$ it is said that {\it the order of the function $f$ is (at most) linear, quadratic, cubic, polynomial, or exponential,} respectively.

The order of the filling function of a metric space $X$ is also called {\it the filling order of $X$}.

\bigskip

{\it A non-principal ultrafilter} is a 
finitely 
additive measure $\omega $ defined on all subsets of $\N$, taking as 
values 
$0$ and $1$ and taking always value $0$ on finite sets. Such a measure 
always 
exists \cite{Dr1*}. For a sequence in a topologic space, $(a_n)$, one can 
define 
{\it the $\omega$-limit} as being the element $a$ with the property 
that for 
every neighborhood ${\mathcal N}(a)$ of $a$, the set $\lbrace n\in \N 
\mid a_n 
\in 
{\mathcal N}(a)\rbrace$ has $\omega $-measure $1$. We denote $a$ by 
$\lim_\omega 
a_n$. 
Any sequence in a compact space has an $\omega $-limit 
\cite[I.9.1]{Bou*}, and 
the $\omega $-limit is unique.

Let $(X,d)$ be a metric space. We fix a sequence $(x_n)$ of points in 
$X$, which 
we call {\it sequence of observation centers}, a sequence of positive 
numbers 
$(d_n)$ diverging to infinity, which we call {\it sequence of 
scalars}, and a 
non-principal ultrafilter $\omega $. Let ${\mathcal C}$ be the set of 
sequences 
$(y_n)$ of points in $X$ with the property that 
$\frac{d(x_n,y_n)}{d_n}$ is 
bounded. We define an equivalence relation on ${\mathcal C}$ :
$$
(y_n)\sim (z_n) \Leftrightarrow \lim_\omega \frac{d(y_n,z_n)}{d_n} =0\, .
$$

The quotient space of ${\mathcal C}$ with respect to this relation, 
which we 
denote by $X_\omega(x_n,d_n)$, is called {\it the asymptotic cone of 
}$X$ with 
respect to the observation centers $(x_n)$, the scalars $(d_n)$ and 
the 
non-principal ultrafilter $\omega $. It is a complete metric space 
with the 
metric 
$$
D(\lbrack y_n\rbrack ,\lbrack z_n \rbrack )=\lim_\omega 
\frac{d(y_n,z_n)}{d_n} 
\; .
$$

We say that the set $A\subset X_\omega (x_n,d_n)$ is {\it the limit set of the 
sets} 
$A_n\subset X$ if 
$$
A=\lbrace \lbrack x_n\rbrack \mid x_n\in A_n 
\; \; \omega-\mbox{ 
almost surely }\rbrace \; .
$$
We denote $A= \lbrack A_n\rbrack $.

In our arguments we shall use the following very simple but important 
remark.

\begin{remark}\label{si}
The map 
$$
I_\alpha:X_\omega(x_n,d_n)\to X_\omega \left( x_n,\frac{1}{\alpha}d_n\right) 
,\; \; I_\alpha( 
\lbrack x_n\rbrack )=\lbrack x_n\rbrack
$$
 is a similarity of factor $\alpha $. 
\end{remark} 

There is a relation between the filling order in the asymptotic cones 
and the 
filling order in the initial space, established by P. Papasoglu ( see \cite[Theorem 2.7]{Dr2*}), who adapted an idea of M. Gromov for this purpose.

\begin{theorem}[P. Papasoglu]\label{crit}
Let $X$ be a metric space. If in every asymptotic cone of $X$ we have that 
$$
A_1(\ell ) \leq C\cdot \ell^p\, ,\forall \ell ,
$$ where $C$ is an universal constant, then in the space $X$ we have 
that for 
every $\varepsilon >0$ there exists $\ell_\varepsilon $ such that 
$$
A_1^X(\ell ) \leq \ell^{p+\varepsilon } \, 
,\forall \ell \geq \ell_\varepsilon\, . 
$$
\end{theorem}

\subsection{Spherical and Euclidean buildings}\label{im}

For this section we refer mainly to \cite{KlL*}, but also to \cite{Ti1*} 
and 
\cite{Dr1*}.

Before discussing about buildings we introduce some 
terminology in Euclidean spaces. For a subset $A$ in an Euclidean space $\E^k$ 
we call {\it affine span} of $A$, and we denote it with $Span\; A$, the minimal affine subspace of $\E^k$ 
containing $A$. Two polytopes of codimension one are called {\it 
parallel} if their affine spans are parallel. A subset $A$ of $\E^k$ is called {\it relatively open} if it is open in $Span\; A$. The {\it relative interior} of a subset $B$ of $\E^k$ is its interior in $Span\; B$.

In an Euclidean sphere $\sph^k$, we call {\it spherical span} of a subset $\sigma$ the trace on the sphere of the affine span of the cone of vertex the origin over $\sigma$. We denote it by $Span\; \sigma$. We say that two subsets $\sigma $ and $\sigma'$ in $\sph^k$ are 
{\it orthogonal to each other} if $Span\; \sigma $ and $Span\; \sigma'$ are orthogonal. In a spherical building  
two subsets $\sigma $ and $\sigma'$ are {\it orthogonal} if they are both contained in the same apartment and are orthogonal. We can also define the distance between a point and a convex set in a spherical building each time they are both contained in an apartment as the spherical distance between them in that apartment.

Let $\Sigma $ be a spherical building. Throughout the whole paper we shall 
suppose that for all spherical buildings, the associated Weyl group 
acts on the associated Coxeter complex without fixed points. In this way we rule out the case of spherical buildings having a sphere as a factor and of Euclidean buildings and symmetric spaces having an Euclidean space as a factor.

An apartment is split by each singular hyperplane in it into two 
halves called {\it half-apartments}. All simplices in $\Sigma 
$ which are not chambers are called {\it walls}, the codimension one simplices are also called {\it panels}. Intersections of singular hyperplanes in an apartment  are called {\it singular subspaces}. A chamber is said to be {\it adjacent} to a singular 
subspace if their intersection is a wall of the same dimension as the subspace. 
We also say, in the previous situation, that the singular subspace {\it 
supports} the chamber.

 We say that a panel {\it separates} a chamber and a point if there is at least one apartment containing the three of them and in each such apartment the spherical span of the panel separates the point and the interior of the chamber.

Two chambers are said to be {\it adjacent} if they have a panel in common and disjoint interiors. A {\it gallery} of chambers is a sequence of chambers such that any two consecutive chambers are adjacent. The number of chambers composing it is called {\it the length of the gallery}. Given a point and a simplex, {\it the combinatorial distance} between them is the minimal length of a gallery of chambers such that the first chamber contains the point and the last contains the simplex. For every such gallery of minimal length between the point and the simplex, the last of its chambers is called {\it the projection of the point on the simplex}. We note that a point may have several projections on a panel if and only if the point and the panel are contained in a singular hyperplane.

For every wall ${\mathcal{M}}$ in the spherical building we call {\it star of} ${\mathcal{M}}$, and we denote it by $Star\; ({\mathcal{M}})$, the set of chambers containing ${\mathcal{M}}$. We use the same name and notation for the union of all the chambers containing ${\mathcal{M}}$. A building is called {\it $c$-thick} if for every panel ${\mathcal{P}}$, the cardinal of $Star\; ({\mathcal{P}})$ is at least $c$.

Every spherical building $\Sigma $ admits a labelling 
\cite[IV.1, 
Proposition 
1]{Br*}. With respect to this labelling one can define a projection of 
the 
building on the model spherical chamber $p:\Sigma \to \Delta_{mod}$. Given a subset ${\mathcal{C}}$ in $\Sigma $ its image under this projection, $p({\mathcal{C}})$, is called {\it the set of slopes of} ${\mathcal{C}}$.

{\it The model Coxeter complex of }$\Sigma $ is the
Coxeter complex ${\mathsf{S}}$ which is isomorphic to any of its apartments ;
for every labelling in  the spherical building there is a
compatible labelling on the model Coxeter  complex and a
compatible projection $p_{\mathsf{S}}:{\mathsf{S}}\to \Delta_{mod}$.

Given an apartment $\mathcal{A}$ in $\Sigma $ and a chamber $\mathcal{W}$ in it one can always define a map retr$_{\mathcal{A}, \mathcal{W}} : \Sigma \to \mathcal{A}$ which preserves labelling and combinatorial distances to $\mathcal{W}$ and diminishes the other combinatorial distances \cite[$\S 3.3-3.6$]{Ti1*}. This map is called {\it the retraction of $\Sigma $ onto $\mathcal{A}$ with center $\mathcal{W}$}.

\begin{definition}\label{ort}
Let $\mathsf{S}$ be a Coxeter complex with a labelling, $\Delta_{mod}$ its
model spherical chamber and $p_{\mathsf{S}} : \mathsf{S}\to \Delta_{mod} $ the projection
corresponding to the labelling.  Let
$\theta \in \Delta_{mod} $ be a given point. {\it The set of orthogonals
to }$\theta $ is the set of points $q\in \Delta_{mod}$ such that any
preimage of $q$ by $p_{\mathsf{S}}$ is at distance $\frac{\pi}{2}$ from a
preimage of $\theta $ by $p_{\mathsf{S}}$. We denote this set by $Ort\;
(\theta )$.   
\end{definition}

\smallskip

We note that if $\mathsf{S}$ is the model Coxeter complex of a spherical building $\Sigma $, endowed with a labelling compatible with the one of $\Sigma $, $Ort\; (\theta )$ also coincides with the set of points $q\in \Delta_{mod}$ such that any
preimage of $q$ by $p:\Sigma \to \Delta_{mod}$ is at distance $\frac{\pi}{2}$ from a
preimage of $\theta $ by $p$.

\medskip

If ${\bf S}$ and 
${\bf S}'$ are two simplicial complexes, their {\it join}, denoted by ${\bf 
S}\circ {\bf S}'$, may be defined abstractly to be a simplicial complex having 
as vertex set the disjoint union of the sets of vertices of ${\bf S}$ and of 
${\bf S}'$ and having one simplex $\sigma \circ \sigma'$ for every pair of 
simplices $\sigma \in {\bf S}$ and $\sigma'\in {\bf S}'$. We also allow the 
possibility for one of the two complexes to be empty, and we make the convention 
that ${\bf S}\circ \emptyset ={\bf S}$.

There is a geometric interpretation of the join for spherical complexes which 
goes as follows. Let ${\bf S}$ and ${\bf S}'$ be spherical simplices in the 
Euclidean spheres $\sph^{k-1}\subset \E^k$ and $\sph^{m-1}\subset \E^m$, respectively. 
The join ${\bf S}\circ {\bf S}'$ is the spherical simplex in $\sph^{k+m-1}$ obtained by embedding ${\bf S}$ together with $\sph^{k-1}$ 
into $\sph^{k+m-1}$ and likewise ${\bf S}'$ together with $\sph^{m-1}$, 
such that the embeddings of $\sph^{k-1}$ and $\sph^{m-1}$ are orthogonal, and considering 
the convex hull of ${\bf S}\cup {\bf S}'$ in $\sph^{k+m-1}$. One may also say that ${\bf S}\circ {\bf S}'$ 
is obtained from ${\bf S}$ and ${\bf S}'$ by gluing the extremities of a quarter 
of a circle to every pair of points $x\in {\bf S}$ and $x'\in {\bf S}'$. Since 
this also makes sense for two spherical complexes and in particular for 
spherical buildings, we thus get a geometric definition of the join $\Sigma 
\circ \Sigma '$ of two spherical buildings $\Sigma $ and $\Sigma '$. B. Kleiner 
and B. Leeb proved in \cite[$\S 3.3$]{KlL*} the following :

\smallskip

$\bullet $ every decomposition as a join of the model chamber of a spherical 
building, $\Delta_{mod} =\Delta_1 \circ \Delta_2 \circ \dots  \circ \Delta_n$, 
or of the associated Coxeter complex $\mathsf{S} =\mathsf{S}_1 \circ \mathsf{S}_2 \circ \dots  \circ \mathsf{S}_n$ 
imply a decomposition of the spherical building $\Sigma =\Sigma_1 \circ \Sigma_2 
\circ \dots  \circ \Sigma_n$ such that $\Delta_i$ and $\mathsf{S}_i$ are the model 
chamber and the associated Coxeter complex of $\Sigma_i$. 

\smallskip

$\bullet $ a spherical building $\Sigma $ is not a join of two nonempty spherical buildings if and 
only if its model chamber $\Delta_{mod}$ has diameter $<\frac{\pi }{2}$ and 
dihedral angles $\leq \frac{\pi }{2}$.

\begin{lemma}\label{fnort}
Let $\Sigma $ be a labelled spherical building and let $q$ be a point in it such 
that for every decomposition of $\Sigma $ as a join, $\Sigma =\Sigma_1 \circ 
\Sigma_2$, $q$ is contained neither in $\Sigma_1$ nor in $\Sigma_2$. Then in every 
chamber containing it, $q$ is not orthogonal to any wall of the chamber.
\end{lemma}

\noindent{\bf Proof.}\quad Let $\sigma $ be 
the unique wall containing $q$ in its interior. Suppose $q$ is orthogonal to 
another wall $\sigma'$ such that $\sigma $ and $\sigma '$ are both contained in 
a chamber $\Delta $. This implies that the diameter of 
$\Delta $ and of the model chamber $\Delta_{mod}$ is $\pi /2$.

If $\sigma $ and $\sigma '$ intersect in a wall, this wall 
being orthogonal to an interior point of $\sigma $ it follows that the diameter 
of $\sigma $ is $>\frac{\pi}{2}$. This is impossible, as the diameter of $\Delta 
$ is at most $\frac{\pi}{2}$. So we may suppose that $\sigma $ and $\sigma '$ do 
not intersect. Let $\mathfrak m$ be the convex hull of $\sigma $ and $\sigma '$ 
in $\Delta $. Since they are orthogonal, ${\mathfrak m} =\sigma \circ \sigma '$.

The building $\Sigma $ is decomposable as a join because $\Delta_{mod}$ has 
diameter $\frac{\pi }{2}$. The maximal decomposition of $\Sigma $ as a join, 
$\Sigma = \Sigma_1 \circ \Sigma_2 \circ \dots  \circ \Sigma_n$, induces a 
maximal decomposition of $\Delta $ as a join, $\Delta = \Delta_1 \circ \Delta_2 
\circ \dots  \circ \Delta_n $, which in its turn induces a decomposition of 
${\mathfrak m}$, ${\mathfrak m}= {\mathfrak m}_1 \circ {\mathfrak m}_2 \circ 
\dots  \circ {\mathfrak m}_n$ (where some of the ${\mathfrak m}_i$ might be 
empty). It follows that $\sigma = {\mathfrak m}_{i_1} \circ {\mathfrak m}_{i_2} \circ 
\dots  \circ {\mathfrak m}_{i_s}$ and $\sigma '= {\mathfrak m}_{j_1} \circ 
{\mathfrak m}_{j_2} \circ \dots  \circ {\mathfrak m}_{j_t},\; \{i_1,i_2,\dots 
,i_s \}\sqcup \{j_1,j_2,\dots ,j_t \}=\{1,2,\dots ,n \}$. This implies that 
$q\in \sigma \subset \Sigma_{i_1} \circ \Sigma_{i_2} \circ \dots  \circ 
\Sigma_{i_s}$, which contradicts the hypothesis.\hspace*{\fill } $\diamondsuit $

\medskip

The following simple remark enlightens us more on the geometry of spherical 
buildings.

\begin{remark}\label{discret}
Let $\theta $ be a point in $\Delta_{mod} $. The set $D(\theta 
)=\{d(x,\sigma )\mid  p(x)=\theta,\; \sigma \mbox{ simplex }\}$ is finite and contains three 
consecutive terms of type $\frac{\pi }{2} -\delta_0 ,\frac{\pi 
}{2}, \frac{\pi }{2} +\delta_0' $  or two consecutive terms of type $\frac{\pi 
}{2} -\delta_0 , \frac{\pi }{2} +\delta_0' $, where $\delta_0 $ and $\delta_0' $ 
depend 
only on $\theta $. 
\end{remark}

\medskip

{\it Notations :} Let $X$ be a CAT(0)-space and its boundary at 
infinity 
$\partial_\infty X$. For every $x\in X$ and $\alpha \in 
\partial_\infty X$, we 
denote by $\lbrack x,\alpha )$ the unique ray having $x$ as origin and 
$\alpha $ as 
point at infinity. For two geodesic segments or rays $\lbrack x,a )$ 
and 
$\lbrack x,b )$ we denote by $\angle_x(a,b)$ the angle between them in 
$x$ (see 
\cite{KlL*} for a definition). We denote the Tits metric on 
$\partial_\infty X$ by 
$d_T$. For every point $x$ and every geodesic ray $\rho $ we denote by $\rho_x$ the ray of origin $x$ asymptotic to $\rho $.

Let now $X$ be a symmetric space or an Euclidean building or a product of a symmetric space with an Euclidean building, $X$ of rank $r$. In the sequel, for simplicity, we call the apartments in Euclidean buildings also maximal flats. For definitions and results in symmetric space theory and in Euclidean building theory
we refer to \cite{He*}, \cite{BH*} and \cite{KlL*}. We only  recall that 

\smallskip

$\bullet $ $m$-flats are isometric copies of the Euclidean 
space 
$\E^m,\; m\leq r$ ; singular $m$-flats(subspaces) are $m$-flats which appear as intersections 
of apartments ; we also 
call 
singular $(r-1)$-flats singular hyperplanes ;

\smallskip

$\bullet $ half-apartments are halves of apartments determined by singular 
hyperplanes ;   

\smallskip

$\bullet $ the faces of the Weyl chambers are called walls ; the codimension 1 walls are also called panels ;

\smallskip

$\bullet $ two Weyl chambers are called adjacent if they have the vertex and a panel in common and disjoint interiors ;

\smallskip

$\bullet $ a singular subspace $\Phi $ is said to be adjacent to a Weyl chamber $W$ or to support $W$ if $\Phi \cap W$ is a wall of the same dimension as $\Phi $ ;

\smallskip

$\bullet $ a gallery of Weyl chambers is a finite sequence of Weyl chambers such that any two consecutive Weyl chambers are adjacent ; its length is the number of Weyl chambers ; a minimal gallery is a gallery of minimal length among the ones which have the same first and last Weyl chambers as itself ; 
 
\smallskip

$\bullet $ for every wall $M$ we define $Star\; (M)$, which we call {\it the star of $M$}, as the set of all Weyl chambers having the same vertex as $M$ and containing $M$ ; we use the same name and notation for the union of all these Weyl chambers ;  
 
\smallskip

$\bullet $ Weyl polytopes are polytopes which appear as intersections of 
half-apartments (they may have dimension smaller than $r$, as we may intersect 
opposite half-apartments) ;
 
\smallskip

$\bullet $ an Euclidean building is called $c$-thick if every singular 
hyperplane is the boundary of at least $c$ half-apartments of disjoint
interiors. 

\medskip

The boundary at infinity, $\partial_\infty X$, is a spherical building (\cite{Mo*}, chapters 15 and 
16, 
\cite{BGS*}, Appendix 5). The model Coxeter complex and chamber of $\partial_\infty X$ are sometimes called {\it the model Coxeter complex and chamber of $ X$}.

If $X$ is a $c$-thick Euclidean building then $\partial_\infty X$ is a $c$-thick spherical building. If $X$ decomposes as $X=X_1\times X_2$, then $\partial_\infty X =\partial_\infty X_1 \circ \partial_\infty X_2$. For every maximal flat $F$ in $X$ we denote its boundary at infinity (which is an apartment) with $F(\infty )$. We likewise denote $\rho(\infty),\, W(\infty ),\, \Phi (\infty )$ the boundary at infinity of a ray $\rho$, a Weyl chamber $W$, a singular subspace $\Phi $. If two maximal flats $F_1,\, F_2$ have $F_1(\infty )=F_2(\infty )$ then $F_1=F_2$. If $F_1(\infty )\cap F_2(\infty )$ is a half-apartment and $X$ is an Euclidean building, then $F_1\cap F_2$ is a half-apartment. A maximal flat $F$ is said to be {\it asymptotic to a ray }$\rho $ if $\rho (\infty )\in F(\infty )$.

If $X$ is an Euclidean building, one can 
define 
{\it the space of directions in a point }$x$, which is the space of 
equivalence 
classes of geodesic segments having $x$ as an endpoint with respect to 
the 
equivalence relation ``angle zero in $x$''. We denote it by 
$\Sigma_xX$, and we 
call its elements {\it directions in} $x$. We denote by $\overline{x 
a}$ the 
direction corresponding to the geodesic segment or ray $\lbrack x,a )$. Given a 
geodesic ray $\rho $, we denote by $\overline{\rho_x}$ the direction in $x$ of the ray of 
origin $x$ asymptotic to $\rho$. For every convex set $\mathcal{C}$ containing $x$ we define the {\it set of directions of $\mathcal{C}$ in $x$, $\mathcal{C}_x$,} as the set of directions $\overline{x a}$ corresponding to all $[x,a)\subset \mathcal{C}$.

With 
respect to 
the metric induced by the angle, $\Sigma_xX$ becomes a spherical 
building. If $X$ is a $c$-thick homogeneous Euclidean building then for every $x$, $\Sigma_xX$ is a $c$-thick spherical building.

 If a segment $[x,b)$ and a convex set $\mathcal{C}$ are both contained, near $x$, in the same apartment, then the distance between  $\overline{x b}$ and $\mathcal{C}_x$ in $\Sigma_xX$ is well defined. We denote it either by $\angle_x(\overline{x b}, \mathcal{C}_x)$ or by $\angle_x(b, \mathcal{C})$.

If $X$ is a homogeneous Euclidean 
building then for every $x$, $\Sigma_xX$ has the same model chamber as 
$\partial_\infty X$, $\Delta_{mod}$, so, with respect to some labelling, one can 
define a 
projection $p_x:\Sigma_x \to \Delta_{mod}$. Moreover one can choose a 
labelling on 
$\Sigma_xX$ compatible with the one on $\partial_\infty X$, that is, 
such that for 
every point at infinity, $\alpha $, we have $p(\alpha 
)=p_x(\overline{x 
\alpha})$. We call this common value {\it the slope of the ray 
}$\lbrack 
x,\alpha )$. This result implies that if $\lbrack x,y\rbrack $ is a non-trivial 
geodesic segment, for every $a\in \lbrack x,y\lbrack ,\; 
p_a(\overline{ay})=p_x(\overline{xy})$. We call this common value {\it the 
slope 
of the segment 
}$\lbrack 
x,y \rbrack $. We note that the slope of the segment $\lbrack 
x,y \rbrack $ is in general not the same as the slope of the segment $\lbrack 
y,x \rbrack $ (since, generically, two opposite points in a spherical building 
do not project on the same point of $\Delta_{mod}$, but on two points sent one onto the other by the opposition involution).

If $X$ is a product of symmetric spaces and Euclidean 
buildings, we call {\it the slope of the ray 
}$\lbrack 
x,\alpha )$ the image of $\alpha $ by the projection $p:\partial_\infty X \to \Delta_{mod}$.

We call the ray $\lbrack x,\alpha )$ {\it regular} ({\it 
singular}) 
if its slope is in $Int\; \Delta_{mod}$ ($\partial \Delta_{mod}$). If
$\theta$ is a slope in $X$ we also call its set of
orthogonals {\it set of
orthogonal slopes}.

For a convex set $\mathcal{C}$ in a Euclidean building, its {\it set of slopes} is the set of all slopes of all segments $[x,y]\subset \mathcal{C}$. It is a set invariant with respect to the opposition involution.

According to \cite[Proposition 4.3.1]{KlL*} every decomposition of $\Delta_{mod}$ as a join, $\Delta_{mod} =\Delta_1 \circ \Delta_2$, corresponds to a decomposition of $X$ as a product, $X=X_1\times X_2$ such that $\Delta_i$ is the model chamber of the factor $X_i,\; i=1,2$.

We say that a slope $\theta $ is {\it parallel to a factor of $X$} if $\Delta_{mod}$ decomposes nontrivially as a join $\Delta_{mod} =\Delta_1 \circ \Delta_2$ and $\theta \in \Delta_1$. Using the previous remark one can verify that a slope is parallel to a factor if and only if one(every) segment or ray of slope $\theta $ is contained in the copy of a factor of $X$.

\begin{lemma}\label{apgcw}
Let $D$ be a half-apartment in an Euclidean building ${\mathbf K}$ and let $W$ 
be a Weyl chamber with a panel in $\partial D$ and with interior disjoint from $D$. Then there 
exists an apartment containing both $D$ and $W$.
\end{lemma}

\noindent{\bf Proof.}\quad By \cite[Proposition 3.27]{Ti1*}, there is an 
apartment $A$ in $\Sigma_x{\mathbf K}$ containing both $W_x$ and $D_x$. The chamber $W_x$ has an opposite chamber $W_x'$ in $D_x$. Let $\rho $ be a regular ray in $W$ and $\rho'$ the 
opposite ray in $W'$. There exists a unique apartment $F$ containing the regular geodesic $\rho \cup \rho' $. It follows that it contains $W$ and $W'$, therefore also $\partial D$ which is the convex hull of the two opposite panels $W\cap \partial D$ and $W'\cap \partial D$. Since $D$ is the convex hull of $\partial D$ and of $W'$ we may conclude.\hspace*{\fill }$\diamondsuit $

\medskip

In particular, two adjacent Weyl chambers are always contained in an apartment, in an Euclidean building.

\begin{definition}
Let $F$ be an apartment in an Euclidean buiding. We say that another apartment 
$F'$ is a ramification of $F$ if either $F'=F$ or $F\cap F'$ is a 
half-apartment. If the case $F'=F$ is excluded, $F'$ is called a strict 
ramification of $F$.
\end{definition}

\begin{cor}\label{ramif}
Let $F$ be an apartment and $W$ a Weyl chamber adjacent to a Weyl chamber $W'\subset F$. Then there exists a ramification $F'$ of $F$ containing $W\cup W'$.
\end{cor}

\subsection{Horoballs and horospheres}\label{h}

Let $X$ be a CAT(0)-space and $\rho $ 
a geodesic 
ray in $X$. {\it The Busemann function associated to $\rho$ } is the 
function $ f_\rho:X\to \R ,\; f_\rho(x)=\lim_{t\to \infty}[d(x,\rho(t))-t]\; .$ This function is well defined and convex. Its level hypersurfaces 
$H_a(\rho ):=\lbrace x\in 
X 
\mid 
f_\rho(x)= a \rbrace$ are called {\it horospheres }, its level sets 
$Hb_a(\rho ):=\lbrace x\in X 
\mid 
f_\rho(x)\leq a \rbrace$ are called {\it closed horoballs }  and their 
interiors, $Hbo_a(\rho ):=\lbrace x\in X 
\mid 
f_\rho(x)< a \rbrace$, 
{\it 
open horoballs}. We use the notations $H(\rho ),\; Hb(\rho ),\; Hbo(\rho )$ for the 
horosphere, the closed and open horoball corresponding to the value $a=0$.

For two asymptotic rays, their Busemann functions differ by a 
constant. Thus the families of horoballs and horospheres are the same and we 
shall call them 
horoballs and 
horospheres {\it of basepoint }$\alpha$, where $\alpha $ is the common 
point at 
infinity of the two rays.

\begin{remarks}\label{prh}
Let $\rho $ be 
a geodesic 
ray in a complete CAT(0)-space $X$ and let $a<b$ be two real numbers.

(a) There is a natural projection $p_{ba}$ of $H_b(\rho )$ onto $H_a(\rho )$ which is a 
surjective contraction.

(b) The distance from a point $x\in H_b(\rho )$ to $p_{ba}(x)$ is $b-a$.  
\end{remarks}

\noindent{\bf Proof.}\quad For every $x\in H_b(\rho )$ it is enough to consider the 
ray $\rho_x$. The point $p_{ba}(x)$ is the intersection of $\rho_x$ with $H_a(\rho )$. For 
every point $y\in H_a(\rho )$ the ray $\rho_y$ may be extended to a geodesic by the 
completeness of $X$, and this geodesic intersects $H_b(\rho )$ in a unique point $x$. 
However the extension itself may not be unique. The property (a) follows by 
the convexity of the distance.\hspace*{\fill }$\diamondsuit $ 

\begin{lemma}\label{3}
Let  $X$ be a product of symmetric spaces and Euclidean
buildings and $\alpha_1 ,\alpha_2 ,\alpha_3$ three distinct
points in $\partial_\infty X$. If there exist three open
horoballs $Hbo_i$ of basepoints $\alpha_i,\; i=1,2,3$, which are mutually disjoint then $\alpha_1 ,\alpha_2 ,\alpha_3$ have the same
projection on the model chamber $\Delta_{mod} $. 
\end{lemma}

\noindent{\bf Proof.}\quad The proof is given in the proof of
Proposition 5.5 \cite{Dr1*}, step (b).

\subsection{Q-rank 1 lattices}\label{r}

A {\it lattice} in a Lie group $G$ is a discrete subgroup $\Gamma $ 
such that 
$\Gamma \backslash G$ admits a finite $G$-invariant measure. We refer 
to 
\cite{Ma*} or \cite{Ra*} for a definition of $\Q$-rank 1 lattices in 
semisimple 
groups. In the introduction we gave some examples of $\Q$-rank 1 
lattices.  In 
the sequel we list the two main properties of $\Q$-rank 1 lattices 
that we use. The first one relates the word metrics to the induced metric.

\begin{theorem}[Lubotzky-Mozes-Raghunathan, \cite{LMR1*}, 
\cite{LMR2*}]\label{lmr}
On any irreducible lattice of a semisimple group of rank at least $2$, 
the word 
metrics and the induced metric are bilipschitz equivalent.
\end{theorem}

By means of horoballs one can construct a subspace $X_0$ of the 
symmetric space 
$X=G/K$ on which the lattice $\Gamma $ acts with compact quotient.

\begin{theorem}[ \cite{Ra*}, \cite{Pr*}]\label{xo}
Let $\Gamma $ be an irreducible lattice of $\Q$-rank one in a 
semisimple group 
$G$. Then there exists a finite set of geodesic rays $\{ \rho_1,\rho_2,\dots 
,\rho_k \}$ such 
that the space $X_0 =X\setminus \bigsqcup_{i=1}^k\bigcup_{\gamma \in \Gamma} Hbo(\gamma 
\rho_i)$ has 
compact quotient with respect to $\Gamma$ and such that any two of the 
horoballs 
$Hbo(\gamma \rho_i)$ are disjoint or coincide.  
\end{theorem}

Lemma \ref{3} implies that 
$p(\lbrace 
\gamma \rho_i(\infty )\mid \gamma \in \Gamma ,\;  i\in
\lbrace1,2,\dots k\rbrace 
\rbrace)$ 
is only one 
point which we denote by $\theta $ and we call {\it the associated 
slope of 
}$\Gamma $ (we recall that $p$ is the projection of the boundary at 
infinity onto 
the model chamber). We have the following property of the associated 
slope :

\begin{proposition}[\cite{Dr1*}, Proposition 5.7]\label{dir}
If $\Gamma $ is an irreducible $\Q$-rank one lattice in a semisimple 
group $G$ 
of $\R$-rank at least $2$, the associated slope, $\theta$, is never
parallel to a factor of $X= G/K$. 
\end{proposition}

In particular, if $G$ decomposes into a product of rank one factors, 
$\theta $ 
is a point in $Int\; \Delta_{mod}$, that is the rays $\gamma \rho_i,\; i\in 
\lbrace 
1,2,\dots 
k\rbrace$, are regular.

Since the action of $\Gamma $ on $X_0$ has compact quotient,
$\Gamma 
$ with the 
word metric is quasi-isometric to $X_0$ with the length metric (the 
metric 
defining the distance between two points as the length of the
shortest curve 
between the two points). Thus, the asymptotic cones of $\Gamma $ are 
bilipschitz 
equivalent to the asymptotic cones of $X_0$. Theorem
\ref{lmr}  implies that 
one may consider $X_0$ with the induced metric instead of the length 
metric. We study the asymptotic cones of $X_0$ with the
induced  metric.

First there is a result on asymptotic cones of symmetric 
spaces and Euclidean buildings.

\begin{theorem}[\cite{KlL*}]\label{cX}
Any asymptotic cone of a product $X$ of symmetric spaces and
Euclidean buildings, $X$ of rank $r\geq 2$, is an  Euclidean 
building ${\mathbf K}$ of rank $r$ which is homogeneous and $\aleph_1$-thick. The apartments of ${\mathbf K}$ appear as limits of sequences of 
maximal flats 
in $X$. The same is true for Weyl chambers and walls, singular subspaces 
and Weyl polytopes of  ${\mathbf K}$. Consequently, $\partial_\infty 
{\mathbf 
K}$ and $\partial_\infty X$ have the same model spherical chamber and 
model 
Coxeter complex.
\end{theorem}

 In the sequel, in any asymptotic cone ${\mathbf K}$ of a
product $X$ of symmetric spaces and Euclidean buildings we
shall consider the labelling on
$\partial_\infty {\mathbf K}$ induced by a fixed labelling on 
$\partial_\infty X$. We denote the projection of $\partial_\infty 
{\mathbf K}$ 
on $\Delta_{mod} $ induced by this labelling by $P$ and the associated Coxeter complex by $\mathsf{S}$.

Concerning the asymptotic cone of a space $X_0$ obtained from a product of symmetric spaces and
Euclidean buildings 
by deleting disjoint open horoballs, we have the following result 
 
\begin{theorem}[\cite{Dr2*}, Propositions 3.10, 3.11]\label{cX0}
Let $X$ be a CAT(0) geodesic metric space and let ${\mathbf K}=X_\omega(x_n,d_n)$ be an 
asymptotic cone of $X$.

(1) If $(\rho_n)$ is a sequence of geodesic rays in $X$ with 
$\frac{d(x_n,\rho_n)}{d_n}$ bounded and $\rho =\lbrack \rho_n\rbrack $ is its limit ray in 
${\mathbf K}$, 
then $H(\rho 
)=\lbrack H(\rho_n)\rbrack $ and $Hb(\rho )=\lbrack Hb(\rho_n)\rbrack$.

(2) If $X_0=X \setminus \bigsqcup_{\rho \in {\mathcal R}}Hbo(\rho )$ and 
$\frac{d(x_n,X_0)}{d_n}$ is bounded then the limit set of $X_0$ (which is 
the same 
thing as the asymptotic cone of $X_0$ with the induced metric)
 is 
$$
{\mathbf K}_0={\mathbf K} \setminus \bigsqcup_{\rho_\omega \in
{\mathcal R}_\omega}Hbo(\rho_\omega )\;  ,\leqno(2.1)
$$ where 
${\mathcal R}_\omega $ is the set of rays $\rho_\omega =\lbrack \rho_n\rbrack 
,\; \rho_n \in 
{\mathcal 
R}$. 
\end{theorem}

We note that if $X$ is a product of symmetric spaces and
Euclidean buildings, the disjointness of 
$Hbo(\rho ),\; \rho \in 
{\mathcal R}$, implies by Lemma \ref{3} that $p(\lbrace \rho (\infty
)\mid \rho \in {\mathcal  R}\rbrace )$ 
reduces 
to one point, $\theta $, if card ${\mathcal R}\neq 2$. Then
$P(\lbrace \rho_\omega (\infty )\mid 
\rho_\omega \in {\mathcal R}_\omega \rbrace ) =\theta $.

We also need the following result. 

\begin{lemma}\label{platr}
Let $X$ be a product of symmetric spaces and
Euclidean buildings, $X$ of rank $r\geq 2$, and ${\mathbf{K}}=X_\omega(x_n,d_n)$ be an asymptotic cone of it. Let $F_\omega $ and $\rho_\omega $ be an apartment and a geodesic ray in ${\mathbf{K}}$, $F_\omega$ asymptotic to $\rho_\omega $. Let $\rho_\omega =[\rho_n ]$, where $\rho_n$ have the same slopes as $\rho_\omega $. Then
\begin{itemize}
\item[(a)] $F_\omega$ can be written as limit set $F_\omega =[F_n]$ with $F_n$ asymptotic to $\rho_n$ $\omega $-almost surely ;
\item[(b)] every geodesic segment $[x,y]$ in $F_\omega \setminus Hbo\, (\rho_\omega)$ may be written as limit set of segments $[x_n,y_n]\subset F_n \setminus Hbo\, (\rho_n)$.
\end{itemize} 
\end{lemma}

\noindent{\bf Proof.}\quad (a) Suppose first that $\rho_\omega $ is regular with slope $\theta $. Then $\rho_n$ are of slope $\theta $. Let $F_\omega =[F_n']$ and let $x_n'=$ proj$_{F_n'}\rho_n(0)$. We may replace in the argument each ray $\rho_n$ with the ray asymptotic to $\rho_n$ of origin $x_n'$. So in the sequel we may suppose that $\rho_n $ has its origin $x_n'$ in $F_n'$. Then $[\rho_n ]\subset [F_n']$, hence there exist rays $\rho_n'\subset F_n'$ of slopes $\theta $ and origin $x_n'$ such that if $x_n''$ is the first point of $\rho_n$ at distance $d_n$ from $\rho_n'$, $\lim_\omega \frac{d(x_n,x_n'')}{d_n}=+\infty $. If such a point $x_n''$ does not exist then $\rho_n \subset {\mathcal{N}}_{d_n}(F_n)$. This implies that $\rho_n \subset F_n$ and we are done. So we suppose that $x_n''$ always exists. For $n$ sufficiently large, $\rho_n(+\infty )$ becomes opposite to $\rho_n^{op}(+\infty )$ in $\partial_\infty X$, where $\rho_n^{op}$ is the ray opposite to $\rho_n'$ in $F_n'$. This happens because up to isometry we may suppose that $F_n$, $x_n'$ and $\rho_n'$ are fixed and then we can use the lower semicontinuity of the Tits metric with respect to the cone topology.

Let $F_n$ be the unique maximal flat containing $\rho_n(+\infty )$ and $\rho_n^{op}(+\infty )$ in its boundary. With an argument up to isometry similar to the previous one we may prove that $d(x_n',F_n)$ is uniformly bounded by a constant $M$. Let $x_n'''=$ proj$_{F_n}x_n''$. Let $W_n$ be the Weyl chamber of vertex $x_n'''$ containing $\rho_n^{op}(\infty )$ in the boundary and let $W_n'$ be the Weyl chamber asymptotic to it of vertex proj$_{F_n'}(x_n'')$. The Hausdorff distance $d_H(W_n, W_n')$ is at most $d_n+M$ hence $d_H([W_n], [W_n'])\leq 1$. On the other hand $\lim_\omega \frac{d(x_n,x_n'')}{d_n}=\lim_\omega \frac{d(x_n,x_n''')}{d_n}=+\infty $ so $[W_n]=[F_n]$ and $[W_n']=[F_n']$. It follows that $d_H([F_n], [F_n'])\leq 1$ which implies $[F_n]= [F_n']=F_\omega $.

The case when $\rho_\omega $ is singular can be reduced to the previous case by choosing $\rho_\omega^0 $ regular in the same Weyl chamber as $\rho_\omega $ and  asymptotic to $F_\omega $ and repeating the previous argument.

\smallskip

(b) The set $F_\omega \setminus Hbo\, (\rho_\omega)$ as well as $\omega $-almost all sets $F_n \setminus Hbo\, (\rho_n)$ are half-flats and the conclusion follows easily. \hspace*{\fill } $\diamondsuit $

\section{Horospheres in Euclidean buildings : intersections with apartments and nondistorsion}\label{horo}

\subsection{Intersection of a horoball with an apartment : global properties}\label{ha}

We first introduce a notation. For every $k$-flat 
$\Phi$ (not necessarily singular) and every geodesic ray $\rho $ such that 
$\inf_{x\in
\Phi}f_\rho(x) > -\infty $ we denote 
$$
Min_\Phi(\rho ):=\lbrace y\in 
\Phi \mid 
f_\rho (y)=\inf_{x\in
\Phi}f_\rho(x)\rbrace \; .
$$

We recall that $f_\rho $ denotes the Busemann function associated to $\rho $.

We use a similar notation for a half-apartment $D$ instead of the flat $\Phi $.

We describe some global features of intersections between 
horoballs/horospheres and apartments.

\begin{proposition}[\cite{Dr1*}, Proposition 3.1, Lemma 3.2, Lemma 
3.8]\label{form}
 Let ${\mathbf K}$ be an Euclidean building, $F$ an apartment
in it and
$\rho
\subset {\mathbf K}$ a geodesic ray of slope $\theta$. Let 
$\{\alpha_1,\alpha_2,\dots \alpha_k \}$ be the points in $F(\infty )$ opposite 
to $\rho(\infty )$, $k\leq q_0$, where $q_0$ is the number of chambers of the 
Coxeter complex associated to ${\mathbf K}$.
\begin{itemize}

\item[(a)] The flat $F$ can be written as $F=\bigcup_{i=1}^m[F\cap F_i]$, where

$\bullet $  each 
$F_i$ is an apartment with $F_i(\infty )$ containing 
$\rho (\infty )$ and $\alpha_i$ ;

$\bullet $  each $F(\infty )\cap F_i(\infty )$ contains a unique point 
$\alpha_{j_i}$ opposite to $\rho (\infty )$.

It is possible that there exist $i_1\neq i_2$ with 
$\alpha_{j_{i_1}}=\alpha_{j_{i_2}}$.   

\item[(b)] The intersection $Hb(\rho )\cap F$, if nonempty, is a
convex polytope in $F$. The hypersurface $H(\rho )\cap F$ has at most $k$ faces.

\item[(c)] If $\inf_{x\in F}f_\rho(x)>-\infty$ then the set
$Min_F(\rho )$ is a Weyl polytope with at most $m$ faces, $m=m({\mathsf{S}})$. All the slopes of
$Min_F(\rho )$ are contained in
$Ort\; (\theta )$.

If $\Phi$ is a $k$-flat (not necessarily singular) contained in $F$, 
$Min_\Phi(\rho )$ is the intersection of a Weyl
polytope with 
$\Phi $.
\end{itemize}
\end{proposition}

We note that (b) follows from (a) since $Hb(\rho )\cap 
F$ is a convex set, $Hb(\rho )\cap F 
=\cup_{i=1}^m[Hb(\rho )\cap F\cap F_i] $, $Hb(\rho )\cap F_i$ is a 
half-apartment and $F\cap F_i$ is a Weyl polytope. In (c) the fact that the number of faces of $Min_F(\rho )$ is uniformly bounded is a consequence of the fact that all faces are supported by intersections of the affine span with singular hyperplanes, and that, by the convexity of $Min_F(\rho )$, there can be at most two parallel faces.

\begin{cor}\label{prf}
 Let ${\mathbf K}$ be an Euclidean building, $\rho
\subset {\mathbf K}$ a geodesic ray of slope $\theta$ and $F$ an apartment
 which intersects $Hb_t(\rho)$. For every $s>t$ the projection in $F$ of 
$H_s(\rho )\cap F$ onto $H_t(\rho )\cap F$ is a contraction. Moreover there 
exists a constant $a>1$ depending only on $\theta $ and on the Coxeter complex $\mathsf{S}$ of 
${\mathbf K}$ such that the distance from a point of $H_s(\rho )\cap F$ to its 
projection onto $H_t(\rho )\cap F$ is at most $a(s-t)$.  
\end{cor}

\noindent{\bf Proof.}\quad The first statement is a straightforward consequence 
of Proposition \ref{form}, (b), and of the fact that $Hb_t(\rho )\subset 
Hb_s(\rho )$.

Let $x$ be a point in $H_s(\rho )\cap F$ and $y$ its projection on $H_t(\rho 
)\cap F$. For every point $z\in Hb_t(\rho )\cap F$, $\angle_y(\overline{yx}, 
\overline{yz})\geq \frac{\pi}{2}$. Let $M$ be the unique wall of vertex $y$ 
containing the segment $[y,x]$ in its interior. Then $M$ cannot have a segment 
in common with $H_t(\rho )\cap F$, otherwise the diameter of the wall $M_y$ in 
the spherical building $\Sigma_y{\mathbf K}$ would be $>\pi /2$. It follows that 
$\angle_y(M_y,\overline{\rho_y})>\frac{\pi}{2}$ and therefore, by Remark 
\ref{discret}, that 
$\angle_y(M_y,\overline{\rho_y})\geq \frac{\pi}{2}+\delta_0$, with $\delta_0$ 
depending only on $\theta$ and on $\mathsf{S}$. In particular 
$\angle_y(\overline{yx},\overline{\rho_y})\geq \frac{\pi}{2}+\delta_0$, which 
implies that $s=f_\rho(x)\geq f_\rho(y)+d(y,x)\sin \delta_0=t+d(y,x)\sin 
\delta_0$. The conclusion follows with $a=\frac{1}{\sin 
\delta_0}$.\hspace*{\fill }$\diamondsuit $

\medskip

A consequence of the previous Corollary is that for an apartment $F$ such that\\ $\inf_{x\in F}f_\rho (x) 
>-\infty $ the set $Min_F(\rho )$ of points minimizing $f_\rho $ on $F$ plays an 
essential part in the shape of $F\cap Hb(\rho )$.

\medskip

\begin{cor}\label{vec}
Let ${\mathbf K}$ be an Euclidean building, $F$ an apartment in it and 
$\rho $ a
ray of slope $\theta $. If $\inf_{x\in F }f_\rho (x)=-m>-\infty $ then there 
exists a constant $a>1$ depending only on the Coxeter complex of ${\mathbf K}$ and on $\theta $ such that 
$$
{\mathcal{N}}_m(Min_F(\rho ))\subset Hb(\rho )\cap F \subset {\mathcal{N}}_{am}(Min_F(\rho ))\; .
$$
\end{cor}

\noindent{\bf Proof.}\quad The first inclusion follows from the fact that the Busemann function $f_\rho $ is a contraction. The second inclusion follows from the previous corollary applied to $H_{-m}(\rho )\cap F$ and each $H_{-s}(\rho )\cap F $ with $0\leq s < m$.\hspace*{\fill }$\diamondsuit $

\medskip

Corollary \ref{vec} shows that the shape of $Min_F(\rho )$ is important for 
the shape of $F\cap Hb(\rho )$. Therefore in the sequel we formulate results 
about the way $Min_F(\rho )$ changes when we change the apartment $F$ with a 
ramification of it.

\begin{lemma}[\cite{Dr1*}, Lemma 3.12]\label{hipp}
Let ${\mathbf K}$ be a $4$-thick
Euclidean building of rank at least $2$, let $\rho $ be a geodesic ray in it, 
not parallel to any factor, let $F$ be an apartment, $D$ a half-apartment and  
$\Phi $ the affine span of $Min_D(\rho )$ in $F$. Let $D_0$ be the 
half-apartment in $F$ opposite to $D$, and $D_1,\; D_2$ two other 
half-apartments of boundary $H=\partial D$ and interiors mutually disjoint and 
disjoint of $F$.

If the singular hyperplane $H$ neither contains $\Phi $ nor is orthogonal to it 
then $\inf_{D_i}f_\rho \geq \inf_{D}f_\rho , \; \forall i\in \{0,1,2 \}$, and 
$Min_{D_i\cup D}(\rho )=Min_{D}(\rho )$ for at least two values of $i$ in 
$\{0,1,2 \}$.
\end{lemma}

\medskip

Lemma \ref{hipp} implies that, by eventually changing a half-apartment, we can 
always make sure that $Min_F(\rho )$ stays in a single half-apartment, and by 
applying the Lemma twice we can even make sure that $Min_F(\rho )$ stays in a 
strip determined by two parallel hyperplanes. The only thing needed to apply 
the Lemma is a singular hyperplane which neither contains nor is orthogonal to 
the affine span of $Min_F(\rho )$. The following results are about the existence 
of such singular hyperplanes.  

\begin{lemma}\label{skew}
Let $X$ be a product of symmetric spaces and Euclidean buildings, $X$ of rank 
$r\geq 2$, let $\rho $ be a geodesic ray which is not parallel to any factor of 
$X$, let $F$ be a maximal flat containing $\rho $ and let $\Phi \subset F$ be a singular 
flat orthogonal to $\rho $. Every Weyl chamber of $F$ has an adjacent singular 
hyperplane which is neither parallel to $\Phi$ (or containing $\Phi$) nor 
orthogonal to $\Phi$. 
\end{lemma}

\noindent{\bf{Proof.}}\quad Since all Weyl chambers are translations of a finite 
set of Weyl chambers with common vertex, we may restrict the problem to Weyl 
chambers having the vertex on $\Phi $. Let $W$ be such a Weyl chamber. First we 
work in the hypothesis that $\Phi $ supports $W$. In this case the result 
follows from Lemma 3.11, \cite{Dr1*}.

Suppose $\Phi $ does not support $W$. If the dimension of $\Phi $ is $k$ then 
the number of singular hyperplanes supporting $W$ and containing $\Phi $ must be 
strictly less than $r-k$. So there are at least $k+1$ hyperplanes supporting $W$ 
not containing $\Phi $. They cannot be all orthogonal to $\Phi $ because if we 
choose vectors orthogonal to each of these hyperplanes we obtain a family of at 
least $k+1$ linearly independent vectors, so the entire family cannot be 
contained in $\Phi $.\hspace*{\fill } $\diamondsuit $

\begin{cor}\label{cskew}
Let $X$ and $\rho $ be as in the previous Lemma. If $F$ is a flat not containing 
$\rho $, ${\mathcal P} $ is a Weyl polytope in $F$ on which $f_\rho$ is constant 
and $\Phi $ is its affine span, the conclusion of the Lemma \ref{skew} still 
holds.
\end{cor}

 \noindent{\bf{Proof.}}\quad By Proposition \ref{form}, (a), an open subset 
$\Omega $ of the Weyl polytope ${\mathcal P} $ shall be entirely included into 
an intersection $F\cap F_i$, where $F_i$  contains a ray $\rho_i$ asymptotic to 
$\rho $ and $F(\infty )\cap F_i(\infty )$ contains a unique point $\alpha $ 
opposite to $\rho(\infty )$. Let $\rho^{op}$ be a ray with origin on $\Omega $ 
and $\rho^{op}(\infty ) =\alpha$. Then $\rho^{op}$ is contained in $F\cap F_i$ 
and in $F_i$ it is opposite to a ray parallel to $\rho_i$. Hence it is 
orthogonal to $\Omega  $ and to $\Phi$ in $F$, on one hand, and on the other 
hand, it is not parallel to a factor in $X$. By applying Lemma \ref{skew} with 
$\rho^{op}$ instead of $\rho $ we obtain the same conclusion.\hspace*{\fill 
}$\diamondsuit $

\medskip

\subsection{Intersection of a horoball with an apartment : local properties}\label{hac}

We look at an intersection $F\cap H(\rho )$ in the neighborhood of one of its points $x$. 
The first and the second part of the following proposition give a meaning to the 
directions orthogonal to the codimension one faces of $F\cap H(\rho )$ through $x$. The third part 
describes a situation in which $F\cap H(\rho )$ has two codimension one faces 
through $x$ symmetric with respect to a singular hyperplane through $x$. This 
result is essential in the argument of ``breaking the faces'' described in Lemma 
\ref{fetze2}.   

\begin{proposition}\label{formlocal}
Let ${\mathbf K}$ be an Euclidean building, $\rho
\subset {\mathbf K}$ a geodesic ray of slope $\theta $ and $F$ an apartment
 which intersects $Hb(\rho)$. Let $x\in F\cap H(\rho )$. 

\begin{itemize}
\item[(a)] There exists $\delta >0$ such that $F\cap B(x,\delta 
)=\bigcup_{i=1}^s[F\cap F_i\cap B(x,\delta )]$, where

$\bullet $  $F_i$ is an apartment and $F_i\cap B(x,\delta )$ contains $\rho_x 
\cap B(x,\delta )$, and

$\bullet $  each set $F\cap F_i\cap B(x,\delta )$ contains only one segment of 
length $\delta $, $[x,a_i)$, of direction $\overline{xa_i}$ opposite  to 
$\overline{\rho_x}$.

It is possible that there exist $i_1\neq i_2$ with $a_{i_1}=a_{i_2}$. 
\item[(b)] Every codimension one face of $F\cap H(\rho )$ through $x$ is orthogonal to 
one of the segments $[x,a_i)$. A segment $[x,a_i)$ is orthogonal to a face of 
$F\cap H(\rho )$ through $x$ if and only if for some component $F\cap F_i \cap 
B(x,\delta )$ containing it either $\angle_x(\overline{\rho_x}, (F\cap 
F_i)_x)<\frac{\pi}{2}$ or $\angle_x(\overline{\rho_x}, (F\cap 
F_i)_x)=\frac{\pi}{2}$ and $(F\cap F_i)_x$ contains a panel orthogonal to 
$\overline{\rho_x}$.

\item[(c)] Let $W^1$ and $W^2$ be two adjacent chambers in $F$, $W^1\cap W^2 = M$. Suppose that $W_x^1$ and $W_x^2$ are at the same combinatorial distance from 
$\overline{\rho_x}$ and that $\angle_x(W^t_x,\overline{\rho_x})<\frac{\pi}{2},\; 
t=1,2$. Suppose no singular hyperplane through $\overline{\rho_x}$ contains $M_x$. Then there exist $i_1\neq i_2$ such that 
$W^t\cap B(x,\delta )\subset F\cap F_{i_t}\cap B(x,\delta ),\; t=1,2$, and the 
segments $[x,a_{i_1}]$ and $[x,a_{i_2}]$ are symmetric with respect to $Span\; M $.
\end{itemize}
\end{proposition}

\noindent{\bf Proof.}\quad (a)  We choose an auxiliary ray $\rho^0$. If $\rho $ 
is regular, then $\rho^0 =\rho_x$ while if $\rho $ is singular we choose a Weyl 
chamber $W^0$ of vertex $x$ containing $\rho_x$, and $\rho_0$ a regular ray in 
$W^0$ of origin $x$. Let $\{\beta_1,\beta_2,\dots \beta_s \}$ be the points in 
$F_x$ opposite to $\overline{\rho^0}$ and let $A_i$ be the unique apartment 
containing $\beta_i$ and $\overline{\rho^0}$. By \cite[Lemma 3.10.2]{KlL*} every 
point in $F_x$ lies in some $A_i$, so $F_x=\bigcup_{i=1}^s [F_x\cap A_i ]$. By 
\cite[Lemma 4.2.3]{KlL*}, $A_i=(F_i)_x$, where $F_i$ is an apartment in the 
Euclidean building containing $x$. Lemma 4.1.2, (1), and Sublemma 4.4.1 of 
\cite{KlL*} allow to conclude that (a) is true.

The situation when there exist $i_1\neq i_2$ with $a_{i_1}=a_{i_2}$ may appear when $\rho $ is singular. 

\smallskip

(b) First we notice that if $Ort\; (\theta)\not\subset \partial \Delta_{mod} $ 
then any open set in it contains a regular point, while if $Ort\; (\theta) 
\subset \partial \Delta_{mod} $ any open set in it contains a point in the 
interior of a panel. For every codimension one face $\mathfrak f$ of $H(\rho )\cap F$ 
through $x$ we consider a point $\eta \in {\mathfrak f}_x$, regular if we are in 
the first case or in the interior of a panel if we are in the second case. The 
point $\eta $ is contained in $F_x\cap A_i$ for some $i\in \{1,2,\dots s \}$.  
It follows that an open subset of $\mathfrak f$ is contained in $F\cap F_i\cap 
B(x,\delta )$, so it is orthogonal to $[x,a_i)$.

We now prove the equivalence stated in (b). We only prove the direct 
implication, since the reciprocal is an easy consequence of previously used 
results. Suppose that $[x,a_i)$ is orthogonal to a codimension one face 
$\mathfrak f$ of $H(\rho )\cap F$ through $x$. We choose a point $\eta $ as previously 
in ${\mathfrak f}_x$. In $\Sigma_x{\mathbf K}$ we consider the geodesic between 
$\overline{xa_i}$ and $\eta $, of length $\frac{\pi }{2}$ and which prolongates 
till it hits $\overline{\rho_x}$. This geodesic contains a regular point $\eta^0 
\neq \eta$ in a chamber $W_x\subset F_x$ intersecting ${\mathfrak f}_x$ either 
in a relatively open set of regular points or in a panel. On the other hand 
$W_x$ is contained in some $(F\cap F_j)_x$ and the geodesic from $\eta $ to 
$\eta^0$ prolongates to the opposite $\overline{xa_j}$. But since in $F_x$ the 
geodesic from $\eta $ to $\eta^0$ prolongates in a unique way, it follows that 
$j=i$. Since $(F\cap F_j)_x$ has in common with ${\mathfrak f}_x$ either a 
relatively open set of regular points or a panel, we may conclude.

\smallskip

(c) Statement (a) implies that there exists $i_t$ such that $W^t\cap B(x,\delta 
)\subset F\cap F_{i_t}\cap B(x,\delta ),\; t=1,2$. Since $W^1_x$ and $W^2_x$ are 
at the same combinatorial distance from $\overline{\rho_x}$ and no singular 
hyperplane through $\overline{\rho_x}$ contains the common panel $M_x$, none of  
$W^1_x$ or $W^2_x$ can be the projection of $\overline{\rho_x}$ on $M_x$ and 
both are separated from $\overline{\rho_x}$ by $M_x$. Also, $W^1_x\neq W^2_x$ implies $i_1\neq i_2$. Moreover, 
$W^1_x$ and $W^2_x$ are at the same combinatorial distance from $(W^0)_x$, which 
implies that retr$_{(F_{i_1})_x, (W^0)_x }(W^2_x)=W^1_x$.

The hypothesis that $\angle_x(W^t_x,\overline{\rho_x})<\frac{\pi}{2}$ implies 
that there exist two points $\xi_t $ in the interiors of $W^t_x$ with 
$\angle_x(\xi_t ,\overline{\rho_x})<\frac{\pi}{2},\; t=1,2,$ and $p_x(\xi_1 
)=p_x(\xi_2)$. In $(F_{i_t})_x$ the geodesic $\gamma_t$ joining 
$\overline{\rho_x}$ to $\xi_t $ prolongates to a geodesic joining 
$\overline{\rho_x}$ to $\overline{xa_{i_t}}$. The previous considerations imply 
that retr$_{(F_{i_1})_x,(W^0)_x }(\xi_2)=\xi_1$, so that 
retr$_{(F_{i_1})_x,(W^0)_x }(\gamma_2)=\gamma_1$ and in particular that 
$\gamma_1$ and $\gamma_2$ intersect $M_x$ in the same point $\zeta $. In $F_x$ 
the geodesics from $\zeta $ to $\xi_1$ and from $\zeta $ to $\xi_2$ are 
symmetric with respect to $Span\; M_x$, so the same is true 
for their prolongations to $\overline{xa_{i_1}}$ and to $\overline{xa_{i_2}}$, 
respectively. In particular $a_{i_1}$ and $a_{i_2}$ are symmetric with respect 
to $Span\; M$, which ends the proof.\hspace*{\fill }$\diamondsuit $

\medskip

\begin{cor}\label{sim}
Let ${\mathbf K}$ be a 3-thick Euclidean building, $x$ a point
in it and $\rho
\subset {\mathbf K}$ a geodesic ray of slope $\theta$. Let $W$ be a Weyl chamber of vertex $x$ such that $W_x$ does not contain $\overline{\rho_x}$ and $\angle_x(\overline{\rho_x}, W_x)<\frac{\pi }{2}$. Let $M$ be a panel such that $M_x$ separates $W_x$ and $\overline{\rho_x}$. Then there exists a Weyl chamber $\widehat{W}$ adjacent to $W$, with $\widehat{W}\cap W=M$, such that in any apartment $F$ containing both $W$ and $\widehat{W}$ the point $x$ is contained in two faces of $F\cap H(\rho )$. 
\end{cor}

\noindent{\bf Proof.}\quad Remark \ref{discret}, implies that $\angle_x(\overline{\rho_x}, W_x)\leq \frac{\pi }{2}-\delta_0$, where $\delta_0 $ depends 
only on $\theta $. The chamber $W$ has a panel $M$ such that $M_x$ separates $W_x$ and $\overline{\rho_x}$. There exist at least two chambers adjacent to $W$ containing $M$ and at least one of these two chambers, $\widehat{W}$, has the property that $W_x$ and $\widehat{W}_x$ are at the same combinatorial distance from $\overline{\rho_x}$. Let $F$ be an apartment containing $W$ and $\widehat{W}$. By Proposition \ref{formlocal}, (c), $x$ is contained into two faces of $F\cap H(\rho )$, orthogonal to two segments $[x,a_1)$ and $[x,a_2)$ symmetric with respect to $Span\; M$. We denote $\widehat{H}:= Span\; M$. By the properties of $M$, $\widehat{H}$ does not contain 
$[x,a_1)$. Also, $\angle_x(\overline{\rho_x} , 
(\widehat{H})_x)\leq \angle_x(\overline{\rho_x} , M_x)\leq 
\angle_x(\overline{\rho_x} , W_x)\leq \frac{\pi }{2}-\delta_0$. Therefore 
$\widehat{H}$ is not orthogonal to $[x,a_1)$ neither. It follows that $[x,a_1)$ and $[x,a_2)$ are not on the same line, which implies that the two corresponding faces are distinct.\hspace*{\fill }$\diamondsuit $
 
\medskip

\subsection{Nondistorsion of horospheres}\label{ndh}

In the proof of our main theorem on filling (Theorem \ref{fillgen}) we shall need the 
following two Lemmata. As a byproduct we also obtain new proofs of 
the Theorems 1.1 and 1.2 in \cite{Dr1*} on the nondistorsion of horospheres in 
Euclidean buildings and symmetric spaces, considerably shorter than the ones 
given in the previously cited paper.

Both Lemmata deal with the possibility of joining two points 
$x,y$ in the intersection of an apartment with a horosphere $F\cap H(\rho )$ by 
a polygonal line of length comparable to $d(x,y).$  Lemma \ref{fetze} states 
that this can easily be done if $x$ and $y$ are contained into non-parallel 
codimension one faces of $F\cap H(\rho )$. Lemma \ref{fetze2} deals with the 
case when $x$ and $y$ are contained in the interiors of parallel faces. The main 
idea in the argument is that, by choosing a singular hyperplane skew to the two 
affine spans and by changing a half-apartment bounded by this hyperplane one is 
able to ``break'' one of the faces into two different faces. The details we give 
in the statement of the Lemma on the different possibilities of ``breaking'' one 
of the two faces will be necessary further on.

\begin{lemma}[non-parallel faces]\label{fetze}
Let ${\mathbf K} $ be an Euclidean building of rank at least $2$ and let $\rho$ 
be a geodesic ray in it of slope $\theta $, not parallel to a rank one factor. 
There exists a constant $C= C(\mathsf{S},\theta )$ such that for every apartment $F$ intersecting $Hb(\rho ) 
$, every two points $x,y\in F\cap H(\rho )$, contained into two distinct 
non-parallel codimension one faces of $F\cap H(\rho )$, can be joined by a 
polygonal line in $F\cap H(\rho )$ of length at most $Cd(x,y)$. Moreover the 
polygonal line may be chosen in a half-plane in $F$ having as boundary the line 
$xy$.

In particular the previous statement applies for every pair of points in $F\cap 
H(\rho )$ one of which is contained into two different codimension one faces of 
$F\cap H(\rho )$.
\end{lemma}

\noindent{\bf{Proof.}}\quad By the hypothesis, $x$ and $y$ are contained into 
two codimension one faces of the convex polytope $F\cap Hb(\rho )$ whose affine 
spans, $H$ and $H'$ respectively, are not parallel. The set of slopes of both 
$H$ and $H'$ is $Ort\; (\theta)$. There is a minimal possible dihedral angle 
between two distinct hyperplanes of set of slopes $Ort\; (\theta)$, which we 
denote by $\varsigma $. So $\angle (H, H')\geq \varsigma $. As $F\cap Hb(\rho )$ 
is a convex polytope, it is entirely contained in a skew quadrant determined by 
$H$ and $H'$, which we denote by ${\mathcal Q}$. On the other hand we have $x\in 
H$ and $y\in H'$. Since $\angle (H, H')\geq \varsigma $, there exists a constant 
$C$ depending only on $\varsigma $ such that $x$ and $y$ may be joined by a 
polygonal line contained in $\partial {\mathcal Q}$ of length at most $Cd(x,y)$. 
The polygonal line of minimal length will be of the form $[x,z]\cup [z,y]$ with 
$z\in H\cap H'$. The plane determined by the points $x,y,z$ in $F$ intersects 
$Hb(\rho )$ in a convex polygon entirely contained in the sector of vertex $z$ 
and sides the rays through $x$ and $y$. In follows that the 
polygonal line joining $x$ and $y$ inside the triangle $xyz$ satisfies the 
conclusion of the lemma.\hspace*{\fill }$\diamondsuit $

\medskip

We note that the hypothesis of $\rho $ not being parallel to a rank one factor 
is necessary for the existence of two non-parallel faces.

\medskip

{\it Notation :} Henceforth for a pair of points $x,y$, in the intersection of 
an apartment with a horosphere $F\cap H(\rho )$, $x,y$ in nonparallel faces, we 
shall denote by ${\mathbf L}_{xy}$ a polygonal line joining 
them in $F\cap H(\rho )$, of length at most $Cd(x,y)$, contained in a half-plane of boundary $xy$, constructed as previously.

\smallskip

\begin{lemma}[breaking parallel faces]\label{fetze2}
Let ${\mathbf K} $ , $\rho$ and $\theta $ be as in the previous Lemma with the 
additional hypothesis that ${\mathbf K}$ is $3$-thick and that $\rho $ is not parallel to any factor of  ${\mathbf 
K} $. Let $F$ be an apartment intersecting $Hb(\rho ) $, and $x$ and $y$ two 
points in $F\cap H(\rho )$. Suppose $x$ 
and $y$ are contained in the interiors of two distinct codimension one faces of 
$F\cap H(\rho )$ with parallel affine spans $H$ and $H'$, respectively. Let 
$M_0$ be the unique wall of vertex $x$ containing $[x,y]$ in its interior.

(1) There exists a constant $C= C(\mathsf{S},\theta )$ and a ramification $F'$ of $F$ containing $M_0$ 
such that $x$ and $y$ may be joined in $F'\cap H(\rho )$ by a polygonal line of 
length at most $Cd(x,y)$.

(2) The choice of the ramification $F'$ can be made as follows :

\hspace{1cm}(a) Suppose $(M_0)_x$ does not contain $\overline{\rho_x}$. Let $W$ 
be a Weyl chamber in $Star\; (M_0)\cap F$ such that $W_x$ is not the projection of $\overline{\rho_x}$ on $(M_0)_x$. Then there exists a ramification $F'$ of $F$ containing 
$W$ such that $x$ is contained into 
two faces of $F'\cap H(\rho )$.

\hspace{1cm}(b) Suppose $(M_0)_x$ contains $\overline{\rho_x}$.

\hspace{1.5cm}(b$_1$) If the connected component of $Star\; (M_0)\cap F \cap 
H(\rho )$ containing $y$ has at least two faces then $F'=F$.

\hspace{1.5cm}(b$_2$) Suppose the connected component of $Star\; (M_0)\cap F 
\cap H(\rho )$ containing $y$ has one face. Let $W$ be a Weyl chamber in $Star\; 
(M_0)\cap F$. For every hyperplane $\widehat{H}$ supporting $W$ which is neither 
orthogonal nor coincident with $H$, there exists a ramification $F'$ of $F$ 
containing $W$ such that $\partial (F'\cap F)=\widehat{H}$ and all the points in 
$H(\rho )\cap \widehat{H}\cap W$ are in two different faces of $F'\cap H(\rho 
)$.      
\end{lemma}

\begin{figure}[!ht]
\centering
\includegraphics{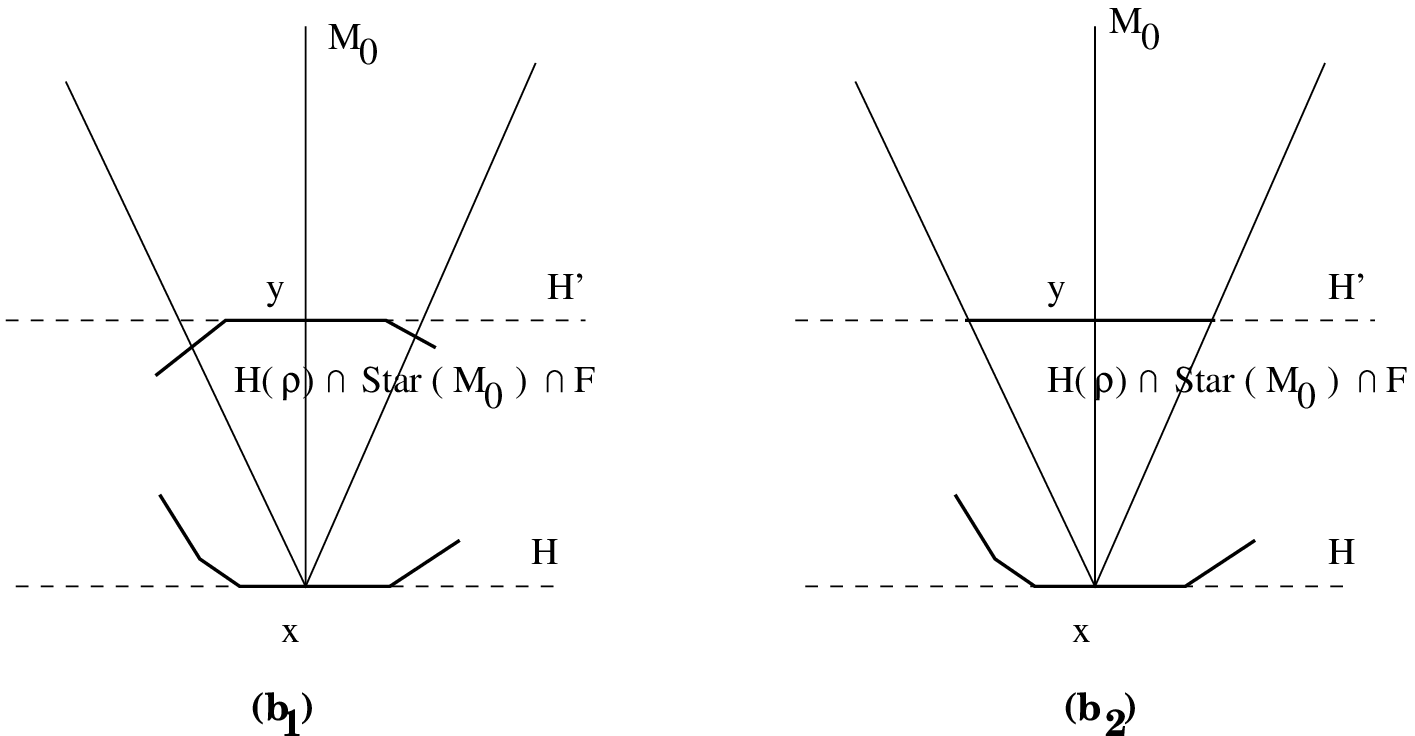}
\caption{}
\label{fig1} 
\end{figure}

\noindent{\bf Proof.}\quad In Figure \ref{fig1} the cases (b$_1$) and (b$_2$) are 
represented.

We prove that the choices of the ramification $F'$ proposed in (2) are always possible and that they satisfy (1). Before proceeding, we make some general remarks.

By hypothesis $(x,y)\subset F\cap Hbo(\rho )$, hence $f_{\rho }$ decreases strictly on $[x,y]$ starting from $x$. It follows that 
$\angle_x(\overline{xy},\overline{\rho_x})<\frac{\pi }{2}$ which, by Remark 
\ref{discret}, implies that $\angle_x((M_0)_x,\overline{\rho_x})\leq \frac{\pi 
}{2}-\delta_0$, $\delta_0 =\delta_0 (\theta , \mathsf{S})$. Therefore, for any Weyl 
chamber $W$ containing $M_0$, $\angle_x(W_x,\overline{\rho_x})\leq \frac{\pi 
}{2}-\delta_0$.

By Proposition \ref{formlocal}, (a), there exists $\delta $ such that $F\cap 
B(x,\delta )=\bigcup_{i=1}^s[F\cap F_i\cap B(x,\delta )]$, where $F_i\cap 
B(x,\delta )$ contains $\rho_x \cap B(x,\delta )$, and each set $F\cap F_i\cap 
B(x,\delta )$ contains only one segment of length $\delta $, $[x,a_i)$, of 
direction $\overline{xa_i}$ opposite  to $\overline{\rho_x}$. The fact that $x$ is in the interior of one face of $H(\rho )\cap F$ and Proposition 
\ref{formlocal}, (b), imply that there exist at most two opposite segments, 
$[x,a_1)$ and $[x,a_2)$, with the property of being contained in components $F 
\cap F_{i_t} \cap B(x,\delta ),\; t=1,2,$ with $\angle_x(\overline{\rho_x}, 
(F\cap F_{i_t})_x)<\frac{\pi }{2}$ or with $\angle_x(\overline{\rho_x}, (F\cap 
F_{i_t})_x)=\frac{\pi }{2}$ and $(F\cap F_{i_t})_x$ containing a panel 
orthogonal to $\overline{\rho_x}$. If $W$ is a Weyl chamber containing $M_0$, 
$W\cap B(x,\delta )$ must be included in one of these components. Suppose it 
is in $F\cap F_{i_1} \cap B(x,\delta )$, containing the segment $[x,a_1)$.

A geodesic from $\overline{\rho_x }$ to $\overline{xy}$ prolongates in a unique 
way in $(F\cap F_{i_1})_x$ and it ends in $\overline{xa_1 }$. Therefore 
$\angle_x (\overline{xy}, \overline{xa_1 } )\geq \frac{\pi }{2}+\delta_0$, on one 
hand, and on the other hand the geodesic between $\overline{xy}$ and $ 
\overline{xa_1 }$ in $(F\cap F_{i_1})_x$  contains a point at distance 
$\frac{\pi }{2}$ from $\overline{\rho_x }$. It follows that the segments $(x,y]$ 
and $(x,a_1]$ are separated by the hyperplane $H$. We conclude that the segment 
$(x,a_1]$ is precisely the segment orthogonal to $H$ and on the other side of 
$H$ than $[x,y]$.

\medskip

(2) (a) Let $W$ be a Weyl chamber in $Star \; (M_0)$ such that $W_x$ is not the projection of $\overline{\rho_x}$ on $(M_0)_x$. In particular $W_x$ does not contain $\overline{\rho_x}$ and it has a panel $M_x$ separating $W_x$ and $\overline{\rho_x}$. This and the fact that 
$\angle_x(W_x,\overline{\rho_x})\leq \frac{\pi }{2}-\delta_0$ imply, by Corollary \ref{sim} and Lemma \ref{apgcw}, that there exists a ramification $F'$ of $F$ with the desired properties.

\medskip

(2) (b) The hypothesis that $\rho $ is not parallel to any factor and Lemma 
\ref{fnort} imply that $\rho_x$ is not perpendicular near $x$ to any wall of any Weyl 
chamber in $Star\; (M_0)\cap F$. In particular the intersection of $H$ with 
$Star\; (M_0)\cap F$ reduces to the point $x$. This implies that the intersection 
of any Weyl chamber in $Star\; (M_0)\cap F$ with $H'$ is a simplex of diameter 
$\leq {\bf c}\cdot d(x,y)$, ${\bf c} ={\bf c} (\theta , \mathsf{S})$. It follows that 
$Star\; (M_0)\cap F \cap H'$ is a convex polytope in $H'$ of diameter $\leq {\bf 
c}'\cdot d(x,y)$, ${\bf c}' ={\bf c}' (\theta , \mathsf{S})$. From this and the fact that 
the connected component of $Star\; (M_0)\cap F \cap H(\rho )$ containing $y$ is 
contained in the truncated polytopic cone determined by $Star\; (M_0)\cap F$ and 
$H'$ we can deduce that the diameter of the connected component of $Star\; 
(M_0)\cap F \cap H(\rho )$ containing $y$, considered with its length metric, is 
at most ${\bf c}''\cdot d(x,y)$, ${\bf c}'' ={\bf c}'' (\theta , \mathsf{S})$. 

\medskip

(b$_1$) By hypothesis there exists $y_1$ in the connected component of $Star\; 
(M_0)\cap F \cap H(\rho )$ containing $y$, which is contained into two faces. The 
previous considerations imply that $y_1$ may be joined to $y$ by a polygonal 
line in $F\cap H(\rho )$ of length at most ${\bf c}''\cdot d(x,y)$. Lemma \ref{fetze} implies that 
$y_1$ may be joined to $x$ in $F\cap H(\rho )$ by a polygonal line in $F\cap H(\rho )$ of length at 
most $C({\bf c}''+1)\cdot d(x,y)$.

\medskip

(b$_2$) In this case the connected component of $Star\; (M_0)\cap F \cap H(\rho 
)$ containing $y$ is simply $Star\; (M_0)\cap F \cap H'$. Let $W$ and 
$\widehat{H}$ be as in the statement. The existence of $\widehat{H}$ is 
guaranteed by Corollary \ref{cskew}. We denote $M:= W\cap\widehat{H}$, $\widehat{D}$ the 
half-apartment bounded by $\widehat{H} $ containing $W$ and $\widehat{D}'$ the 
opposite half-apartment in $F$. We have that $H(\rho )\cap \widehat{H} \cap W = 
H'\cap M $. It is enough to prove the statement in (b$_2$) 
for some point $y_1\in Int\; M\cap H(\rho )=Int\; M\cap H'$. Suppose $y_1$ is also in the 
interior of the face of span $H'$, otherwise we are in the case (b$_1$) and we 
may take $F'=F$.

By hypothesis $\overline{\rho_x} \in (M_0)_x$ so $\rho_x \cap 
B(x,\varepsilon )\subset M_0\cap B(x,\varepsilon )$ for a small $\varepsilon 
>0$. Let $d$ be the line in $F$ containing the segment $\rho_x \cap 
B(x,\varepsilon )$. Then $d$ is orthogonal to $H$ and $H'$ and $d\cap M_0$ is a 
ray. The fact that $d\cap M_0$, $[x,y]$ and $M$ are contained in the same Weyl chamber 
implies that the skew quadrant determined by $\widehat{H}$ and $H'$ containing 
$[x,y]$ has a dihedral angle $\alpha <\frac{\pi }{2}$ (Figure 2). The same 
follows for the opposite quadrant, contained in $\widehat{D}'$. Let ${\mathcal 
D}H'$ and ${\mathcal D}\widehat{H}$ be the half-hyperplanes bounding this 
quadrant. For a small $\varepsilon_1$, $H(\rho )\cap \widehat{D}' \cap B(y_1, 
\varepsilon_1 )={\mathcal D}H'\cap B(y_1, \varepsilon_1 )$ and it makes with 
${\mathcal D}\widehat{H}$ the dihedral angle $\alpha $.

Let $\widehat{D}''$ be another half-apartment bounded by $\widehat{H}$ whose 
interior is disjoint from $F$. We now consider the apartment $F'=\widehat{D}\cup 
\widehat{D}''$. Suppose by absurd that $y_1$ is in the interior of a face also 
in $F'$. Let $H''$ be the affine span in $F'$ of $H'\cap \widehat{D}$ and let 
${\mathcal D}H''$ be its half which is opposite to $H'\cap \widehat{D}$. As 
previously we may deduce that for a small $\varepsilon_2$, $H(\rho )\cap 
\widehat{D}'' \cap B(y_1, \varepsilon_2 )={\mathcal D}H''\cap B(y_1, \varepsilon_2 )$ 
and it makes with ${\mathcal D}\widehat{H}$ the dihedral angle $\alpha $.

\begin{figure}[!ht]
\centering
\includegraphics{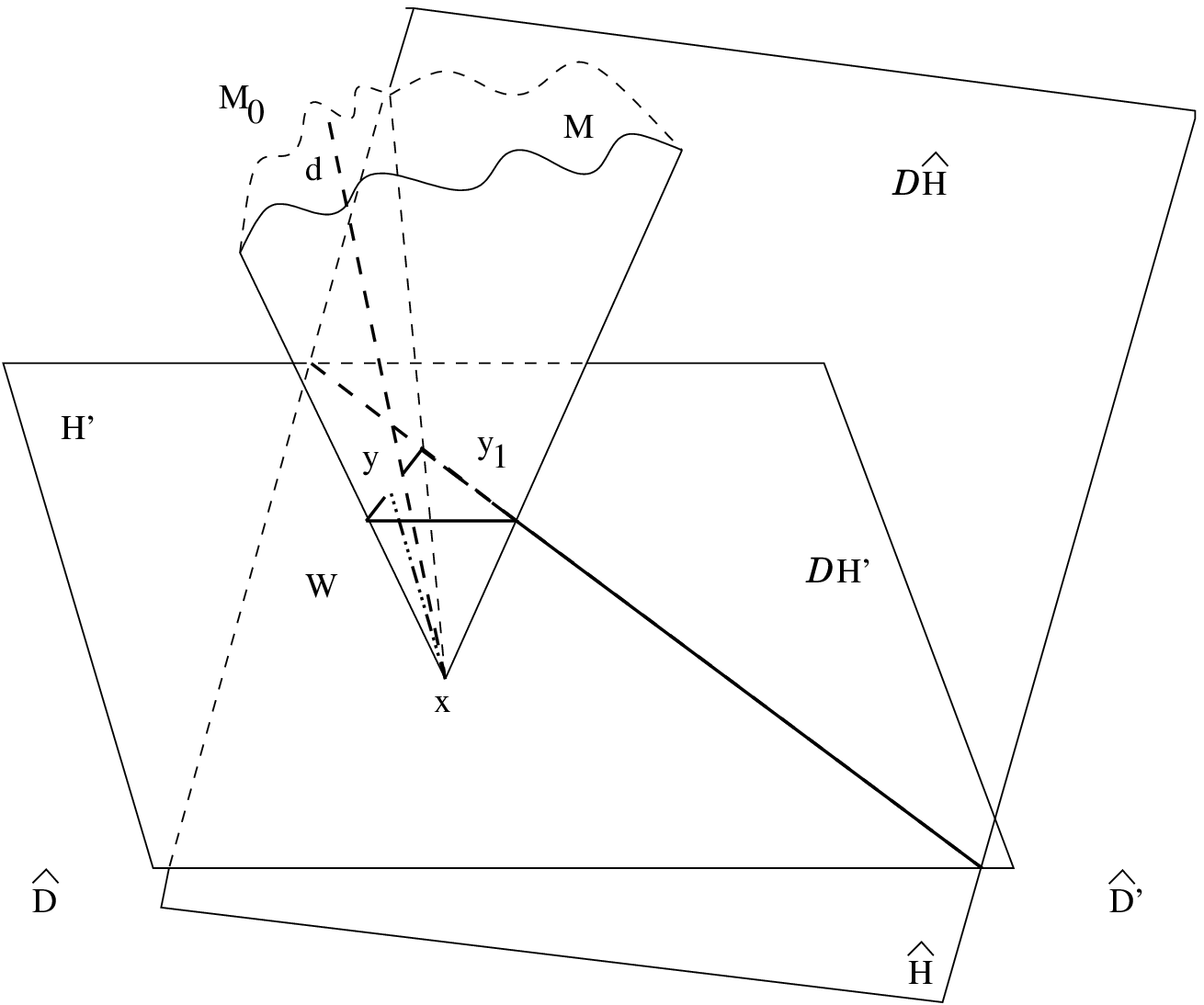} 
\caption{} 
\label{fig2}
\end{figure}

Then 
in the apartment $\widehat{D}'\cup \widehat{D}''$ we have the quadrant 
${\mathcal Q}$ bounded by ${\mathcal D}H'$ and ${\mathcal D}H''$, containing 
${\mathcal D}\widehat{H}$ in its interior, of dihedral angle $\leq 2\alpha < \pi 
$ and such that $\partial {\mathcal Q}\cap B(y_1, \varepsilon )\subset H(\rho )$ 
for $\varepsilon =\min (\varepsilon_1, \varepsilon_2 )$. It follows, by the convexity of 
the Busemann function $f_\rho $, that for a smaller $\varepsilon' $, ${\mathcal 
D}\widehat{H}\cap B(y_1, \varepsilon' )\subset {\mathcal Q}\cap B(y_1, \varepsilon' ) 
\subset Hb(\rho )$. This gives a contradiction in $F$, in which $Hb(\rho )$ is 
on the opposite side of $H'$ than ${\mathcal D}\widehat{H}$. \hspace*{\fill }$\diamondsuit $ 

\medskip

We recall that if a subset $Y$ in a metric space $(X,d_X)$ is endowed with its 
own metric $d_Y$, the metric space $(Y,d_Y)$ is said to be {\it nondistorted} in
 $(X,d_X)$ if the set $\left\{ \frac{d_Y(y_1,y_2)}{d_X(y_1,y_2)} 
\mid y_1,y_2\in Y,\; y_1\neq y_2 \right\}$ has a finite upper bound and a positive 
lower bound. An upper bound of this set is called {\it a nondistorsion constant 
of $Y$ in $X$}. An upper bound of a set of the form $\left\{ 
\frac{d_Y(y_1,y_2)}{d_X(y_1,y_2)} \mid y_1,y_2\in Y,\;  
d_X(y_1,y_2)\geq D \right\}$ is called {\it a nondistorsion constant of $Y$ in 
$X$ for sufficiently large distances}.

\medskip

The previous two Lemmata imply the following theorem 

\begin{theorem}[\cite{Dr1*}, Theorem 1.2]\label{nondist1}
Let ${\mathbf K}$ be a 3-thick Euclidean building of rank at least 2 and $\rho $ 
a geodesic ray in it. The following three properties are equivalent :

${\mathbf (P_1)}$ The horosphere $H(\rho )$ endowed with its length metric is 
nondistorted in ${\mathbf K}$ ;

${\mathbf (P_2)}$ $H(\rho )$ is connected ;

${\mathbf (P_3)}$ $\rho $ is not parallel to a rank one factor of ${\mathbf K}$.
\end{theorem}

\noindent{\bf Proof.}\quad The implications ${\mathbf (P_1)} \Rightarrow 
{\mathbf (P_2)}$ and ${\mathbf (P_2)} \Rightarrow {\mathbf (P_3)}$ are obvious. 
We show  ${\mathbf (P_3)} \Rightarrow {\mathbf (P_1)}$. We may suppose that 
$\rho $ is not parallel to any factor of ${\mathbf K}$. For if ${\mathbf 
K}={\mathbf K}_1\times {\mathbf K}_2$ and if there exists a ray $r$ in ${\mathbf 
K}_1$ and a point $x$ in ${\mathbf K}_2$ such that $\rho(t)=(r(t),x)$, then 
$H(\rho )=H(r)\times {\mathbf K}_2$, so it will be enough to prove ${\mathbf 
(P_1)}$ for $r$ in ${\mathbf K}_1$, by hypothesis ${\mathbf K}_1$ being of rank 
at least 2.

But if $\rho $ is not parallel to any factor of ${\mathbf K}$, then Lemma 
\ref{fetze} and Lemma \ref{fetze2}, (1), allow to conclude that ${\mathbf 
(P_1)}$ is true.\hspace*{\fill }$\diamondsuit $

\medskip

\begin{remark}\label{C0}
Let $H(\rho )$ be a horosphere in a 3-thick Euclidean building, ${\mathbf K}$, 
of rank at least 2, and suppose $\rho $ is not parallel to a rank one factor. 
\begin{itemize}
\item[(a)] Lemmata \ref{fetze} and \ref{fetze2} imply that we may find a 
nondistorsion constant $C_0$ of $H(\rho )$ in ${\mathbf K}$ which may be 
effectively computed when given the Coxeter complex and the slope $\theta $ of 
$\rho $. More precisely, $C_0$ may be computed given :

- the minimal dihedral angle between two distinct hyperplanes of sets of slopes $Ort\; 
(\theta )$ ;

-the dihedral angles between one hyperplane $H$ of set of slopes $Ort\; (\theta 
)$ and each of the singular hyperplanes supporting a Weyl chamber $W$ which contains a line orthogonal to $H$ (both $H$ and $W$ being considered in the same apartment).  
\item[(b)] Two points $x,y\in H(\rho )\cap F$, where $F$ is an apartment, may be 
joined either by a polygonal line of type ${\mathbf L}_{xy}$ contained in $F$ or 
a ramification of it, or by the union of a planar polygonal line in $H(\rho 
)\cap F$ joining $y$ to a point $y_1$ with a line ${\mathbf L}_{xy_1}$ contained 
in $F$ or a ramification of it. This and Proposition \ref{form}, (b), imply in 
particular that $x$ and $y$ may always be joined by a polygonal line in $H(\rho 
)$ with at most $2q_0$ edges and of length $\leq C_0d(x,y)$.
\end{itemize}
\end{remark}

\medskip

Theorem \ref{nondist1} implies the following.

\begin{theorem}[\cite{Dr1*}, Theorem 1.3]\label{nondist2} 
Let $X$ be a product of symmetric spaces and Euclidean buildings, $X$ of rank at 
least 2 and $\rho $ a geodesic ray in it. The following two properties are 
equivalent :

${\mathbf (P_1^*)}$ The horosphere $H(\rho )$ is nondistorted ;

${\mathbf (P_2^*)}$ $\rho $ is not parallel to a rank one factor of $X$.
\end{theorem} 
 
In \cite[$\S 4$]{Dr1*} there is a proof of the fact that Theorem \ref{nondist1} 
implies Theorem \ref{nondist2}.

In the sequel we shall relate the nondistorsion constant of a horosphere $H(\rho 
)$ in a product $X$ of symmetric spaces and Euclidean buildings to the Coxeter 
complex of $\partial_\infty X$ and to $\theta =P(\rho (\infty ))$, similarly to 
the case of Euclidean buildings. Also we shall show that any two points $x$, 
$y$, in a horosphere $H(\rho )$ or in the bigger subspace $X\setminus Hbo(\rho 
)$ may be joined by a curve ``almost polygonal'', with at most $3q_0$ edges.

First we introduce some new notions. Let $X$ and $\rho $ be as in the previous 
theorem, $\rho $ satisfying property ${\mathbf (P_2^*)}$. Let $x$, $y$ be two 
points in $X\setminus Hbo(\rho )$. We consider a sequence of points $x_0=x,\, 
\bar{x}_0,\, x_1,\, \bar{x}_1,\, \dots , x_p,\, \bar{x}_p=y$ satisfying the 
following properties :
\begin{itemize}
\item[(1)] for every $i\in \{0,1,\dots p-1 \}$ there exists a maximal flat $F_i$ such that $[\bar{x}_i,x_{i+1}]\subset F_i\setminus Hbo(\rho )$ ;
\item[(2)] $\Sigma_{i=0}^{p-1} d(\bar{x}_i,x_{i+1})\leq 2C_0 d(x,y)$ ;
\item[(3)] $p\leq 3q_0$ and $d(x,x_i),\; d(x, \bar{x}_i)\leq 2C_0 d(x,y),\; \forall i\in \{0,1,\dots p \}$.
\end{itemize}

Let ${\mathfrak L}_i$ be a curve of minimal length between $x_i$ and $\bar{x}_i$ 
in $X\setminus Hbo(\rho )$. The curve joining $x$ and $y$ obtained as 
${\mathfrak L}_0\cup [\bar{x}_0,x_1]\cup {\mathfrak L}_1\cup [\bar{x}_1,x_2]\cup 
\dots \cup {\mathfrak L}_{p-1}\cup [\bar{x}_{p-1},x_p]\cup {\mathfrak L}_p$ is called 
an {\it almost polygonal curve joining $x$ and $y$}. The points $x_0=x,\, 
\bar{x}_0,\, x_1,\, \bar{x}_1,\, \dots , x_p,\, \bar{x}_p=y$ are called {\it the 
vertices of the almost polygonal curve}.

 We denote 
 $$
 \epsilon(x,y):=\inf \{ \max_{i\in \{0,1,\dots p \} }d(x_i, \bar{x}_i) \mid p\in 
\N\cap [1,3q_0] \mbox{ and }x_0=x,\, \bar{x}_0,\, x_1,\, \bar{x}_1,\, \dots 
x_p,\, \bar{x}_p=y
 $$ 
\begin{center}
 vertices of an almost polygonal curve$\}$.
\end{center}

We denote 
$$
\epsilon(d):= \sup \{\epsilon(x,y)\mid x,\, y\in X\setminus Hbo(\rho ),\, 
d(x,y)=d \}.
$$

The definition of an almost polygonal curve implies that $\epsilon(x,y)\leq 4C_0 
d(x,y)$ and that $\epsilon(d)\leq 4C_0 d$.

\begin{lemma}\label{eod}
 We have that $\epsilon(d)=o(d).$
\end{lemma}

\noindent {\bf Proof.}\quad  We reason by contradiction and suppose that there 
exists a sequence of pairs of points $x_n,y_n \in X\setminus Hbo(\rho )$ with 
$d(x_n,y_n)=d_n$ and $\epsilon(x_n,y_n)>\delta d_n,\; \delta >0$. Without loss of generality we may suppose that $x_n,y_n \in H(\rho )$. In the 
asymptotic cone $X_\omega (x_n,d_n)$ let $x_\omega =[x_n],\; y_\omega =[y_n]$ 
and $\rho_\omega =[\rho ]$. According to Remark \ref{C0}, (b), $x_\omega $ and $ 
y_\omega $ may be joined in $H(\rho_\omega )$ by a polygonal line with at 
most $2q_0$ segments and of length at most $C_0$. Let $x^0_\omega 
=x_\omega ,\; x^1_\omega ,\dots , x^p_\omega =y_\omega $ be the vertices of 
this line, $p\leq 2q_0$. Each segment $[x^i_\omega, x^{i+1}_\omega]$ is 
contained in an apartment $F^i_\omega $ asymptotic to $\rho_\omega $. By Lemma 
\ref{platr}, we may suppose that $F^i_\omega =[F^i_n]$, where each $F^i_n$ is 
asymptotic to $\rho $, and that the segment $[x^i_\omega, x^{i+1}_\omega]$ 
appears as a limit of a sequence of segments $[\bar{x}^i_n, x^{i+1}_n ]\subset F^i_n 
\setminus Hbo(\rho )$. The sequence $x_n^0=x_n,\; \bar{x}^0_n, x^{1}_n, 
\bar{x}^1_n, x^{2}_n ,\dots , x^{p}_n, \bar{x}^p_n =y_n$ satisfies the 
properties (1), (2), (3) in the definition of an almost polygonal curve $\omega 
-$ almost surely. Also $\lim_\omega \frac{d(x^i_n, \bar{x}^i_n)}{d_n}=0,\; 
\forall i\in \{ 0, 1,2,\dots p \}$. This contradicts the fact that 
$\epsilon(x_n,y_n)>\delta d_n,\; \forall n\in \N $. \hspace*{\fill 
}$\diamondsuit $

\medskip

An immediate consequence of Lemma \ref{eod} is an improvement of Theorem 
\ref{nondist2}. In the theorem only the nondistorsion is stated without any 
specification on nondistorsion constants. Let ${\bf C}$ be such a constant. For two 
points $x$ and $y$ in $X\setminus Hbo(\rho )$ with $d(x,y)=d$ we have that the 
length distance in $X\setminus Hbo(\rho )$, $d_\ell (x,y)\leq 2C_0d +3q_0{\bf 
C}\epsilon(d)$. The previous inequality and Lemma \ref{eod} imply that for $d$ 
sufficiently large $d_\ell(x,y)\leq 3C_0 d(x,y)$. Thus $3C_0$ is a nondistorsion 
constant of $X\setminus Hbo(\rho )$ in $X$ for sufficiently large distances.

\begin{cor}
Let $X$ be a product of symmetric spaces and Euclidean buildings such that each 
factor is of rank at least 2. Then there exists a constant ${\bf C}_0={\bf 
C}_0({\mathsf{S}})$ such that for each geodesic ray $\rho $ in $X$ ${\bf C}_0$ is a 
nondistorsion constant of $X\setminus Hbo(\rho )$ and of $H(\rho )$ in $X$ for 
sufficiently large distances.
\end{cor}

\noindent {\bf Proof.}\quad The conclusion of the argument preceding the 
statement of the Corollary was that for each geodesic ray $\rho $ in $X$, $3C_0$ 
is a nondistorsion constant of $X\setminus Hbo(\rho )$ in $X$ for sufficiently 
large distances. We recall that $C_0=C_0(\mathsf{S},\theta )$ where $\theta =P(\rho 
(\infty ))$. The dependence of $C_0$ on $\theta $ is made explicit in Remark 
\ref{C0}, (a). From this dependence it follows that the function associating to 
each $\theta $ in $\Delta_{mod}$ the constant $C_0$ is continuous. In particular 
it has a lowest upper bound ${\bf C}_0'$. Then ${\bf C}_0 = 3{\bf C}_0'$ is a 
nondistorsion constant of $X\setminus Hbo(\rho )$ in $X$ for sufficiently large 
distances, for every geodesic ray $\rho \subset X$. \hspace*{\fill 
}$\diamondsuit $

\medskip

\begin{remark}
We can obviously generalize the previous arguments to the case of a space $X_0=X\setminus \bigsqcup_{\rho \in {\mathcal R}}Hbo(\rho )$ with ${\mathcal R}$ finite and all rays in ${\mathcal R}$ of the same slope $\theta $. The only difference in the definition of an almost polygonal curve is that  we must replace condition $p\leq 3q_0$ in (3) by $p\leq 3q_0\cdot $card ${\mathcal R}$. Lemma \ref{eod} remains true.
\end{remark}

\medskip

{\it Notations :} In the sequel we shall deal only with spaces $X_0=X\setminus \bigsqcup_{\rho \in {\mathcal R}}Hbo(\rho )$ with all $\rho \in {\mathcal R}$ of the same slope satisfying property ${\mathbf (P_2^*)}$. Suppose ${\mathcal R}$ is finite and $x,y$ are two points in $X_0$. If the ambient space $X$ is an Euclidean buiding, we shall 
denote ${\mathcal L}_{xy}$ a curve joining $x$ and $y$ in $X_0$ obtained from $[x,y]$ by replacing each subsegment $[x',y']=[x,y]\cap Hb(\rho )$ with a curve as 
described in Remark \ref{C0}, (b). If $X$ has a symmetric space as a factor, we 
shall denote ${\mathcal L}_{xy}$ an almost polygonal curve joining $x$ and $y$ with $\max_{i}d(x_i,\bar{x}_i)\leq 2\epsilon(x,y)$. We shall call such a curve a {\it minimizing almost polygonal curve joining $x$ and $y$} (though the appropriate name should probably be ``almost minimizing almost polygonal curve'').

\section{Filling in Euclidean buildings and symmetric spaces
with deleted open horoballs}\label{k}

In this section we shall prove the following result

\begin{intro}\label{fillgen}
Let $X$ be a product of symmetric spaces
and Euclidean buildings, $X$ of rank at least $3$, and $X_0$ a subset 
which can 
be written as 
$$
X_0=X\setminus \bigsqcup_{\rho \in {\mathcal R}}Hbo(\rho )\; .
$$

Suppose $X_0$ has the following properties

({\cal P}$_1$) for every point $x$ of $X_0$ there exists a
maximal flat $F\subset {\mathcal{N}}_d(X_0)$ such that $x\in {\mathcal{N}}_d(F)$,
where $d$ is a universal constant ;

({\cal P}$_2$) all rays $\rho \in {\mathcal 
R} $ have the same slope, and the common slope is 
not parallel to a rank one factor or to a rank
two factor of $X$. Then

\begin{itemize}
\item[(a)] the filling order in $X_0$ is asymptotically quadratic, that is 
$$
\forall \varepsilon >0, \; \exists \ell_\varepsilon \mbox{ such that 
}A_1(\ell 
)\leq \ell^{2+\varepsilon} \, , \, \forall \ell \geq 
\ell_\varepsilon \; ;
$$ 
\item[(b)] if the set of rays ${\mathcal R}$ is finite then the filling order in $X_0$ is quadratic, that is 
$$
A_1(\ell 
)\leq C\ell^{2} \, , \, \forall \ell , \;  
$$ where the constant $C$ depends on $X$ and on the cardinal of ${\mathcal R}$.
\end{itemize}
\end{intro}

\begin{rintro}\label{fact}
 We may suppose that the common slope of all rays $\rho \in
{\mathcal  R} $ is not parallel to any factor.
\end{rintro}

\noindent{\bf Proof.}\quad Let 
$X=X_1\times X_2$ be a decomposition as a product, let $\Delta_{mod}=\Delta_{mod}^1\circ \Delta_{mod}^2$ be the corresponding decomposition of the model chamber and suppose $\theta \in \Delta_{mod}^1$. Then $X_0=\left( X_1\setminus 
\bigsqcup_{\rho_1 \in {\mathcal R}_1} Hbo(\rho_1) \right) \times X_2$. Any 
loop can be projected onto a loop entirely contained into a copy of the 
factor $ X_1\setminus \bigsqcup_{\rho_1\in {\mathcal R}_1}
Hbo(\rho_1)$, and between the two loops there is a filling cylinder of quadratic surface. Thus, it suffices to prove the result in $ X_1\setminus \bigsqcup_{\rho_1\in {\mathcal R}_1}
Hbo(\rho_1)$.\hspace*{\fill }$\diamondsuit $

\medskip

We note that actually we can specify the order of the
filling in
$X_0$ in many of the cases when 
$\theta $ is parallel to a factor of rank at most $2$. With and argument as in the proof of the previous remark we reduce the problem to the case when $X$ is itself of rank at most $2$. Then the filling order in $X_0$ is quadratic if $X=\hip_{\R }^n,\; n\geq 3$, or $X=\hip_{\C }^n, \; n\geq 3$, or cubic if $X=\hip_{\C }^2$, exponential if $X$ is a symmetric space  of rank
two, and linear if $X$ is a tree or $X=\hip_{\R }^2$.

\smallskip

Before proving Theorem \ref{fillgen} we prove the following intermediate result.

\begin{intro}\label{fillK0}
Let ${\mathbf K}$ be a $4$-thick
Euclidean building of rank at least $3$ and ${\mathbf
K}_0$ a subspace of the form 
$$
{\mathbf K}_0={\mathbf K} \setminus \bigsqcup_{\rho \in
{\mathcal R}_\omega}Hbo(\rho )\; .\leqno(4.1)
$$

Suppose ${\mathbf K}_0$ has the following properties

(${\mathcal P}_1$) through every point of ${\mathbf K}_0$
passes an apartment entirely contained in ${\mathbf K}_0$ ;

(${\mathcal P}_2$) all rays $\rho \in
{\mathcal  R}_\omega $ have the same slope $\theta $ which is 
not parallel to a rank one factor or to a rank
two factor of ${\mathbf K}$. Then
\begin{itemize}
\item[(a)] the filling order in ${\mathbf K}_0$ is at most cubic, that is 
$$
A_1(\ell )\leq C\cdot \ell^3 \, , \, \forall \ell>0\;  ;
$$ 
\item[(b)] if the set of rays ${\mathcal R}_\omega$ is finite every loop $\frak C$ in ${\mathbf K}_0$ composed of at most $m$ segments and of length $\ell $ has 
$$
A_1({\frak C})\leq C\cdot \ell^2 \, , 
$$
 where $C=C(\theta , \mathsf{S}, m)$. 
\end{itemize} 
\end{intro}

The proof of Theorem \ref{fillK0} is done in several steps. First we
look at loops contained in just one apartment of the Euclidean building.

{\it Notation}: Henceforth for a curve $\frak L$ without self-intersection and two points $x,y$ on it we shall denote by ${\frak L}_{xy}$ the arc on $\frak L$ of extremities $x$ and $y$.

\subsection{Loops contained in one apartment}\label{cbpl}

\begin{proposition}\label{polit3}
Let ${\mathbf K}$ be a $4$-thick
Euclidean building of rank at least $3$, $F$ an apartment in it, $\rho $ a
ray of slope $\theta $ not parallel to a rank two factor and ${\frak C} :\sph^1\to F\setminus Hbo(\rho )$ a loop of length $\ell $. 
Then there exists a positive constant $L=L(\theta ,\mathsf{S})$ such that the filling area of the loop ${\frak C}$ in  ${\mathbf K}\setminus Hbo(\rho )$ satisfies
$$
A_1(\frak C)\leq L\cdot \ell^2\; .
$$    
\end{proposition}

\noindent{\bf Proof.}\quad We show that we may suppose ${\frak C}(\sph^1)\subset F\cap H(\rho )$. If a filling disk obtained by joining a fixed point of ${\frak C}(\sph^1)$ with all the other points does not intersect $F\cap Hbo(\rho )$ then we are done. We suppose therefore the contrary which implies that ${\frak C}(\sph^1)$ is in the $2\ell $ neighborhood of $F\cap Hb(\rho )$. By Corollary \ref{prf} we may project $\frak C$ on $F\cap H(\rho )$ and obtain a curve ${\frak C}' : \sph^1 \to F\cap H(\rho )$ of length at most $\ell $. Since ${\frak C}(\sph^1) \subset {\mathcal{N}}_{2\ell} (F\cap Hb(\rho )) $ the segments along which we project the curve $\frak C$ on $\frak C'$ form a cylinder in $F\setminus Hbo(\rho )$ with area of order $\ell^2$. So it suffices to fill ${\frak C}'$. Thus we may suppose from the beginning that ${\frak C}(\sph^1)\subset F\cap H(\rho )$.

With an argument as in the proof of Remark \ref{fact} we may reduce to the case when $\rho $ is not parallel to any factor. By hypothesis, we may suppose that the rank of ${\mathbf K}$ is not $2$. If ${\mathbf K}$ has rank one then it is an $\R$-tree, any horosphere is totally disconnected in it so any loop in $H(\rho )$ reduces to a point. In the sequel we suppose that ${\mathbf K}$ is of rank at least $3$.

If $\inf_{x\in F}f_\rho(x)=-\infty $ the result has been 
proven in \cite{Dr2*}, in the proof of Proposition 4.3, in which it appears as 
case (1). So in the sequel we suppose $\inf_{x\in F}f_\rho(x)=-m > - \infty $. 
Let $\Phi$ be the affine span of $Min_F(\rho )$. The idea of the proof is to fix 
the loop but eventually to change the apartment $F$ with a more convenient one, 
in which $Min_F(\rho )$ is at a Hausdorff distance of order $\ell $ from a 
polytope of codimension $3$. More precisely, we prove that by eventually 
changing the apartment $F$ we may suppose that either $Min_F(\rho )$ is of codimension $3$ or $Min_F(\rho )$ is of codimension $2$ but contained in a codimension one $k\ell$ -strip or $Min_F(\rho )$ is of codimension $1$ but contained in a codimension two $(\epsilon, k\ell)$ -strip, where $k$ and $\epsilon$ are constants depending only on $\theta 
$ and $\mathsf{S}$. See Definition \ref{strip} for the different notions of strips. This already happens if $\Phi$ has codimension at 
least $3$. So in the sequel we suppose $\Phi $ has codimension one or two. By 
Corollary \ref{cskew} there exists a singular hyperplane $H_1$ which neither contains 
nor is orthogonal to $\Phi $. For technical reasons we suppose $H_1$ intersects 
$\frak C(\sph^1)$. For every $x\in \frak C(\sph^1)$ we consider $H_1(x)$ the singular hyperplane  
through $x$ parallel to $H_1$, and we consider the strip ${\mathcal 
S}:=\bigcup_{x\in \frak C(\sph^1)}H_1(x)$. There are two cases : the strip  ${\mathcal 
S}$ intersects or does not intersect the relative interior of $Min_F(\rho )$. We 
show that in both cases we may change $F$ in such a way that $Min_F(\rho )$ either has one dimension less or is contained in a codimension one strip.

If the strip ${\mathcal S}$ does not intersect the relative interior of 
$Min_F(\rho )$ then there exists a hyperplane $H_1'$ parallel to 
$H_1$ such that

$\bullet$  $\frak C(\sph^1)$ is contained in one half-appartment $D$ of $F$ bounded by 
$H_1'$ while $Min_F(\rho )$ is contained in the opposite half-apartment $D'$ ;

$\bullet$  $H_1'\cap Min_F(\rho )$ has codimension at least $1$ in $\Phi$.

By Lemma \ref{hipp} there exists a half-apartment $D_1$ of boundary $H_1'$ and 
interior disjoint of $F$ such that $Min_{D'\cup D_1}(\rho )=Min_{D'}(\rho )$. 
Then $Min_{D\cup D_1}(\rho )=Min_{H_1'}(\rho )$ has codimension one in $\Phi $. 
By replacing the initial flat $F$ with $D\cup D_1$ we may therefore increase the 
codimension of $ Min_F(\rho )$ with 1.

 Suppose the strip ${\mathcal S}$ intersects the relative interior of 
$Min_F(\rho )$. Let $H_1^a$ and $H_1^b$ be the extremal hyperplanes 
of this strip. The distance between them is at most $\ell $. By applying Lemma 
\ref{hipp} twice we may suppose that $Min_F(\rho )$ is entirely contained in 
the strip. That is, $Min_F(\rho )$ is contained in the strip determined by 
$H_1^a\cap \Phi $ and $H_1^b\cap \Phi $ in $\Phi $. Since there is a finite 
number of possibilities for the dimension of $\Phi$ and the angle between $\Phi 
$ and $H_1$, there exists a constant $\kappa =\kappa(\theta , \mathsf{S})$ such that $H_1^a\cap \Phi $ and $H_1^b\cap \Phi $ are at 
distance at most $\kappa \ell$ from each other. Thus $Min_F(\rho )$ is contained in a codimension one $k\ell$ -strip.

If $\Phi $ has codimension $2$, this finishes the proof. So the only case left is 
when $\Phi $ has codimension $1$, so it is a singular hyperplane. To finish the 
proof in this case it is enough to find another singular hyperplane $H_2$ in $F$ 
which neither contains nor is orthogonal to $\Phi $ and such that $\Phi \cap 
H_1\cap H_2$ has codimension $2$ in $\Phi $. We note that the hyperplane $H_1$ 
remained in the flat $F$ even after one of the two possible changes presented previously was 
performed, so the latter condition makes sense.

Let $W$ be a Weyl chamber adjacent both to $H_1$ and to $H_1\cap \Phi$ and let 
$M=W\cap H_1\cap \Phi$. If there exists a hyperplane $H'$ supporting $W$ which 
doesn't contain $M$ and which is not orthogonal to $\Phi$ we take $H_2=H'$. We 
suppose then that all hyperplanes supporting $W$ and not containing $M$ are 
orthogonal to $\Phi$. By Corollary \ref{cskew} one of these hyperplanes, $H''$, is 
not orthogonal to $Span\, M = H_1\cap \Phi $. Since $H''$ is orthogonal to the 
hyperplane  $\Phi $, it follows that $H''$ is not 
orthogonal to $H_1$. In this case we take as $H_2$ the image of $H_1$ by 
orthogonal symmetry with respect to $H''$. Since $H_1$ is not orthogonal to 
$\Phi $ and does not contain $\Phi$, the same is true for $H_2$. Also $H_1\cap 
H_2 $ coincides with $ H_1\cap H''$ so it differs from $H_1\cap \Phi$.

By repeating with $H_2$ the argument done previously with $H_1$ we obtain a flat $F'$ containing the image of the loop $\frak C$ such that $Min_{F'}(\rho )$ is contained in a codimension two $(\epsilon, k\ell)$-strip in $\Phi $.

By 
Corollary 
\ref{vec}, ${\mathcal{N}}_m(Min_F(\rho ))\subset Hb(\rho )\cap F \subset {\mathcal{N}}_{am}(Min_F(\rho 
)).$ Lemma \ref{proj}, (b), in the Appendix implies that the projection $p$ of $H(\rho )\cap F$ onto $\partial {\mathcal{N}}_m(Min_F(\rho ))$ is bilipschitz with respect to the 
length metrics, the bilipschitz constant depending only on $a$. According to Proposition \ref{lpolit} in the Appendix the filling area 
of the loop $p\circ \frak C$ is at most $L(1+k+k^2)\ell^2$, so the filling area of $\frak C$ is 
bounded by $L'\ell^2$, with $L'$ depending only on the constants $a$ and $k$, so on the 
Coxeter complex $\mathsf{S}$ of ${\mathbf K}$ and on $\theta $.\hspace*{\fill }$\diamondsuit $

\medskip

\subsection{Loops contained in two apartments}\label{two}

We shall show that loops in the exterior of a horoball which are composed of two arcs, each arc contained in one apartment, have a quadratic filling area outside the horoball. We also have to add a condition on the two apartments : they must have at least a Weyl chamber in common. It is possible that a similar result is true for loops contained in the union of two apartments without any other restriction. But for our purpose the previous result suffices.

First we prove an intermediate result.

\begin{lemma}\label{demi}
Let ${\mathbf K}$ be a $4$-thick
Euclidean building of rank at least $2$, let $F, \; F'$ be two apartments having a half-apartment $D_0$ in common and let $\rho $ be a
ray of slope $\theta $ not parallel to a rank two factor. Let $x$ and $y$ be two points in $D_0\setminus Hbo(\rho )$ and ${\frak C}_{xy}\subset F\setminus Hbo(\rho )$, ${\frak C}_{xy}'\subset F'\setminus Hbo(\rho )$ two arcs without self-intersections joining $x$ and $y$, of length at most $\ell $.

Then there exists a positive constant $L=L(\theta ,\mathsf{S})$ such that the filling area of the loop ${\frak C}={\frak C}_{xy}\cup {\frak C}_{xy}'$ in  ${\mathbf K}\setminus Hbo(\rho )$ satisfies
$$
A_1(\frak C)\leq L\cdot \ell^2\; .
$$ 
\end{lemma}

\noindent{\bf Proof.}\quad {\it Notations }: We denote $H:=\partial D_0$. For every two points $a,b\in {\frak C}_{xy}$ we denote ${\frak C}_{ab}$ the arc of ${\frak C}_{xy}$ between $a$ and $b$. We use a similar notation for ${\frak C}_{xy}'$.

\smallskip

{\sc Step 1.}\quad First we show that we can reduce to the case where $x,y\in D_0\cap H(\rho )$ and ${\frak C}_{xy},\, {\frak C}_{xy}'$ are the shortest polygonal curves joining $x$ and $y$ in $ F\cap H(\rho )$ and $ F'\cap H(\rho )$, respectively. Let $[x',y']=[x,y]\cap Hb(\rho)$ and let $\overline{{\frak C}}_{x'y'},\, \overline{{\frak C}}_{x'y'}'$ be the shortest polygonal curves joining $x'$ and $y'$ in $ F\cap H(\rho )$ and $ F'\cap H(\rho )$, respectively. The loop ${\frak C}_{xy}\cup [x,x']\cup \overline{{\frak C}}_{x'y'} \cup [y',y]$ is contained in $F\setminus Hbo(\rho )$ and has length of order $\ell $. By Proposition \ref{polit3} it can be filled with an area of order $\ell^2$. Likewise for the loop ${\frak C}_{xy}'\cup [x,x']\cup \overline{{\frak C}}'_{x'y'} \cup [y',y]$. So it suffices to fill the loop $\overline{{\frak C}}_{x'y'}\cup \overline{{\frak C}}_{x'y'}'$ with an area of order $\ell^2$. 

\smallskip

{\sc Step 2.}\quad We show that we may reduce the problem to the case where $x,\; y \in H$.

If one of the two arcs ${\frak C}_{xy},\;  {\frak C}_{xy}'$ is entirely contained in $D_0\setminus Hbo(\rho )$ Proposition \ref{polit3} allows to conclude. So in the sequel we may suppose that both have points in the half-apartments opposite to $D_0$. Let $D$ and $D'$ be the half-apartments opposite to $D_0$ in $F$ and $F'$ respectively. Let $x_1$ be the nearest point to $x$ of ${\frak C}_{xy}\cap H$.

  Suppose there exists a point $x'\in {\frak C}_{x,x_1}\cap Int\; D_0$ contained into two faces of $F\cap H(\rho )$. We consider the nearest $x'$ to $x$. Then ${\frak C}_{xy}=[x,x']\cup {\frak C}_{x'y}$. By eventually adding some area of order $\ell^2$, we may replace ${\frak C}_{x'y}$ with the planar polygonal line ${\bf L}_{x'y}$. We may suppose that ${\bf L}_{x'y}$ intersects $H$ otherwise we finish by applying Proposition \ref{polit3}. We can construct a copy of ${\bf L}_{x'y}$ in $F'$, ${\bf L}_{x'y}'$. Both lines intersect $H$ in the same point $x''$. Let ${\bf L}$ and ${\bf L}'$ be the arcs of ${\bf L}_{x'y}$ and ${\bf L}_{x'y}'$, respectively,  of extremities $x''$ and $y$. Proposition \ref{polit3} implies that the loop ${\bf L}\cup {\bf L}'\subset [D\cup D']\setminus Hbo(\rho )$ may be filled with an area of order $\ell^2$. Then the same follows for the loop ${\bf L}_{x'y}\cup {\bf L}_{x'y}'$. So now it suffices to fill the loop $[x,x']\cup {\bf L}_{x'y}'\cup {\frak C}_{xy}'$. But this loop is contained in $F'\setminus Hbo(\rho )$ so we use Proposition \ref{polit3} once more.

 Thus, in the sequel we may suppose there exists no point $x'\in {\frak C}_{x,x_1}\cap Int\; D_0$ contained into two faces of $F\cap H(\rho )$. This means that either $x\in H$ or ${\frak C}_{x,x_1}=[x,x_1]$ which is a segment contained in the unique face of $F\cap H(\rho )$ through $x$. Let $x_1'$ be the nearest point to $x$ of ${\frak C}_{xy}'\cap H$. With a similar argument we conclude that we may suppose that either $x\in H$ or ${\frak C}_{x,x_1'}'=[x,x_1']$ is contained in the unique face through $x$. If $x\not\in H$ then by joining $x_1$ and $x_1'$ with a segment in the face containing both we separate the loop ${\frak C}$ into an Euclidean triangle contained in the face and another loop of length of order $\ell $. It follows that, up to adding a term of order $\ell^2$, we may suppose $x\in H$. A similar argument allows to conclude the same for $y$.

\smallskip

{\sc Step 3.}\quad We show that by eventually adding some area of order $\ell^2$, we may suppose ${\frak C}_{xy}$ is entirely contained into one half apartment bounded by $H$. The same argument done for ${\frak C}_{xy}'$ and Proposition \ref{polit3} imply then the conclusion of the Lemma.

Let $x_2$ be the nearest point to $x$ on ${\frak C}_{xy}$ contained in a face non-parallel to a face through $y$. If $x_2=x$, we may replace ${\frak C}_{xy}$ with a planar polygonal line ${\bf L}_{xy}$ and we are done. So we may suppose in the sequel that $x_2\neq x$. Then it follows that $x$ and $y$ are in the interior of two parallel faces and that ${\frak C}_{xx_2}$ coincides with the segment $[x,x_2]$ entirely contained in the face through $x$. We may replace ${\frak C}_{x_2y}$ by ${\bf L}_{x_2y}$. If ${\bf L}_{x_2y}$ intersects $H$ only in $y$ then the curve ${\frak C}_{xy}$ thus modified is contained in one half-apartment. Suppose ${\bf L}_{x_2y}\cap H =\{x_3,\, y \}$. We now look at the convex polytope $H\cap Hb(\rho )$ whose boundary contains the points $x,y,x_3$. The point $x_3$ may be joined either to $x$ or to $y$ by a polygonal line of length of order $\ell $. This is because $H\cap Hb(\rho )$ cannot have three distinct parallel faces. Thus we replace either ${\frak C}_{xx_3}$  or ${\frak C}_{x_3y}$ with a polygonal line in $H\cap H(\rho )$. In both cases we obtain a new curve contained in one half-apartment.\hspace*{\fill }$\diamondsuit $

\medskip

We now prove the more general statement we need.

\begin{proposition}\label{camera}
Let ${\mathbf K}$ be a $4$-thick
Euclidean building of rank at least $3$ and let $\rho $ be a
ray of slope $\theta $ not parallel to a rank two factor. Let $F, \; F'$ be two apartments having at least a Weyl chamber $W_0$ in common. Let $[x,y]$ be a regular segment with extremities in $W_0\setminus Hbo(\rho )$ and ${\frak C}_{xy}\subset F\setminus Hbo(\rho )$, ${\frak C}_{xy}'\subset F'\setminus Hbo(\rho )$ two arcs joining $x$ and $y$ of length at most $\ell $.

Then there exists a positive constant $L=L(\theta ,\mathsf{S})$ such that the filling area of the loop ${\frak C}={\frak C}_{xy}\cup {\frak C}_{xy}'$ in  ${\mathbf K}\setminus Hbo(\rho )$ satisfies
$$
A_1(\frak C)\leq L\cdot \ell^2\; .
$$ 
\end{proposition}

\noindent{\bf Proof.}\quad To simplify we may suppose that $W_0$ has vertex $x$ and that $[x,y]\subset W_0$. As in Step 1 of the previous proof, we show that we may suppose that $x,y\in W_0\cap H(\rho ),\, (x,y)\subset Hbo(\rho ),$ and that ${\frak C}_{xy}$ and ${\frak C}_{xy}'$ are polygonal lines of minimal length joining $x,y$ in $F\cap H(\rho )$ and $F'\cap H(\rho )$, respectively.

 We shall prove the conclusion as follows. Let $W^{op}_0$ be the Weyl chamber of vertex $x$ opposite to $W_0$ in $F'$. We shall choose a minimal gallery of Weyl chambers in $F'$ stretched between $W_0$ and $W^{op}_0$, $\overline{W}_0=W_0,\, \overline{W}_1,\dots , \overline{W}_{p_0}=W^{op}_0$ and correspondingly we shall construct a sequence of apartments $\overline{F}_0=F, \overline{F}_1,\dots , \overline{F}_{p_0}=F'$ such that each $\overline{F}_{i+1}$ is a ramification of $\overline{F}_i$ containing $\overline{W}_{i+1}$. Such a sequence of apartments always exists by Corollary \ref{ramif}. The choice of the gallery will ensure that in each $\overline{F}_i$ the points $x$ and $y$ may be joined by a polygonal line with length of order $d(x,y)$. This and the repeated application of Lemma \ref{demi} will imply the conclusion.

The inclusion $(x,y)\subset Hbo(\rho )$ implies $\angle_x(\overline{\rho_x}, (W_0)_x)<\angle_x(\overline{\rho_x}, \overline{xy})<\frac{\pi }{2}$. There are several cases, which we shall denote as in Lemma \ref{fetze2}. 

\smallskip

$(a)$ Suppose $(W_0)_x$ does not contain $\overline{\rho_x}$. By Corollary \ref{sim} there exists a Weyl chamber $\widehat{W}$ adjacent to $W_0$ with the property that in any apartment containing $W_0\cup \widehat{W}$ $x$ is contained into two faces of the trace of $H(\rho )$. According to Corollary \ref{ramif} there exists a ramification $F''$ of $F'$ containing $W_0\cup \widehat{W}$. We apply Lemma \ref{demi} to the loop composed of ${\frak C}_{xy}'$ and of a polygonal line ${\bf L}_{xy}\subset F''$ and conclude that, up to an additional quadratic filling area, we may suppose from the beginning that $F'$ contains $W_0\cup \widehat{W}$. We choose our gallery from $W_0$ to $W_0^{op}$ such that $\overline{W}_1=\widehat{W}$. This implies that each apartment $\overline{F}_i$ for $i>0$ contains $W_0\cup \widehat{W}$ so $x$ is contained into two faces of $\overline{F}_i\cap H(\rho )$. In particular in each $\overline{F}_i$ the points $x$ and $y$ may be joined outside $Hbo(\rho )$ by a polygonal line ${\bf L}_{xy}^i$.

\smallskip

$(b_1)$ Suppose $\overline{\rho_x}\in (W_0)_x$ and the connected component of $W_0\cap H(\rho )$ containing $y$ has at least two faces. Then we may choose any gallery of Weyl chambers and corresponding sequence of apartments. Any apartment $\overline{F}_i$ will contain $W_0$ and by the previous hypothesis and Lemma \ref{fetze2} $x$ and $y$ shall be joined in $\overline{F}_i\cap H(\rho )$ by a polygonal line of length of order $d(x,y)$.

\smallskip

$(b_2)$ Suppose $\overline{\rho_x}\in (W_0)_x$ and the connected component of $W_0\cap H(\rho )$ containing $y$ has only one face. Lemma \ref{fetze2}, $(b_2)$, implies that in a ramification $F''$ of $F'$ all the points in the intersection of a panel $M$ of $W_0$ with $H(\rho )$ are in two different faces of $F''\cap H(\rho )$. By eventually aplying once Lemma \ref{demi}, we may suppose that $F'$ has this property itself. We choose our gallery from $W_0$ to $W_0^{op}$ such that $\overline{W}_1\cap W_0 = M$. Since each apartment $\overline{F}_i$ for $i>0$ contains $W_0\cup \overline{W}_1$, in each of these apartments $\overline{F}_i$ the points in $Int\; M\cap H(\rho )$ are in two different faces of $\overline{F}_i\cap H(\rho )$. According to Lemma \ref{fetze2}, $(b_2)$, in each apartment $\overline{F}_i$ $x$ and $y$ may be joined in $\overline{F}_i\cap H(\rho )$ by a polygonal line of length of order $d(x,y)$. \hspace*{\fill }$\diamondsuit $    

\subsection{Proof of Theorem \ref{fillK0}}

An important tool in our argument is 
the set of ``good
slopes'' with respect to the slope of a ray defining a
horoball. These ``good slopes'' are the slopes which are
transverse to the horosphere defined by the given ray.

\begin{lemma}[\cite{Dr2*}, Lemmata 4.9, 4.10]\label{beta}
(1)  Let $\Sigma$ be a labelled spherical building, and $\theta
$ a point in $ \Delta_{mod} $. For every small $\delta_1$ there
exists a continuum
 of points $\beta \in \Delta_{mod} $ such that $d(\beta
,Ort\; (\theta ))>\delta_1$.

(2)  Let ${\mathbf K}$ be an Euclidean building and $\rho $ a
ray of slope $\theta $. If a geodesic segment $\lbrack
x,y\rbrack $ has slope $\beta $ as in (1), then $f_\rho $ decreases or
increases on it with a rate at least equal to $\sin \delta_1$. 
\end{lemma}

 If $\theta \in \Delta_{mod} $ and $\delta_1 >0$ are fixed, we call slopes 
$\beta \in 
\Delta_{mod}
$ verifying the condition in Lemma \ref{beta}, (1), $\delta_1$-{\it good slopes 
with respect to} $\theta $. Slopes
which moreover verify 
$d(\beta
,\partial \Delta_{mod})>\delta_1$ are called $\delta_1 $-{\it good regular slopes with 
respect 
to} $\theta $. Whenever there is no possibility of confusion we omit $\theta $ and $\delta_1$.

\bigskip

\noindent{\bf Proof of Theorem 
\ref{fillK0}.}\quad Since the arguments for the proofs of (a) and (b) follow the same lines, we shall prove them simultaneously, specifying the differences each time it is the case. We shall say that the case when ${\mathcal R}_\omega $ is finite is the finite case and that otherwise we are in the general case. If ${\mathcal R}_\omega $ is finite then let ${\frak C}:\sph^1\to {\mathbf K}_0$ be 
a lipschitz loop of length $\ell $ constituted by at most $m$ segments. Otherwise by ${\frak C}:\sph^1\to {\mathbf K}_0$ we mean 
an arbitrary lipschitz loop of length $\ell $.

\medskip

{\sc Step 1.}\quad First we construct a filling disk of quadratic 
area in ${\mathbf K}$, which we 
shall deform afterwards in order to avoid the deleted horoballs. More precisely, 
we shall construct a 
cylinder over ${\frak C}$ in a regular direction and cut it at a height of 
order $\ell $, making sure that the top of the cylinder is entirely contained in 
${\mathbf K}_0$. The advantage of this construction with respect to others 
(like, for 
instance, joining one fixed point on the image of ${\frak C}$ with all the others) is 
that when deforming the filling disk to avoid the deleted horoballs we shall 
have to deal with segments which have all the same slope, the 
choice of which is left to us.

First we shall choose a finite set of points on $\frak C (\sph^1)$. Let $P_0$ be a point on $\frak C (\sph^1)$ which we choose to be extremity of a segment in the finite case. Let $F$ be an 
apartment through 
$P_0$ entirely contained in ${\mathbf{K}}_0$ and $b$ a point in $F(\infty )$ 
such that its projection on $\Delta_{mod} $ is a $\delta_1$-good 
regular 
slope 
$\beta $ with respect to $\theta $. Then $b$ is contained into an unique 
spherical chamber $\Delta_0$.

\begin{figure}[!ht]
\centering
\includegraphics{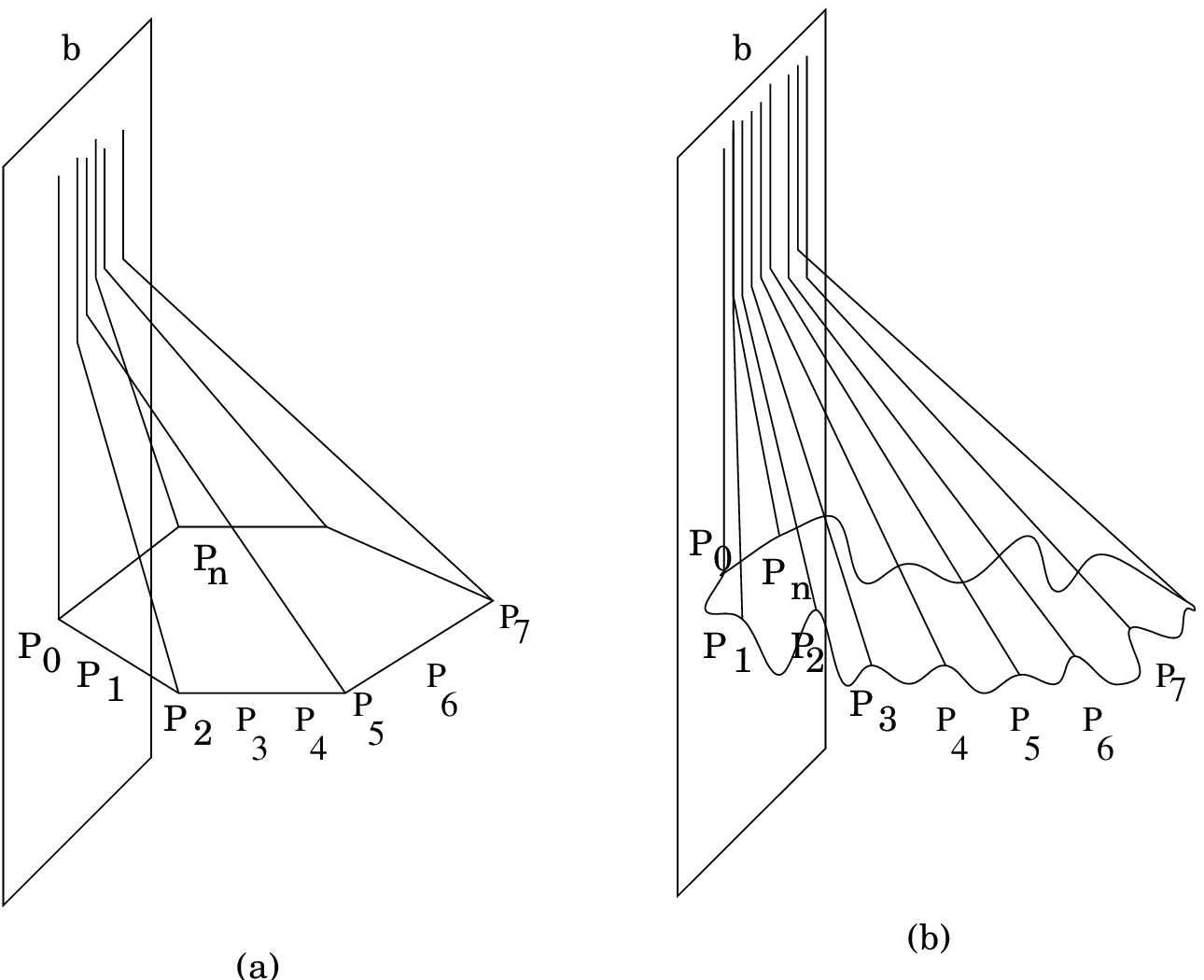}
\caption{} 
\label{fig3}
\end{figure}

In the general case we fix a small $\lambda >0$ and choose a 
finite sequence $P_0,P_1,\dots P_n$ of points, with $n\leq \frac{2\ell 
}{\lambda }$, which 
determin a partition of the image of $\frak C$ into arcs of length at most $\lambda 
$. We consider the rays $r_k=\lbrack P_k ,b)$ (Figure \ref{fig3}) and apartments $F_k$ containing 
these 
rays, where $F_0=F$. Each $F_k$ has at least a Weyl chamber in common with 
$F$ (this Weyl chamber having as boundary at infinity $\Delta_0$).

In the finite case let  $Q_0=P_0,Q_1,\dots Q_j,\; j\leq m$, be the extremities of the segments. Each segment $[Q_i,Q_{i+1}]$ is contained in an apartment $F_i'$. By applying Proposition \ref{form}, (a), to the flat $F_i'$ and to the ray $[P_0,b)$ we obtain a partition of $[Q_i,Q_{i+1}]$ into at most $q_0$ segments, each contained in an apartment asymptotic to $[P_0,b)$. In the end we obtain $P_0,P_1,\dots P_n,\; n\leq mq_0$, points on $\frak C$ such that each $[P_k,P_{k+1}]\subset {\frak C}$ is contained in an apartment asymptotic to $[P_0,b)$ (Figure \ref{fig3}). We denote the apartment corresponding to $[P_k,P_{k+1}]$ by $F_k$. It contains both $r_k=[P_k,b)$ and $r_{k+1}=[P_{k+1},b)$. We note that in this case $n$ is uniformly bounded while in the general case it was of order $\ell $. 

We know that $d(P_k,P_0)\leq \ell ,\;
\forall  k\in \{1,2,\dots n \}$, so 
$d(P_k,F_0)\leq \ell ,\; \forall k$. By Lemma 4.6.3, \cite{KlL*}, all 
points $M$ of $r_k$ with $d(M,P_k)\geq \frac{d(P_k,F_0)}{\sin 
\delta_1}$ are 
contained in $F_0$. So $r_k$ is contained in $F_0$ at least starting 
from the 
point 
$M_k$ with $d(M_k,P_k)=\frac{\ell }{\sin \delta_1}$. We note that 
$d(P_0,M_k)\leq 
\ell \left(1+\frac{1}{\sin \delta_1} \right)$. Then the hyperplane in $F_0$ 
orthogonal to 
$r_0$ at distance $\ell \left(1+\frac{1}{\sin \delta_1} \right)$ from 
$P_0$ 
intersects each $r_k$ in a point $N_k$. The convexity of the distance function 
implies 
that 
$d(N_i,N_j)\leq d(P_i,P_j),\forall i\neq j$. The closed polygonal curve ${\frak C}'$ 
with 
vertices 
$N_0,N_1,N_2,\dots , N_n$ has length at most $\ell $ and is entirely 
contained in 
$F_0$ so in ${\mathbf K}_0$. The filling area of ${\frak C}'$ in
${\mathbf  K}_0$ is 
then quadratic. So, to end the proof, it suffices to find a ``filling 
cylinder'' 
between the loops ${\frak C}$ and ${\frak C}'$ of the desired area in
${\mathbf  K}_0$. We note 
that each 
segment $\lbrack P_k,N_k\rbrack$ has length smaller than 
$d(P_k,P_0)+d(P_0,N_0)+d(N_0,N_k)\leq \ell \left(3+\frac{1}{\sin 
\delta_1} 
\right)$ 
and it has its ends in ${\mathbf K}_0$, but it is not necessarily entirely 
contained in 
${\mathbf K}_0$.

\bigskip

\noindent{\sc Step 2.}\quad An auxiliary construction is needed in the general case. In this case $\lbrack P_k ,N_{k} \rbrack $ and $\lbrack P_{k+1} ,N_{k+1} \rbrack 
$ are not in the same apartment asymptotic to $r_0$. We remedy this by constructing, for each segment $\lbrack P_{k+1} ,N_{k+1} \rbrack 
$ a ``copy'' of it contained in $F_k$ and coinciding with $\lbrack 
P_{k+1} ,N_{k+1} \rbrack $ on most of its length.

Let 
$R_{k+1}$
be the first point in which the ray $r_{k+1}$ meets the apartment 
$F_k$.
We have that $d(P_{k+1},R_{k+1})\leq c\lambda $ by Lemma 4.6.3, \cite{KlL*}. If 
$R_{k+1}$ is 
contained
in
$F_k\cap {\mathbf K}_0$, then we denote it by
$P_{k+1}'$. If it is contained in some $Hbo(\rho )\cap F_k$, then, as 
$f_\rho$ is a contraction and $f_\rho(P_{k+1})\geq 0$, $f_\rho(R_{k+1})\geq 
-c\lambda $. 
Since the ray
$r_{k+1}$ has a good  slope,
by continuing towards $b$ or by going in $F_k$ in the opposite 
direction a
distance at most $c_1\lambda$ away, we meet $H(\rho )\cap F_k$. We denote 
this point
of intersection of the ray $r_{k+1}\cap F_k$ (or of its opposite in $F_k$) with 
$H(\rho )\cap F_k$ by
$P_{k+1}'$. In both cases we have $d(P_{k+1},P_{k+1}')\leq 
c_2\lambda$, with $c_2=c+c_1$, which implies that $d(P_{k},P_{k+1}')\leq 
c_3\lambda $, with $c_3=c_2+1$. To simplify ulterior arguments we make the convention that in the finite case, $P_{k+1}'=P_{k+1}$.

\bigskip

\noindent{\sc Step 3.}\quad We replace simultaneously the pair of segments $\lbrack P_k ,N_k\rbrack $ and $\lbrack 
P_{k+1}' ,N_{k+1}\rbrack $, both included in $F_k$, by curves in ${\mathbf K}_0$ with length of order 
$\ell$. We do this by deforming them in order to avoid each horoball they intersect. There are three possible situations : a horoball $Hb(\rho )$ intersects both segments into subsegments $[x,y]\subset \lbrack P_k ,N_k\rbrack $ and $[x',y']\subset \lbrack P_{k+1}' ,N_{k+1}\rbrack $ or it intersects only $\lbrack P_k ,N_k\rbrack $ or it intersects only $\lbrack P_{k+1}' ,N_{k+1}\rbrack $. The last two situations being symmetric we shall discuss only one of them.

First we consider the first situation. We shall always suppose that $x$ and $x'$ are the extremities 
which are nearest to $P_k$ and $P_{k+1}'$, respectively.

\medskip

Lemma \ref{fetze} implies that $x$ and $x'$ may be joined by a polygonal line ${\mathbf{L}}_{xx'}$ in 
$F_k\cap H(\rho )$ of length at most $Cd(x,x')$. This 
is because $x$ and $x'$ cannot be contained in the interiors of two distinct and 
parallel faces, as the spans of these faces mustn't separate points in the set 
$[x,y]\cup [x',y']$, by the convexity of $F_k\cap Hb(\rho )$. Similarly, $y$ and 
$y'$ may be joined by a polygonal line ${\mathbf{L}}_{yy'}$ in $F_k\cap H(\rho )$ of length at most 
$Cd(y,y') $. In the general case $d(x,x')$ and $d(y,y')$ are both of order $\lambda $. In the finite case they are both of order $\ell $.  

We have two possibilities.

(1) There exist $x_1\in \{x,x' \}$ and $y_1\in \{y,y' \}$ contained into two 
non-parallel codimension 1 faces of $H(\rho )\cap F_k$. 

(2)  Every two points $x_1\in \{x,x' \}$ and $y_1\in \{y,y' \}$ are contained in 
the interiors of two parallel codimension 1 faces of $H(\rho )\cap F_k$.

 In the situation (1) according to Lemma \ref{fetze} we may join $x_1$ and $y_1$ by a polygonal line of length at most $Cd(x_1,y_1)$. This and the previous remark on 
the possibility of joining $x,x'$ and $y,y'$ respectively, imply that in this 
way we may join both $x, y$ and $x',y'$ with polygonal lines of lengths 
comparable to the distances in the general case and to $\ell $ in the finite case. 

 In the situation (2) it follows that $x,x'$ are in the interior of the same 
codimension 1 face, and the same for $y,y'$, and the spans of the two faces are 
two parallel hyperplanes $H$ and $H'$, respectively. Let $W$ and $W'$ 
be the Weyl chambers of vertices $x$ and $x'$ respectively and boundary at 
infinity $\Delta_0$. The Weyl chamber $W$ contains the 
segment $(x,y]$ in its interior and so does $W'$ for the segment $(x',y']$. As $f_\rho $ decreases on $[x,y]$ of slope $\beta $ it follows that $\angle_x 
(\overline{\rho_x}, W_x)\leq \angle_x(\overline{\rho_x}, 
\overline{xy})<\frac{\pi }{2}-\delta_1$. In $\Sigma_x {\bf K}$ a geodesic from $\overline{\rho_x} $ to $\overline{xy}$ extends in $(F_k)_x$ to a direction $\overline{x\alpha }$ opposite to $\overline{\rho_x} $, $\alpha \in F_k(\infty )$. Since $\overline{\rho_x} $ is orthogonal to $H_x$, the same is true for $\overline{x\alpha }$. The geodesic in $(F_k)_x$ from $\overline{xy}$ to $\overline{x\alpha }$ contains a point at distance $\frac{\pi }{2}$ from $\overline{\rho_x}$. It follows that $\overline{xy}$ and $\overline{x\alpha }$ are separated by $H_x$, hence $[x,y]$ and $[x,\alpha )$ are separated by $H$. We conclude that a germ of $W$ is contained in the same apartment $\widetilde{F}$ with a germ of $\rho_x$ and a germ of the ray $[x,\alpha )$ in $F_k$ which is orthogonal to $H$ and on the other side of $H$ than $W$. Likewise we obtain that a germ of $W'$ is contained in the same apartment $\widetilde{F}'$ with a germ of $\rho_{x'}$ and a germ of the ray $[x',\alpha ')$ in $F_k$ which is orthogonal to $H$ and on the other side of $H$ than $W'$. Since $W$ and $W'$ are asymptotic, $\alpha' =\alpha $. It follows that the 
relative position of $W$ with respect to $[x,\alpha )$ is the same as the one of 
$W'$ with respect to $[x',\alpha )$. Then the same is true for the relative 
position of $W$ with respect to $\rho_x$ in $\tilde{F}$ and that of $W'$ with 
respect to $\rho_{x'}$ in $\tilde{F}'$. In particular , the chambers $W_x$ and $W_{x'}'$ both contain or do not contain $\overline{\rho_x}$ and $\overline{\rho_{x'}}$, respectively.

Suppose we are in the second case, that is, $W_x$ and $W_{x'}'$ do not contain $\overline{\rho_x}$ and $\overline{\rho_{x'}}$, respectively. If the panel $M$ of $W$ has the property that  $M_x$ separates $\overline{\rho_x} $ and $W_x$, the same is true for the panel $M'$ of $W'$ asymptotic to $M$, with respect to $\overline{\rho_{x'}}$. Such a panel exists in the considered case, so 
we fix one, $M$. Either $M$ or $M'$ has the property that $W\cup W'$ is on 
the same side of its affine span. Suppose it is $M$ and let $\widehat{H}$ be its 
affine span. By Lemma \ref{fetze2}, (2), (a), there exists a ramification $F_k'$ 
of $F_k$ with $\partial (F_k'\cap F_k) =\widehat{H} $, containing $W$ (and 
consequently $W'$) and such that $x$ is in at least two faces of $F_k'\cap 
H(\rho )$. It follows that, by replacing $F_k$ with $F_k'$ the point $x$ may be joined with a polygonal line ${\bf L}_{xy}$ 
 to $y$.

Suppose now that the chambers $W_x$ and $W_x'$ contain $\overline{\rho_x}$ and 
$\overline{\rho_{x'}}$, respectively. Suppose the connected component of $W\cap 
H(\rho )$ containing $y$ has at least two faces. The proof of Lemma \ref{fetze2}, 
$(b_1)$, implies that there exists a point in it, $y_1$, which may be joined to 
$y$ by a polygonal line in $H(\rho )\cap F$ of length $\leq {\bf c}''d(x,y)$ and 
to $x$ by a polygonal line ${\mathbf L}_{xy_1}$. The same is true for $y'$ and 
$x'$ if the connected component of $W'\cap H(\rho )$ containing $y'$ has at 
least two faces.

Finally, if both $W\cap H(\rho )$ and $W'\cap H(\rho )$ have only one face, its 
affine span must be $H'$. We then choose a hyperplane $\widehat{H}$ supporting 
$W$ neither orthogonal nor coincident to $H$ and $\widehat{H}'$ parallel to it 
and supporting $W'$. Either $\widehat{H}$ or $\widehat{H}'$ have $W\cup W'$ on 
the same side. Suppose it is $\widehat{H}$. By Lemma \ref{fetze2}, $(b_2),$ 
there exists a ramification $F_k'$ of $F_k$ containing $W\cup W'$ and with 
$\partial (F_k'\cap F_k)=\widehat{H}$ such that the points in $W\cap 
\widehat{H}\cap H(\rho )$ are in two different faces. Therefore we may join $y$ 
to one of these points, $y_1$, by a segment in $H(\rho )\cap W$ of length $\leq {\bf 
c}''d(x,y)$ and we may join $y_1$ to $x$ in $H(\rho )\cap F_k'$ by a polygonal line ${\mathbf 
L}_{xy_1}$.

\medskip

Now we consider the case of a horoball $Hb(\rho )$ that intersects only $\lbrack P_k ,N_k\rbrack $ in a segment $[x,y]$. In the general case, due to the fact that $d(P_k, P_{k+1}')$ and $d(N_k, N_{k+1})$ are of order $\lambda $, we have, according to \cite[Lemma 4.11]{Dr2*}, that $d(x,y)\leq c\lambda $.

If $x$ and $y$ are contained into nonparallel faces of $H(\rho )\cap F_k$ then by Lemma \ref{fetze} they may be joined by a line ${\bf L}_{xy}$. If not, let $H$ and $H'$ be the affine spans of the faces containing $x$ and $y$ respectively in their interiors. Let $x_1 = H\cap [P_{k+1}', N_{k+1}]$ and $y_1 = H'\cap [P_{k+1}', N_{k+1}]$ (eventually $x_1$ and $y_1$ might be on a prolongation of the segment $[P_{k+1}', N_{k+1}]$). Since $\lbrack P_k ,N_k\rbrack $ and $[P_{k+1}', N_{k+1}]$ both have slope $\beta $ their angle with $H$ and $H'$ is at least $\delta_1$. It follows that $d(x,x_1),\; d(y,y_1)\leq c' d(P_k,P_{k+1}')$. Then the polygonal line $[x, x_1]\cup [x_1,y_1]\cup [y_1,y]$ has length of order $\ell $ in the finite case and $\lambda $ in the general case, and it is outside $Hbo(\rho )$. By the convexity of $F\cap Hb(\rho )$ it follows that this line projects onto a polygonal line in $F\cap H(\rho )$ joining $x$ and $y$ and of the same order of length.

\medskip

 We denote the curve between $P_k$ and $N_k$ thus obtained 
by $\gamma_k$ and the curve between $P_{k+1}'$ and $N_{k+1}$ 
by $\gamma_{k+1}'$. We denote by ${\frak C}_k$ the loop formed
by the curves $\gamma_k$, $\gamma_{k+1}$, ${\mathcal L}_{P_kP_{k+1}}$ 
and $[N_k, N_{k+1}]$. We denote ${\frak C}_k':=\gamma_k \cup 
\gamma_{k+1}'\cup {\mathcal L}_{P_kP_{k+1}'} \cup [N_kN_{k+1}]$ 
and ${\frak C}_{k+1}'':= \gamma_{k+1}'\cup \gamma_{k+1}\cup {\mathcal 
L}_{P_{k+1}P_{k+1}'}$. We note that in the finite case although $P_{k+1}'=P_{k+1}$ the curves $\gamma_{k+1}'$ and $\gamma_{k+1}$ do not coincide : one is a curve contained in $F_k$ or in a ramification of it while the other is contained in $F_{k+1}$ or in a ramification of it. Also, ${\mathcal L}_{P_kP_{k+1}}={\mathcal L}_{P_kP_{k+1}'}=[P_k, P_{k+1}]$ and ${\mathcal 
L}_{P_{k+1}P_{k+1}'}=\{ P_{k+1} \}$. In the general case the shapes of these lines do not really matter. Since the notion of area we work with is discrete all that matters is that the 
sets $\{ P_k,\; P_{k+1}',\; P_{k+1}\}$ have diameters 
of order $\lambda $. Also, in the general case we may replace the loop $\frak C$ with the loop  $\bigcup_{k=0}^{n-1}{\mathcal L}_{P_kP_{k+1}}\cup {\mathcal L}_{P_nP_0}$. Since the two loops are at a Hausdorff distance of order $\lambda $ one from the other the replacement can be done up to adding a linear filling area.

We construct the filling cylinder by filling in all the 
loops ${\frak C}_k$. We do this in two steps : first we fill all the loops ${\frak C}_k'$, 
then all the loops ${\frak C}_{k+1}''$.

\bigskip

\noindent{\sc Step 4.}\quad We fill the loop ${\frak C}_k'$. First, for every pair of segments $[x,y]=[P_k,N_k ]\cap Hb(\rho )$ and $[x',y']=[P_{k+1}',N_{k+1}]\cap Hb(\rho )$, we fill the loops obtained by joining $x$ to $x'$, $y$ to $y'$,$x$ to $y$ and $x'$ to $y'$ as in Step 3. According to the construction in Step 3, the subarc $\gamma_{xy}$ of $\gamma_k $ of extremities $x,y$ and the subarc $\gamma_{x'y'}'$ of $\gamma_{k+1}'$ of extremities $x',y'$ eventually differ at extremities where one is eventually obtained from the other 
by adding either ${\mathbf L}_{xx'}$ or ${\mathbf L}_{yy'}$ or both. Thus the 
loop $\gamma_{xy} \cup {\mathbf L}_{xx'} \cup \gamma_{x'y'}' \cup 
{\mathbf L}_{yy'}$ is reduced to an arc and there is no area needed to fill it.

What remains to be filled is a set of loops of ``vertices'' $y,\; y',\; \bar{x}',\; \bar{x}$, where $y,y'$ are the upper extremities of a pair of segments as previously and $\bar{x},\; \bar{x}'$ are the lower extremities of the next pair.

In the general case the arc $\gamma_{y\bar{x}}$ is at Hausdorff distance of order $\lambda $ from the segment $[y, \bar{x}]$ and the same is true for the arc $\gamma_{y'\bar{x}'}$ with respect to the segment $[y', \bar{x}']$. On the other hand, the segments $[y, \bar{x}]$ and $[y', \bar{x}']$ are at Hausdorff distance of order $\lambda $ one from the other. It follows that the loop $\gamma_{y\bar{x}}\cup {\mathbf L}_{\bar{x}\bar{x}'}\cup \gamma_{y'\bar{x}'}\cup 
{\mathbf L}_{yy'}$ may be filled with an area of order $d(y, \bar{x})$.

In the finite case the loop $\gamma_{y\bar{x}}\cup {\mathbf L}_{\bar{x}\bar{x}'}\cup \gamma_{y'\bar{x}'}\cup 
{\mathbf L}_{yy'}$ has length of order $\ell $ and by Proposition \ref{polit3} it may be filled with an area of order $\ell^2$.

We conclude that to fill the loop ${\frak C}_k'$ we need an area of order 
$\Sigma d(y, \bar{x})$, so of order at most $\ell $ in the general case, and an area of order $\ell^2$ in the finite case (we recall that in this case the family of horoballs is finite).

\bigskip

\noindent{\sc Step 5.}\quad Now we fill the loop ${\frak C}_{k+1}''$.

First we consider the finite case. In this case $P_{k+1}'=P_{k+1}$ and the two arcs $\gamma_{k+1}$ and $\gamma_{k+1}'$ composing the loop differ only between pairs of points $x',\, y'$ such that $[x',y']=[P_{k+1},N_{k+1}]\cap Hb(\rho )$. The arc of $\gamma_{k+1}'$ between $x'$ and $y'$, which we denote by $\gamma_{x'y'}'$, is in $F_k$ or a ramification of it, while the arc of $\gamma_{k+1}$ between $x'$ and $y'$, $\gamma_{x'y'}$, is in $F_{k+1}$ or a ramification of it. Both arcs have length of order $\ell $. The apartments $F_k$ and $F_{k+1}$ or their respective ramifications have in common a Weyl chamber of boundary $\Delta_0$ and vertex $x'$. Proposition \ref{camera} implies that the loop $\gamma_{x'y'}'\cup \gamma_{x'y'}$ can be filled with an area of order $\ell^2$. Since there is a uniformly bounded number of such loops along $[P_{k+1},N_{k+1}]$, composing ${\frak C}_{k+1}''$, we may conclude that ${\frak C}_{k+1}''$ can be filled with an area of order $\ell^2$.

 In the general case $\gamma_{k+1}$ and $\gamma_{k+1}'$ differ between pairs of points $x',\, y'$ such that $[x',y']=[P_{k+1}',N_{k+1}]\cap Hb(\rho )=[P_{k+1},N_{k+1}]\cap Hb(\rho )$ and they may also differ near $P_{k+1}'$ and $P_{k+1}$, respectively. For a pair of points $x',y'$ as previously, we may reason as in the finite case. The only difference is that the lengths of $\gamma_{x'y'}$ and of $\gamma_{x'y'}'$ are of order $d(x',y')$ so the area needed to fill $\gamma_{x'y'}'\cup \gamma_{x'y'}$ is of order $d(x',y')^2$.

We now look at what happens near $P_{k+1}'$ and $P_{k+1}$. We recall that we denoted 
$R_{k+1}$ the first point in which the ray $r_{k+1}$ meets the apartment 
$F_{k}$. We have that $d(P_{k+1},R_{k+1})\leq c\lambda $. When we chose 
$P_{k+1}'$ in Step 2 we
 had three cases :

(1)  If $R_{k+1}$ is 
contained
in
$F_k\cap {\mathbf K}_0$, then $P_{k+1}'=R_{k+1}$ ;

(2) Suppose $R_{k+1}$ is 
contained
in some $Hbo(\rho )\cap F_k$. Then one possibility would be that $P_{k+1}'\in [R_{k+1}, N_{k+1}]$ ;

(3) The other possibility when $R_{k+1}$ is contained in some $Hbo(\rho )\cap F_k$ is that $R_{k+1}\in [P_{k+1}', N_{k+1}]$ ; then $R_{k+1}\in 
[P_{k+1}', y_0']\cap [x_0',y_0']$ where the previous two segments are the intersections of $Hb(\rho )$ with $[P_{k+1}', N_{k+1}]$ and $[P_{k+1}, N_{k+1}]$, respectively.

In the cases (1) and (2) $\gamma_{k+1}$ contains with respect to $\gamma_{k+1}'$ an extra-arc with length of order $\lambda $ which we may ignore.

In the case (3), $\gamma_{k+1}'$ and $\gamma_{k+1}$ differ between $ 
P_{k+1}', y_0'$ and $x_0',y_0'$, respectively, and $\gamma_{k+1}$ again contains the 
extra-arc between $P_{k+1}$ and $x_0'$ which we may likewise ignore, as being of 
order $\lambda $. Let $\gamma_{P_{k+1}', y_0'}$ be the curve which replaces the segment 
$ [P_{k+1}', y_0']$ in $\gamma_{k+1}'$ and $\gamma_{x_0',y_0'}$ the curve which replaces 
the segment $ [x_0', y_0']$ in $\gamma_{k+1}$. The difference with the situation of 
the curves $\gamma_{x',y'}$ and $\gamma_{x',y'}'$ studied before is that the extremities of these two 
curves do not coincide. But a simple argument allows to apply however Proposition
\ref{camera}. The argument is as follows.
 Besides the horosphere $H(\rho )$ we also consider the horosphere 
$H_{-c_{0}\lambda }(\rho )$ containing $R_{k+1}$, where $0<c_0\leq c$ by Step 2. Let 
$y_0''$ be the second intersection point of $[ R_{k+1}, N_{k+1}]$ with 
$H_{-c_{0}\lambda }(\rho )$.

By Corollary \ref{prf} the curves $\gamma_{P_{k+1}', y_0'}$ and $\gamma_{x_0',y_0'}$ project onto two curves ${\bf {\frak c}}_1$ and ${\bf {\frak c}}_2$ 
of smaller length, contained in $H_{-c_{0}\lambda }(\rho )\cap F_k$ and  $H_{-c_{0}\lambda }(\rho )\cap F_{k+1}$ respectively, at Hausdorff distance $C\cdot c_0\lambda $ from the 
initial curves. By eventually prolongating ${\bf {\frak c}}_i,\; i=1,2$, with arcs of 
length of order $\lambda $ one may suppose that both have $R_{k+1}$ and $y_0''$ as 
extremities. By Proposition
\ref{camera} one needs a $\lambda$-filling area of order $d(R_{k+1},y_0'')^2$ to fill ${\bf 
{\frak c}}_1\cup {\bf {\frak c}}_2 $. By the properties of the horospheres ennounced in Section 
\ref{h}, the filling area previously found for ${\bf {\frak c}}_1\cup {\bf {\frak c}}_2 $ gives 
a $\delta $-filling area for $\gamma_{P_{k+1}', y_0'}\cup \gamma_{x_0',y_0'}\cup {\mathcal 
L}_{x_0',P_{k+1}'}$, with $\delta =\lambda (1+6C\cdot c_0)$. For $\lambda $ 
sufficiently small $\delta$ is smaller than $1$.

In the end we obtain that in the general case the loop ${\frak C}_{k+1}''$ can be filled with an area of order $\Sigma_{x',y'}d(x',y')^2+d(P_{k+1}', y_0')^2 \leq \ell^2 $.

\bigskip

\noindent{\sc Step 6.}\quad We conclude that to fill the loops ${\frak C}_k'$ 
and 
${\frak C}_{k+1}''$, so to fill the loop ${\frak C}_k$, we need an area of order 
$\ell^2$. By 
summing over $\{1,2,\dots n\}$ we obtain a filling area of order $n\ell^2$ for the initial loop. In the finite case this gives an area of order $\ell^2 $ while in the general case this gives an area of order 
$\ell^3$.

\hspace*{\fill}$\diamondsuit $

\begin{remark}\label{rc}
We note that the cubic order in the general case comes from the fact that to fill the loops ${\frak C}_k''$ one must sometimes spend an area of order
$\ell^2$. But the important thing is
that all the bricks of length at most $1$ used to fill each ${\frak C}_k''$ are more or less boundaries of small Euclidean triangles entirely
contained in a polytopic surface $H(\rho )\cap F$.

Thus for a generic loop of length
$\ell$ we have obtained a $1$-partition composed of $k_1\ell^3$ boundaries of small Euclidean triangles entirely
contained in polytopic surfaces of type $H(\rho )\cap F$
and of $k_2\ell^2$ bricks on which nothing special can be said. 
\end{remark}

\subsection{Quadratic filling order in solvable groups}\label{qsol}

 By means of Theorem \ref{fillK0}, (b), we prove Theorem \ref{fillgen}, (b). A corollary of it is the fact that the  filling order in solvable groups acting on horospheres in symmetric spaces of rank at least $3$ is quadratic.

First we prove an intermediate result. 

\begin{proposition}\label{ap}
Let $X$ be a product of symmetric spaces
and Euclidean buildings, $X$ of rank at least $3$, and $X_0$ a subset of it which can 
be written as 
$$
X_0=X\setminus \bigsqcup_{\rho \in {\mathcal R}}Hbo(\rho )\; .
$$

Suppose the set of rays ${\mathcal R}$ is finite and suppose $X_0$ has the properties ({\cal P}$_1$) and ({\cal P}$_2$) formulated in Theorem \ref{fillgen}. For every $m\geq 4$ there exists a constant $\frak K$ depending on $m$, on $X$, on the cardinal of ${\mathcal R}$, on the constant $d$ appearing in property ({\cal P}$_1$) and on the slope $\theta $ appearing in property ({\cal P}$_2$) such that for every loop $\frak C$ composed of at most $m$ minimizing almost polygonal curves, 
 $$
A_1({\frak C})\leq {\frak K}\ell^{2} \, . \leqno(4.2) 
$$
\end{proposition}

\noindent{\bf Proof.}\quad It suffices to prove that (4.2) is satisfied for loops ${\frak C}$ with length at least $\ell_0$ for $\ell_0$ sufficiently large. We reason by contradiction. Suppose that in $X$ there exists a sequence of subsets $X_0^n=X\setminus \bigsqcup_{\rho \in {\mathcal R}_n}Hbo(\rho )$ with card ${\mathcal R}_n\leq N$, and $\theta $ the common slope of all rays in ${\mathcal R}_n$, such that $X_0^n$ has properties ({\cal P}$_1$) and ({\cal P}$_2$), where the constant $d$ in ({\cal P}$_1$) is independent of $n$, and in each $X_0^n$ there is a loop ${\frak C}_n$ of length $\ell_n \geq \ell_0$ composed of at most $m$ minimizing almost polygonal curves, with 
$$
A_1({\frak C}_n)\geq n\ell_n^{2} \, .\leqno(4.3)
$$
In each $X_0^n$ we consider a loop ${\frak C}_n$ of minimal length with the previous properties. Inequality $(4.3)$ implies that $\ell_n$ must diverge to $+\infty $. 
 Let $x_n$ be a point on the image of 
${\frak C}_n$ and let ${\bf K}=X_\omega (x_n, \frac{\ell_n }{10})$ and ${\bf K}_0=[ X_0^n ] $. We may write ${\bf K}_0={\bf K}\setminus \bigsqcup_{\rho_\omega \in {\mathcal R}_\omega }Hbo(\rho_\omega )$, where all rays $\rho_\omega $ have slope $\theta $. Since card ${\mathcal R}_n\leq N,\; \forall n\in \N $, and since $\omega $ chooses one out of a finite number of possibilities , card ${\mathcal R}_\omega \leq N$. Also ${\bf K}_0$ has properties ({\cal P}$_1$) and ({\cal P}$_2$) formulated in Theorem \ref{fillK0}. The limit set of the sequence of loops ${\frak C}_n$ is a loop ${\frak C}$ of length $10 $ composed of at most $3q_0mN$ segments. According to Theorem \ref{fillK0}, (b), $\frak C$ may be filled with an area of at most $100C$. Moreover each of the $100C$ bricks composing the filling disk is an Euclidean triangle contained either in the exterior of a set of polytopes in a maximal flat, $F\setminus \bigsqcup_{\rho_\omega\in {\mathcal R}_\omega} Hbo(\rho_\omega )$, or in a face of a polytope $F'\cap H(\rho_\omega ),\; \rho_\omega\in {\mathcal R}_\omega $, so in a hyperplane $F\cap H(\rho_\omega )$, where $F$ asymptotic to $\rho_\omega$. Thus each of the bricks $B_i,\; i\in \{1,2,\dots  ,100C\}$, is the limit of a sequence of  Euclidean triangles $B_i^n$ contained either in flat sets of type $F_i^n \setminus \bigsqcup_{\rho \in {\mathcal R}_n } Hbo(\rho )$ or in intersections $F_i^n\cap H(\rho_i^n)$ with $\rho_i^n \in {\mathcal R}_n$ and $F_i^n$ asymptotic to $\rho_i^n$.

Let $E_i$ and $E_j$ be the edges of two bricks $B_i$ and $B_j$ that coincide, and have length $\Lambda \leq 1$. Let $E_i^n$ and $E_j^n$ be the sequences of edges of $B_i^n$ and $B_j^n$, respectively, so that $[E_i^n]=E_i=E_j=[E_j^n]$. The Hausdorff distance $\delta_{ij}^n$ between $E_i^n$ and $E_j^n$ has the property that $\lim_\omega\frac{\delta_{ij}^n}{\ell_n}=0$. By joining the pairs of extremities of $E_i^n$ and $E_j^n$ which give the same limit point with minimizing almost polygonal curves we obtain a loop, ${\frak C}_{ij}^n$. The pair of ``opposite sides'' $E_i^n$ and $E_j^n$ of ${\frak C}_{ij}^n$ have lengths of order $\ell_n$ while the other pair of ``opposite sides'' have lengths of order $\delta_{ij}^n=o(\ell_n)$. We may divide the loop ${\frak C}_{ij}^n$ into approximately $\frac{\Lambda \ell_n}{\delta_{ij}^n }$ loops composed of $2$ minimizing almost polygonal curves and of two subsegments of $E_i^n$ and $E_j^n$ respectively, and with lengths of order $\delta_{ij}^n$. Since $\ell_n$ was the minimal length of a loop satisfying $(4.3)$, each of these loops has filling area at most $Kn \left(\delta_{ij}^n \right)^2$, where $K$ is an universal constant. It follows that $A({\frak C}_{ij}^n)\leq \kappa' n\ell_n \delta_{ij}^n$.

We may fill ${\frak C}_n$ by filling each of the bricks $B_i^n\subset F_i^n \setminus \bigsqcup_{\rho \in {\mathcal R}_n } Hbo(\rho )$ and $B_i^n\subset F_i^n\cap H(\rho_i^n)$ and each of the loops ${\frak C}_{ij}^n$. Thus $A_1({\frak C}_n)\leq 100C \cdot \kappa \frac{\ell_n^2}{100} + 200C\kappa' n\ell_n \delta_{ij}^n$. It follows that $n\ell_n^{2}\leq C \cdot \kappa \ell_n^2 + 200C\kappa' n\ell_n \delta_{ij}^n$. If we divide the inequality by $n\ell_n^{2}$ and we consider the $\omega $-limit, we obtain $1\leq 0$ so a contradiction. \hspace*{\fill } $\diamondsuit$

\bigskip

\noindent{\bf Proof of Theorem \ref{fillgen}, (b).}\quad We shall proceed by induction. First we need some constants. According to Theorem \ref{nondist2} $X_0$ is undistorted in $X$. Let ${\bf C}$ be a nondistorsion constant. 
Then a minimizing almost polygonal curve joining two points $x,y$ in $X_0$ has length at most $2C_0d(x,y)+6q_0N{\bf C}\epsilon(x,y)$, where $N=$ card ${\mathcal R}$. It follows that its length is at most ${\frak k}d(x,y)$ with ${\frak k}= 2C_0+24q_0N{\bf C}C_0$. Let $b$ be an integer which is very large compared to ${\frak k}$. Then there exists an integer $M$ between $b({\frak k}+1)$ and $\frac{b^2}{2}$. Let ${\frak K}$ be the constant provided by Proposition \ref{ap} for loops composed of at most $M$ minimizing almost polygonal curves. Let also $C\geq 2{\frak K}{\frak k}^2$. We show by induction the following statement

\begin{center}

(I$_n$) If $\ell \leq b^n$ then $A_1(\ell )\leq C\ell^2$.

\end{center}       

(I$_0$) is satisfied if $C$ is big enough. Suppose (I$_n$) is satisfied and let us prove (I$_{n+1}$). Let $\frak C$ be a loop of length $\ell \in (b^n, b^{n+1}]$. We divide the loop $\frak C$ into $M$ arcs and we join the extremities of each of these arcs by almost polygonal curves. We thus obtain M loops ${\frak c}_1,\, {\frak c}_2,\dots ,{\frak c}_M$, of lengths $\frac{\ell }{M}(1+{\frak k})$ and one loop ${\frak C}_0$ of length at most ${\frak k}\ell$ composed of $M$ minimizing almost polygonal curves.

Since $M\geq b({\frak k}+1)$, $\frac{1+{\frak k}}{M}\leq \frac{1}{b}$ and, by (I$_n$), $A({\frak c}_i)\leq C\frac{\ell^2}{b^2},\; \forall i$. By Proposition \ref{ap} $A({\frak C}_0)\leq {\frak K}{\frak k}^2\ell^2$. It follows that $A({\frak C})\leq CM\frac{\ell^2}{b^2}+{\frak K}{\frak k}^2\ell^2\leq C\ell^2$.\hspace*{\fill }$\diamondsuit $

\medskip

\begin{remark}
To obtain Theorem \ref{T1} in the introduction from Theorem \ref{fillgen}, (b), we only need to verify that if the space $X_0$ is the exterior of an open 
horoball, then it verifies property $({\mathcal P}_1)$. This is proved in \cite[proof of Corollary 4.16]{Dr2*}.
\end{remark}

\subsection{Asymptotically quadratic filling order in lattices}\label{aq}

We show that a space $X_0$ on which a $\Q$-rank one lattice acts cocompactly, endowed with the 
induced
metric, has asymptotically quadratic filling order. We note that the space $X_0$ satisfies properties $({\mathcal P}_1)$ and 
$({\mathcal P}_2)$ formulated in Theorem \ref{fillgen} (see \cite[Propositions 5.5 and 5.7]{Dr1*} and \cite[Lemma 8.3]{Mo*}). Hence it suffices to prove Theorem \ref{fillgen}, (a).

\bigskip

\noindent{\bf Proof of Theorem \ref{fillgen}, (a).}\quad By Theorem \ref{fillK0}, (a), the filling order in 
every 
asymptotic cone ${\mathbf
K}_0$ of
$X_0$ is at most cubic, that is $A_1(\ell )\leq k\ell^3$. We shall
show  that it is 
actually quadratic. We have similarities between
 asymptotic cones (Remark \ref{si}) which 
allow to 
conclude, first, that $A_\lambda(\ell )\leq
k\left( \frac{\ell}{\lambda } \right)^3,\; \forall \lambda >0$. Now we recall that $A_1(\ell )\leq k\ell^3 $ means that we
need at  most 
$k\ell^3$ bricks of length at most $1$ to fill a loop of
length 
$\ell $. 
But by construction $k_1\ell^3 $ of these bricks bound small Euclidean triangles entirely
contained into intersections  of 
horospheres with apartments, while $k_2\ell^2$ bricks have shapes on which we
know  nothing. If a 
small loop ${\frak c}_1$ of length less than 1 bounds a small Euclidean triangle entirely contained into
the  intersection 
of a horosphere with an apartment, then its filling area in ${\mathbf 
K}_0$ is 
$A_\lambda({\frak c}_1)\leq k\frac{1}{\lambda^2},\; \forall \lambda \in [0,1[$. If a loop ${\frak c}_2$
has length  at most 
one and is arbitrary, at least $A_\lambda({\frak c}_2)\leq 
k\frac{1}{\lambda^3}$. Thus, 
for a generic loop ${\frak C}$ of length at most $\ell$ we have 
$$
A_\lambda ({\frak C})\leq k_1\ell^3 \cdot k\frac{1}{\lambda^2 }+ 
k_2\ell^2\cdot k\frac{1}{\lambda^3}\; .\leqno(4.4)
$$

If we replace $\lambda$ by $\frac{\ell}{M}$, we obtain 
$$
P\left( \ell,\frac{\ell}{M}\right) \leq k_1'\ell
M^2+\frac{k_2'}{\ell}M^3\; .
$$

If one takes $\sqrt{M}\leq \ell \leq 2\sqrt{M}$, one obtains 
$$
P\left(\ell,\frac{\ell}{M}\right)\leq CM^{2.5}\; .\leqno(4.5)
$$

But since, by similarities and changing cone (Remark \ref{si}), one can modify the 
length, the 
relation $(4.5)$ holds for every length $\ell $ and every $M$ in every asymptotic cone. We 
may conclude that the filling order is at most $2.5$ in all asymptotic cones.  

In the same way one can show that the filling order is
quadratic. Suppose  that the 
minimal order of filling common to all asymptotic cones is 
$2+\varepsilon $. 
Then in all asymptotic cones we have $A_1(\ell )\leq 
k\ell^{2+\varepsilon }$. By 
means of similarities we may conclude that $A_\lambda(\ell )\leq 
k\left( \frac{\ell}{\lambda }\right)^{2+\varepsilon }$. Then we
can modify the  relation 
$(4.4)$ and write 
$$
A_\lambda ({\frak C})\leq k_1\ell^3 \cdot k\frac{1}{\lambda^2 }+ 
k_2\ell^2\cdot k\frac{1}{\lambda^{2+\varepsilon}}\; ,
$$
 which implies that in all asymptotic cones we have
$$
P\left( \ell,\frac{\ell}{M}\right) \leq k_1'\ell 
M^2+\frac{k_2'}{\ell^\varepsilon}M^{2+\varepsilon}\; .
$$

For lengths $\ell \in \lbrack 
M^{\frac{\varepsilon}{2}},2M^{\frac{\varepsilon}{2}}\rbrack $, we 
obtain

$$
P\left(\ell,\frac{\ell}{M}\right)\leq 
CM^{2+\varepsilon-\frac{\varepsilon^2}{2}}\; .
$$

The previous inequality can be generalized to all lengths, by similarities and changing cone. It finally gives a filling order $2+\varepsilon-\frac{\varepsilon^2}{2}$ in
all asymptotic cones. This  contradicts the minimality of the
order
$2+\varepsilon
$.

Thus we get a quadratic filling order in all asymptotic cones of $X_0$, so 
an 
asymptotically quadratic filling order in $X_0$, by Theorem \ref{crit}.
\hspace*{\fill}$\diamondsuit $

\section{Appendix}

We prove several useful results of Euclidean geometry.

\begin{lintro}\label{simplu}
Let ${\mathcal P}$ be a convex polytope in the Euclidean space $\E^n ,\, n\geq 2,$ and let 
$\Phi \subsetneq \E^n$ be its affine span. Let $k$ be the common dimension of ${\mathcal P}$ and $\Phi ,\; k\leq n-1$. For every $p\in \Phi $ we denote 
$\Phi_{p}^{\perp}$ the subspace orthogonal to $\Phi $ through $p$.

Suppose the polytope ${\mathcal P}$ is either of codimension at least $2$ or of 
codimension 1 and different from $\Phi $. Let the hypersurface $\partial 
{\mathcal{N}}_R({\mathcal P}),\; R>0$, be endowed with the length metric $d_\ell$. For two arbitrary points $ x,y $ in $\partial 
{\mathcal{N}}_R({\mathcal P})$ we denote by $x_0,y_0$ their respective projections on 
${\mathcal P}$ and by $x',y'$ their respective projections on $\Phi $. Let $\alpha_x 
:=\angle_{x_0}(x,\Phi_{x_0}^{\perp}),$ $ \alpha_y 
:=\angle_{y_0}(y,\Phi_{y_0}^{\perp})$ and $\beta_{xy}$ the angle between the 
projection of $[x_0,x]$ on $\Phi_{x_0}^{\perp}$ and the translation of vector 
$\overrightarrow{y_0x_0}$ of the projection of $[y_0,y]$ on 
$\Phi_{y_0}^{\perp}$. If one of the two previous projections is a point then we 
take $\beta_{xy}=0$. Let $z\in \partial {\mathcal P}$ be such that 
$d(x_0,z)+d(z,y_0)=\inf_{t\in \partial {\mathcal P}}[ d(x_0,t)+d(t,y_0)]$.

Suppose that the segment $[x',y']$ intersects ${\mathcal P}$. Then there exist constants $k_1, k_2, \varkappa_1, \varkappa_2$ independent of 
$x,\; y,$ ${\mathcal P}$ and $R$ such that

\begin{itemize}
\item[(a)]  if $\Phi $ is either of codimension at least 2 or of codimension 1 
and not separating $x$ and $y$ then 
$$
k_1d(x,y)\leq d_\ell (x,y)\leq k_2d(x,y)\leqno(5.1)
$$
 and 
$$
\varkappa_1(d(x_0,y_0)+R(\alpha_x +\alpha_y +\beta_{xy} ))\leq d_\ell (x,y)\leq 
\varkappa_2(d(x_0,y_0)+R(\alpha_x +\alpha_y +\beta_{xy} ))\; ;\leqno(5.2)
$$

\item[(b)]  if $\Phi $ is of codimension 1 and separating $x$ and $y$ then 
$$
k_1(d(x,y)+d(x_0,z)+d(z,y_0))\leq d_\ell (x,y)\leq 
k_2(d(x,y)+d(x_0,z)+d(z,y_0))\leqno(5.3)
$$
 and 
$$
\varkappa_1(d(x_0,z)+d(z,y_0)+R(\alpha_x +\alpha_y +\pi ))\leq d_\ell (x,y)\leq 
\varkappa_2(d(x_0,z)+d(z,y_0)+R(\alpha_x +\alpha_y +\pi ))\; .\leqno(5.4)
$$
\end{itemize}
\end{lintro}

\noindent{\bf{Proof.}}\quad Since ${\mathcal P}$ is a convex 
set, $d(x,y)\geq d(x_0,y_0)$.

We have several cases.

\medskip

{\bf{(A)}}.\quad Suppose that the segments $[x,x_0]$ and $[y,y_0]$ are 
orthogonal to $\Phi $. Then $[x,x_0]\subset \Phi_{x_0}^{\perp}$ and 
$[y,y_0]\subset \Phi_{y_0}^{\perp}$. Let $y_1$ be the image of $y$ by the 
translation of vector $\overrightarrow{y_0x_0}$. Then $y_1\in 
\Phi_{x_0}^{\perp}\cap \partial {\mathcal{N}}_R({\mathcal P})$. Also, $y_1$ is the 
projection of $y$ on $\Phi_{x_0}^{\perp}$. This implies that $d(x,y_1)\leq 
d(x,y)$.

Suppose $\Phi $ is of codimension at least $2$. Since $x$ and $y_1$ are both points of $\Phi_{x_0}^{\perp}\cap 
\partial {\mathcal{N}}_R({\mathcal P})$, which is an Euclidean sphere of dimension at least 
$2$, they can be joined by an arc on this sphere, of length comparable to 
$d(x,y_1)$. We consider the curve ${\frak C}_{xy}$ composed of this arc and the segment $[y_1,y]$. Its length is $d(x_0,y_0)+R\beta_{xy}$ and it is contained in $\partial {\mathcal{N}}_R({\mathcal P})$. Thus we obtain the second inequality in $(5.2)$. The first inequality as well as $(5.1)$ are easy to obtain.

\smallskip

Suppose $\Phi $ is of codimension $1$ and not separating $x$ and $y$. Then $y_1=x$ and all inequalities are obvious.

\smallskip

Suppose $\Phi $ is of codimension $1$ and it separates $x$ and $y$. Let $z$ be as in the statement of the Lemma and $x_2$ and $y_2$ the respective images of $x$ and $y$ by the 
translations of vectors $\overrightarrow{x_0z}$ and $\overrightarrow{y_0z}$. Then $x$ and $y$ can be joined in $\partial {\mathcal{N}}_R({\mathcal P})$ by joining $x$ to $x_2$ and $y$ to $y_2$ respectively by segments, and $x_2$ to $y_2$ by a half circle of length $\pi R$. Thus the right-hand inequalities in $(5.3)$ and $(5.4)$ are satisfied. In order to prove the left-hand inequalities it suffices to prove that $d_\ell (x,y)\geq d(x_0,z)+d(z,y_0)$. Suppose we choose a system of coordinates such that $Ox_1\dots x_{n-1}=\Phi $ and $Ox_n$ is orthogonal to $\Phi $. Every curve joining $x$ and $y$ in $\partial {\mathcal{N}}_R({\mathcal P})$ has at least one point $\zeta $ with the $n$-th coordinate zero. Then $d_\ell (x,y)\geq d(x,\zeta )+d(\zeta ,y)\geq d(x_0,\zeta )+d(\zeta ,y_0) \geq d(x_0,z)+d(z,y_0)$.  

\medskip

\begin{figure}[!ht]
\centering
\includegraphics{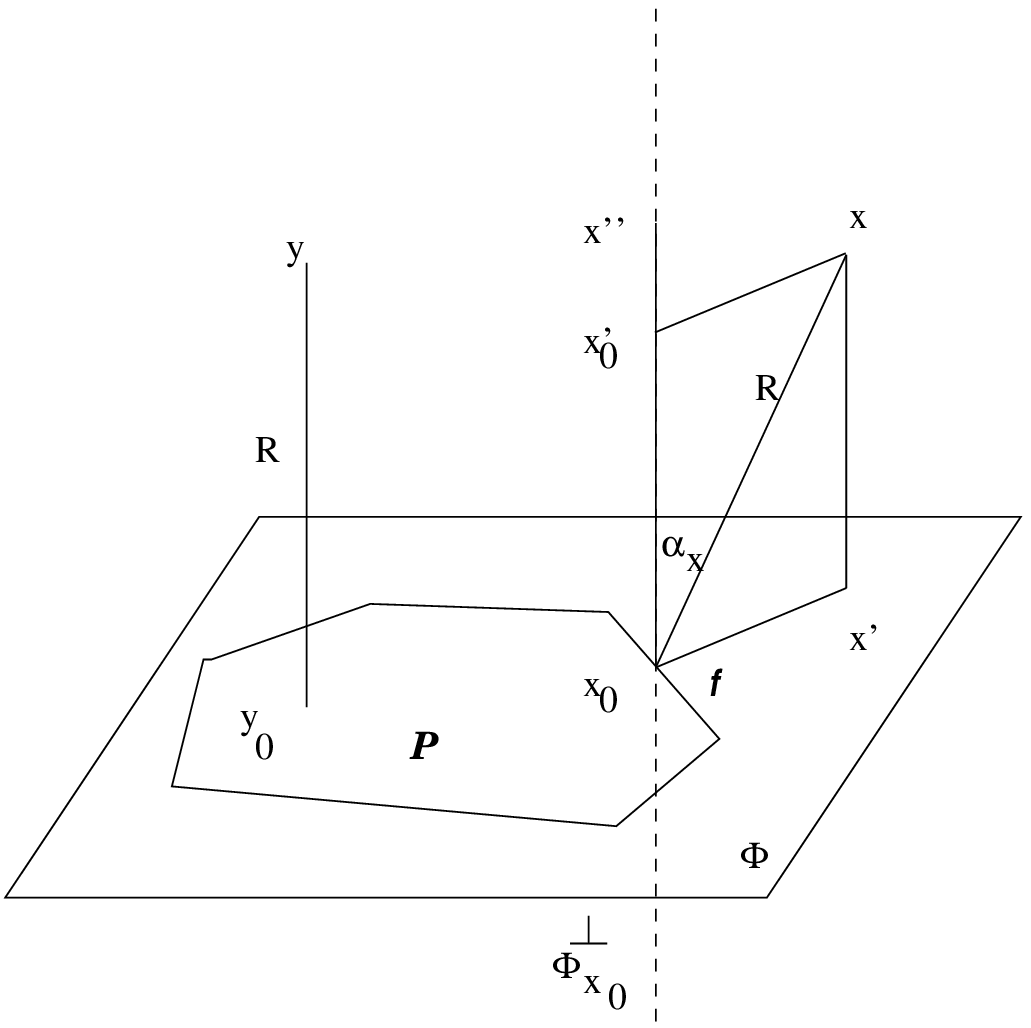} 
\caption{}
\label{fig4} 
\end{figure}

{\bf{(B)}}.\quad Suppose that $[y,y_0]$ is orthogonal to $\Phi $ but $[x,x_0]$ 
is not. Then $x_0$ is contained in a face $\mathfrak f$ of ${\mathcal P}$. Since $x_0$ is also the projection of 
$x'$ on ${\mathcal P}$, $d(x',x_0)\leq d(x',y_0) \leq d(x,y)$. Let $x_0'$ be the projection of $x$ on 
$\Phi_{x_0}^{\perp }$ (Figure \ref{fig4}).

We have the equalities $d(x', x_0)=d(x, x_0')=R\sin \alpha_x$. Let now $x''$ be the intersection point of the ray of origin $x_0$ through $x_0'$ with $\partial 
{\mathcal{N}}_R({\mathcal P})$. If $x_0=x_0'$ then we chose as $x''$ the image of $y$ by the 
translation of vector $\overrightarrow{y_0x_0}$. The arc in $\partial {\mathcal{N}}_R({\mathcal P})$ joining the points 
$x$ and $x''$, which is actually an arc on the sphere $S(x_0, R)$, is of length 
$R\alpha_x = d(x', x_0)\frac{\alpha _x}{\sin 
\alpha_x }\leq cd(x',x_0)\leq cd(x,y)$. So to end the argument it suffices now to 
join the points $x''$ and $y$ by an arc in $\partial {\mathcal{N}}_R({\mathcal P})$ with 
length of order $O(d(x'',y))$. This can be done as in the previous case 
{\bf{(A)}}.

\medskip

{\bf{(C)}}.\quad Suppose that neither $[x,x_0]$ nor $[y,y_0]$ are orthogonal to 
$\Phi $. Let $x_0'$ and $x''$ be obtained as previously from $x$. If $x_0=x_0'$ then we consider $y$ instead of $x$ and its associated points $y_0'$ and $y''$. If $y_0=y_0'$ also, then we return to $x$ and we choose as $x''$ any point in $\partial 
{\mathcal{N}}_R({\mathcal P})\cap \Phi_{x_0}^{\perp }$. In the sequel we suppose, for the sake of simplicity, that we are in the point $x$ and work with $x_0'$ and $x''$ which are chosen one way or the other.

Since $[x',y']$ intersects ${\mathcal P}$, $d(x',y')\geq d(x_0,x')$, so $d(x,y)\geq d(x_0,x')$. On the other hand $d(x_0,x')=R \sin \alpha_x$. As in case 
{\bf{(B)}}, we can join $x$ and $x''$ with an arc in $\partial 
{\mathcal{N}}_R({\mathcal P})$ of length $R\alpha_x\leq cd(x',x_0)\leq cd(x,y)$. It suffices then to join $x''$ and $y$ in $\partial {\mathcal{N}}_R({\mathcal P})$ with an arc of length of order $O(d(x'',y))$. This can be done, by case 
{\bf{(B)}}.\hspace*{\fill }$\diamondsuit $

\medskip

\begin{lintro}\label{proj}
Let ${\mathcal P}$ be a convex polytope in the Euclidean space $\E^n,\, n\geq 2,$ and let 
$\Phi \subsetneq \E^n$ be its affine span. Suppose the polytope ${\mathcal P}$ is either of codimension at least $2$ or of 
codimension 1 and different from $\Phi $. Let the hypersurface $\partial 
{\mathcal{N}}_R({\mathcal P}),\; R>0$, be endowed with the length metric $d_\ell$ and let $a\geq 1$ be a fixed constant.
\begin{itemize}
\item[(a)] The projection from $\partial {\mathcal{N}}_{aR}({\mathcal P})$ onto $\partial  
{\mathcal{N}}_R({\mathcal P})$ is bilipschitz with respect to the length metrics, the constant of the bilipschitz equivalence depending only on $a$.
\item[(b)] Let $\Re $ be a convex polytope in 
$\E^n$. If 
$$
{\mathcal{N}}_R({\mathcal{P}})\subset \Re \subset {\mathcal{N}}_{aR}({\mathcal{P}})\; 
$$
then the projection of the hypersurface $\partial \Re$ onto $\partial 
{\mathcal{N}}_R({\mathcal{P}})$ is bilipschitz with respect to the length metrics, the 
bilipschitz constant depending only on $a$.
\end{itemize}
\end{lintro}

\noindent{\bf Proof.}\quad First we show that two arbitrary points in $\partial 
{\mathcal{N}}_R({\mathcal P})$ may be joined by a curve in $\partial 
{\mathcal{N}}_R({\mathcal P})$ of a special kind and of length comparable to the length distance. Let $x,y$ be two such points and let $\frak g$ be a geodesic joining $x$ and $y$ in $\partial {\mathcal{N}}_R({\mathcal P})$. The projection of ${\frak g}$ on ${\mathcal P}$ is covered by a finite set of faces ${\frak f}_1,\; {\frak f}_2, \;  \dots ,\, {\frak f}_k$ of ${\mathcal P}$ (where ${\mathcal P}$ itself is considered a face). We fix the face ${\frak f}_i$. Let $p$ be the first point of ${\frak g}$ with projection $p_0$ on ${\mathcal P}$ contained in ${\frak f}_i$ and $[p,p_0]$ orthogonal to ${\frak f}_i$. If such a point $p$ does not exist we may ignore the face ${\frak f}_i$. Let $q$ be the last point of ${\frak g}$ with the same properties as $p$ and let $q_0$ be its projection on ${\mathcal P}$. If $p\equiv q$ we may likewise ignore the face ${\frak f}_i$. There are two cases : either ${\frak f}_i$ is of codimension at least $2$ in $\E^n$ or ${\frak f}_i\equiv {\mathcal P}$ is of codimension one.

Suppose ${\frak f}_i$ is of codimension at least $2$ in $\E^n$. Let $p'$ be the image of $p$ by the translation of vector $\overrightarrow{p_0q_0}$. Then $p$ and $q$ may also be joined in $\partial {\mathcal{N}}_R({\mathcal P})$ by the curve composed of the segment $[p,p']$ and of the arc of sphere between $p'$ and $q$ contained in $S(q_0,R)\cap (Span\; {\frak f}_i)_{q_0}^{\perp }$. The length of this curve is $\leq cd(p,q)$, with $c$ a universal constant, so $\leq cd_\ell (p,q)$.

Suppose ${\frak f}_i\equiv {\mathcal P}$ is of codimension one. In this case according to Lemma \ref{simplu} and its proof, $p$ and $q$ may be joined in $\partial {\mathcal{N}}_{R}({\mathcal P})$ either by the segment $[p,q]$ or by the union of two segments with an arc of length $\pi R$. In both cases the curve previously described, joining $p$ and $q$ in $\partial {\mathcal{N}}_{R}({\mathcal P})$, has length at most $\frac{1}{\varkappa_1} d_\ell (x,y)$.

We repeat the previous argument for each face ${\frak f}_i$. We obtain in the end a curve $\frak C$ joining $x$ and $y$, with length at most $c'd_\ell (x,y)$, which decomposes as ${\frak C}={\frak C}_1\cup {\frak C}_2 \cup \cdots \cup {\frak C}_k$, where each ${\frak C}_i$ corresponds to a face ${\frak f}_i$.

\medskip

(a)\quad Let $x,y$ be two points in $\partial {\mathcal{N}}_{aR}({\mathcal P})$ and let $x',y'$ be their respective projections on $\partial {\mathcal{N}}_{R}({\mathcal P})$. Since ${\mathcal{N}}_{R}({\mathcal P})$ is a convex set, $d_\ell(x,y)\geq d_\ell (x',y')$. In order to prove the converse inequality, we join $x'$ and $y'$ in $\partial {\mathcal{N}}_{R}({\mathcal P})$ with a curve ${\frak C}={\frak C}_1\cup {\frak C}_2 \cup \cdots \cup {\frak C}_k$ as previously. For every curve ${\frak C}_i$ we construct the curve ${\frak C}_i'$ by considering for each point $x_i\in {\frak C}_i$ its projection $x_i^0$ on ${\mathcal P}$ and the intersection of the segment $[x_i^0,x_i]$ extended behind $x_i$ with $\partial {\mathcal{N}}_{aR}({\mathcal P})$. It is obvious that the length of ${\frak C}_i'$ is less than the length of ${\frak C}_i$ multiplied by $a$. Then $d_\ell (x,y)\leq \sum_{i=1}^k {\rm{length }}({\frak C}_i')\leq a\sum_{i=1}^k {\rm{length }}({\frak C}_i)=a\cdot {\rm{ length }}({\frak C})\leq a c' d_\ell(x',y')$.

\medskip

(b)\quad Let $x$ and $y$ be two distinct points on $\partial 
\Re$ and let $x'$ and $y'$ be their respective projections on $\partial 
{\mathcal{N}}_R({\mathcal{P}})$. Obviously, $d_\ell (x,y)\geq d_\ell (x',y')$. For an inequality in the other 
way we first note that $x'$ is on the segment $[x,x_0]$, where $x_0$ is the projection 
of $x$ on ${\mathcal{P}}$. Let $x''$ be the intersection point of 
$\partial {\mathcal{N}}_{aR}({\mathcal{P}})$ with the ray through $x$ of origin $x_0$. The 
point $y''$ is 
obtained from $y$ in the same way. By (a) $d_\ell (x'',y'')\leq ac'd_\ell (x',y')$. Thus, it suffices to show that $d_\ell (x,y)\leq kd_\ell (x'',y'')$. In the proof of the inequality the following remark is essential. By the convexity of $\Re $, every hyperplane $H$ in $\E^n$ containing a face of $\partial \Re $ has the property that the distance from ${\mathcal P}$ to $H$ is at least $R$. In particular for every point $\alpha_0 \in {\mathcal P}$ and every nontrivial segment $[\alpha ,\beta ]\subset \partial \Re $, the distance from $\alpha_0$ to the line $\alpha \beta$ is at least $R$.

We join $x''$ and $y''$ in $\partial {\mathcal{N}}_{aR}({\mathcal{P}})$ as in the beginning of the proof, with a curve ${\frak C}={\frak C}_1\cup {\frak C}_2 \cup \cdots \cup {\frak C}_k$ of length at most $c'd_\ell(x'',y'')$. Each ${\frak C}_i$ is joining two points $p_i$ and $q_i$, whose projections on ${\mathcal{P}},\, p_i^0$ and $q_i^0$, are contained in a face ${\frak f}_i$ of ${\mathcal{P}}$, and such that $[p_i,p_i^0]$ and $[q_i,q_i^0]$ are orthogonal to ${\frak f}_i$. We note that $q_i\equiv p_{i+1},\; \forall i$. Let $\bar{p}_i = [p_i,p_i^0]\cap \partial \Re$ and $\bar{q}_i = [q_i,q_i^0]\cap \partial \Re$. It suffices to prove that $d_\ell(\bar{p}_i,\bar{q}_i)\leq k \cdot {\rm{length }}({\frak C}_i)$ for every $i$, where $k$ is a universal constant.

\begin{figure}[!ht]
\centering
\includegraphics{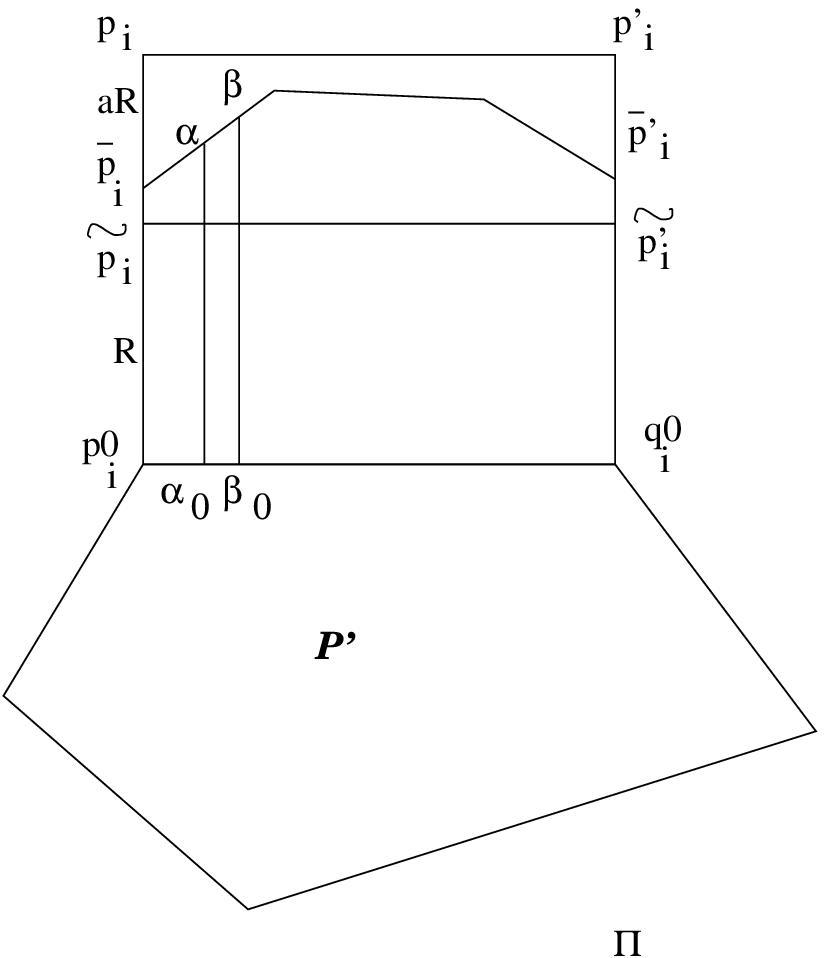} 
\caption{} 
\label{fig5} 
\end{figure}

Suppose ${\frak f}_i$ is of codimension at least $2$. Then ${\frak C}_i$ is composed of a segment $[p_i,p_i']$ and of an arc of sphere between $p_i'$ and $q_i$. Let $\bar{p}_i' = [p_i',q_i^0]\cap \partial \Re ,\; \tilde{p}_i = [p_i,p_i^0]\cap \partial {\mathcal{N}}_{R}({\mathcal{P}})$ and $\tilde{p}_i' = [p_i',q_i^0]\cap \partial {\mathcal{N}}_{R}({\mathcal{P}})$. The plane $\Pi $ determined by $p_i, p_i', p_i^0$ and $q_i^0$ intersects ${\mathcal{P}}$ in a convex polygon ${\mathcal{P}}'$ having $[p_i^0,q_i^0]$ in the boundary. The intersection $\Pi \cap \partial {\mathcal{N}}_{aR}({\mathcal{P}})$ contains the segment $[p_i,p_i']$, $\Pi \cap \partial {\mathcal{N}}_{R}({\mathcal{P}})$ contains the segment $[\tilde{p}_i,\tilde{p}_i']$ and the intersection $\Pi \cap \Re $ is a polygon containing $\bar{p}_i$ and $\bar{p}_i'$ in its boundary (Figure \ref{fig5}). 

Let $[\alpha, \beta ]$ be a segment contained by the polygonal line $\Pi \cap \partial \Re $ and also by the interior of the quadrangle $\tilde{p}_i - \tilde{p}_i' - p_i' - p_i$. Let $[\alpha_0, \beta_0 ]$ be its projection onto $[p_i^0,q_i^0]$. Since the distance from $\alpha_0$ to the line $\alpha \beta$ is at least $R$, while $d(\alpha_0, \alpha )\leq aR$, it follows that $\angle_\alpha (\alpha_0, \beta )\geq \arcsin \frac{1}{a}$, therefore that $d(\alpha, \beta )\leq ad(\alpha_0 , \beta_0 )$. This implies that $d_\ell (\bar{p}_i,\bar{p}_i')\leq ad(p_i,p_i')=ad_\ell(p_i,p_i')$.

\begin{figure}[!ht]
\centering
\includegraphics{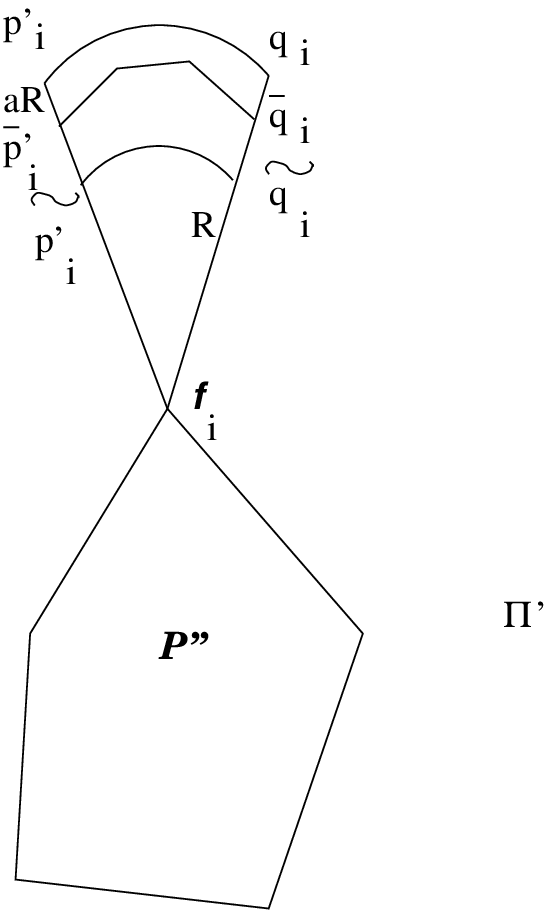}
\caption{}
\label{fig6}
\end{figure}

Let $\tilde{q}_i=[q_i^0,q_i]\cap \partial {\mathcal{N}}_{R}({\mathcal{P}})$. The plane $\Pi'$ determined by $p_i', q_i^0$ and $q_i$ intersects ${\mathcal{P}}$ also in a convex polygon. The intersections $\Pi'\cap \partial {\mathcal{N}}_{aR}({\mathcal{P}})$ and $\Pi'\cap \partial {\mathcal{N}}_{R}({\mathcal{P}})$ contain the arcs of the circles of center $q_i^0$ and radius $aR$ and $R$, respectively, joining $p_i',\, q_i$ and $\tilde{p}_i',\, \tilde{q}_i$, respectively (Figure \ref{fig6}). The fact that $d_\ell(\bar{p}_i',\bar{q}_i) \leq cd_\ell(p_i',q_i)$, with $c$ a universal constant, follows from \cite[Lemma 3.5]{Dr1*}.

We may conclude that $d_\ell(\bar{p}_i,\bar{q}_i)\leq d_\ell (\bar{p}_i,\bar{p}_i') +d_\ell(\bar{p}_i',\bar{q}_i) \leq a d_\ell(p_i,p_i') + cd_\ell(p_i',q_i) \leq k \cdot {\rm{length }}{\frak C}_i$, where $k=\max (a,c)$.

Suppose ${\frak f}_i\equiv {\mathcal{P}}$ is of codimension $1$. In this case ${\frak C}_i$ consists either of the segment $[p_i,q_i]$ or of the union of two segments with an arc of circle of length $\pi aR$. In the first situation, by looking at the plane determined by $p_i,p_i^0, q_i, q_i^0$ and repeating the first part of the previous argument we may conclude. In the second situation we split ${\frak C}_i$ into its three components. For each of the two segments we repeat the first part of the previous argument, for the arc of circle we repeat the second part of the previous argument. \hspace*{\fill  } $\diamondsuit $

\bigskip

Now we prove the key Euclidean geometry result about filling in hypersurfaces. We use the following terminology. Let $\Phi \subset\E^n$ be a linear subspace and let $F:\E^n \to \R $ be a linear form such that $\Phi \not\subset \ker F$. 

\begin{dintro}\label{strip}
We call {\it codimension one $d$-strip in $\Phi $} a set of the form $\{ x\in \Phi \mid a\leq F(x)\leq b \}$ such that its boundary hyperplanes in $\Phi $ are distance $d$ apart. We call {\it codimension two $(\epsilon ,\, d)$-strip in $\Phi $} the intersection of two codimension one  $d$-strips such that the dihedral angles between the boundary hyperplanes of one strip and the boundary hyperplanes of the other are greater than $\epsilon $.
\end{dintro}    

\begin{pintro}\label{lpolit}
Let ${\mathcal P}$ be a convex polytope in the Euclidean space $\E^n,\, n\geq 3,$ and let 
$\Phi \subsetneq \E^n$ be its affine span. Suppose that the polytope ${\mathcal P}$ has at most $m$ faces and either has codimension at least $3$ or has codimension $2$ and is contained in a codimension one  $\delta$-strip in $\Phi $ or has codimension one and is contained in a  codimension two $(\epsilon ,\, \delta )$-strip in $\Phi $.

Then there exists a 
constant $L$ depending on $m$ and $\epsilon $ such that for every $R>0$ the 
filling area of any loop 
$\frak C$ of length $\ell $ in $\partial {\mathcal{N}}_R({\mathcal P})$ is 
$$
A_1(\frak C)\leq L\cdot ( \ell^2+ \ell \delta +\delta^2 )\; .\leqno(5.5)
$$
\end{pintro}

\noindent{\bf Proof.}\quad Let $\frak C :\sph^1 \to \partial {\mathcal{N}}_R({\mathcal P})$ be a loop of length $\ell $, parametrized proportionally to the arc length and let ${\frak C}':\sph^1 \to \Phi $ be its  composition to the left with the projection on $\Phi $. It suffices to construct a filling disk of the loop $\frak C$ with the desired area outside $\breve{{\mathcal{N}}}_R({\mathcal P})$. We do this construction in the sequel. The hypothesis implies that there exists an affine subspace $\Phi_1\subset \Phi $ such that ${\mathcal P}\subset {\mathcal{N}}_{c\delta }(\Phi_1)$ and $\Phi_1$ has codimension at least $3$ in $\E^n$. We choose as $\Phi_1$ one of the extremal subspaces of the strip containing ${\mathcal P}$.

We again separate into three different cases. The first two cases roughly discuss the situation when $\frak C$ goes around ${\mathcal P}$ in $\partial {\mathcal{N}}_R({\mathcal P})$ (though there are some easier situations also treated in these cases), while the third case is easier and it corresponds to the situation when $\frak C$ ``stays on the same side of ${\mathcal P}$''.

\medskip

{\bf{(1)}}.\quad Suppose ${\frak C}'(\sph^1)\cap {\mathcal P}\neq \emptyset $. Let then $x$ be a fixed point in ${\frak C}(\sph^1)$ such that its projection $x_0$ on $\Phi $ is in ${\mathcal P}$. Let $y$ be a generic point in ${\frak C}(\sph^1)$, $y'$ its projection on $\Phi $ and $y_0$ its projection on ${\mathcal P}$.  We describe a canonical way to join $x$ to $y$ outside $\breve{{\mathcal{N}}}_R({\mathcal P})$ by a curve ${\frak C}_{xy}$ with length of order $d(x,y)+\delta $.

Let $x_1$ and $y_1$ be the respective projections of $x$ and $y$ on $\Phi_1$. Let $x_1'$ be the image of $x$ by translation of vector $\overrightarrow{x_0x_1}$.

(a)\quad Suppose $y'=y_0$. Let $y_1'$ be the image of $y$ by translation of vector $\overrightarrow{y_0x_1}$. We consider ${\frak C}_{xy}$ composed of $[x,x_1'],\; [y,y_1']$ and of the arc of minimal length joining $x_1'$ and $y_1'$ in $\Phi_{x_1}^{\perp }\cap S(x_1,R)$ if $\Phi $ has codimension at least $2$ or in $(\Phi_1)_{x_1}^{\perp }\setminus \breve{{\mathcal{N}}}_R({\mathcal P})$ if $\Phi $ has codimension $1$. Let $\beta_{xy}:=\angle_{x_1}(x_1',y_1')$. The length of the arc considered previously is $R\beta_{xy}$. It is easy to see that $\beta_{xy}$ is also equal to $\angle_{x_0}(x,y'')$ where $y''$ is the projection of $y$ on $\Phi_{x_0}^{\perp }$. Since the minimal arc joining $x$ and $y''$ in $\Phi_{x_0}^{\perp }\cap S(x_0,R)$ is also of length $R\beta_{xy}$ and since $d(x,y'')\leq d(x,y)$, it follows that $R\beta_{xy}$ is of order $d(x,y)$. Also, $d(x,x_1')\leq \delta $ and $d(y,y_1')=d(y_0,x_1)\leq d(y_0,y_1)+d(x_1,y_1)\leq \delta +d(x,y)$. Therefore in this case the length of ${\frak C}_{xy}$ is of order $\delta +d(x,y)$.

\begin{figure}[!ht]
\centering
\includegraphics{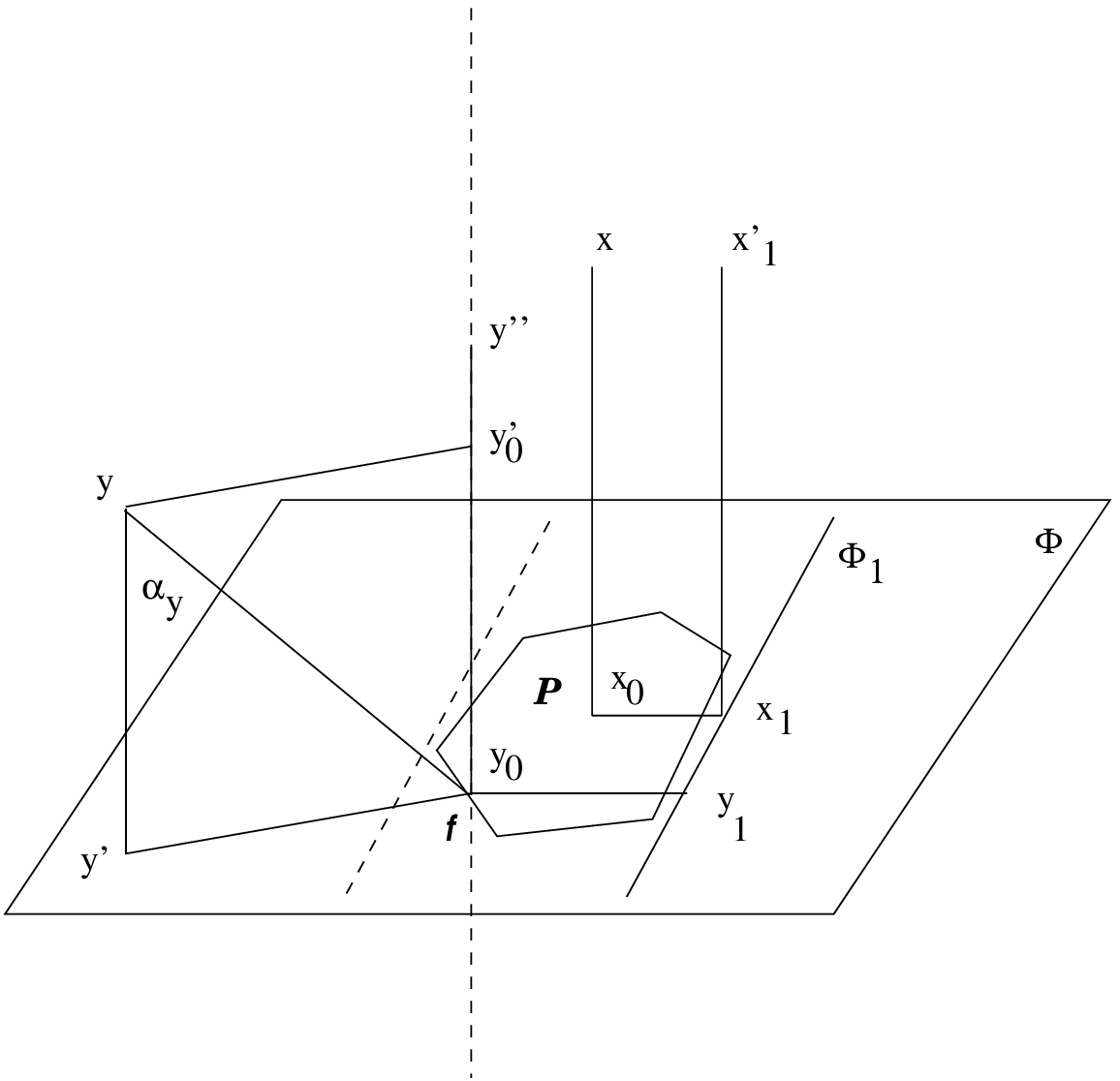} 
\caption{}
\label{fig7} 
\end{figure}

(b)\quad Suppose $y'\neq y_0$. Then $y_0$ is contained in a face $\frak f$ of ${\mathcal P}$. Since $y_0$ is also the projection of 
$y'$ on ${\mathcal P}$, $d(y',y_0)\leq d(y',x_0) \leq d(x,y)$. Let $y_0'$ be the projection of $y$ on 
$\Phi_{y_0}^{\perp }$ (Figure \ref{fig7}).

Let $\alpha_y :=\angle_{y}(y',y_0)$ if $y\neq y'$ and $\alpha_y :=\frac{\pi }{2}$ otherwise. Then $d(y', y_0)=d(y, y_0')=R\sin \alpha_y$. Let now $y''$ be the intersection point between the ray of origin $y_0$ through $y_0'$ and $\partial 
{\mathcal{N}}_R({\mathcal P})$. If $y_0=y_0'$ then we chose as $y''$ the image of $x$ by the 
translation of vector $\overrightarrow{x_0y_0}$. The arc in $\partial {\mathcal{N}}_R({\mathcal P})$ joining the points 
$y$ and $y''$, which is actually an arc on the sphere $S(y_0, R)$, is of length 
$R\alpha_y =d(y, y_0')\frac{\alpha_y }{\sin \alpha_y }=d(y', y_0)\frac{\alpha _y}{\sin 
\alpha_y }\leq cd(y',y_0)\leq cd(x,y)$ if $y\neq y'$. If $y=y'$ then the arc is of length $R\frac{\pi }{2}$ and $R\leq d(x,y)$. So to end the argument it suffices now to 
join the points $y''$ and $x$ by an arc outside $\breve{{\mathcal{N}}}_R({\mathcal P})$ with 
length of order $O(d(x,y'')+\delta )$. This can be done as in (a).

We note that if $y_1$ and $y_2$ are two points of ${\frak C}(\sph^1)$ close enough, the curves ${\frak C}_{xy_1}$ and ${\frak C}_{xy_2}$ are close. Thus, to end the proof in case {\bf{(1)}}, it suffices to consider $y_1,y_2,\dots , y_s$ on ${\frak C}(\sph^1)$ such that $x$ is between $y_1$ and $y_s$ and $x$ together with $y_1,y_2,\dots , y_s$ partition the loop into $s+1$ arcs of length at most $1$. It follows that $s=O(\ell )$. The curves ${\frak C}_{xy_1},\, {\frak C}_{xy_2}, \dots ,\, {\frak C}_{xy_s}$ form a discrete filling disk for ${\frak C}$ outside $\breve{{\mathcal{N}}}_R({\mathcal P})$. The length of each ${\frak C}_{xy_i}$ is of order $\ell +\delta $. It follows that the area of the disk is of order $\ell (\ell +\delta )$.

\medskip

{\bf{(2)}}.\quad Suppose ${\frak C}'(\sph^1)\cap {\mathcal P} = \emptyset $ and suppose there exist two points $x',y'\in {\frak C}'(\sph^1)$ such that $[x',y']\cap {\mathcal P} \neq \emptyset $.

Let $x,y \in {\frak C}(\sph^1)$ be two points whose projections on $\Phi $ are $x',y'$ respectively. Let $x_0,y_0$ be their respective projections on ${\mathcal P}$. According to Lemma \ref{simplu}, $x$ and $y$ may be joined in $\partial {\mathcal{N}}_R({\mathcal P})$ by a curve whose length is of order $d(x,y)+\inf_{t\in \partial {\mathcal P}}[d(x_0,t)+d(t,y_0)]$, so of order $\ell + \delta $. Moreover, by construction this curve contains a point whose projection on $\Phi $ coincides with the projection on ${\mathcal P} $. Each of the two arcs determined by the pair $x,y$ on $\frak C$ forms with this curve a loop with length of order $\ell +\delta $. We apply the argument in case {\bf{(1)}} to each of the two loops and we obtain filling areas of order $\ell^2+\ell \cdot \delta +\delta ^2$. Consequently we have a filling area of the same order for the initial loop $\frak C$.

\medskip

{\bf{(3)}}.\quad Suppose ${\frak C}'(\sph^1)\cap {\mathcal P} = \emptyset $ and suppose that for every two points $x',y' \in {\frak C}'(\sph^1)$, $[x',y']\cap {\mathcal P} = \emptyset $.

Let ${\frak C}_0 : \sph^1 \to  {\mathcal P}$ be the composition to the left of ${\frak C}$ with the projection on ${\mathcal P}$. Its image is entirely contained in $\partial {\mathcal P}$. By joining with segments one fixed point of ${\frak C}'(\sph^1)$ with the other points of ${\frak C}'(\sph^1)$ we can obtain a filling disk for ${\frak C}'$ with area of order $\ell^2$ and which is in the exterior of ${\mathcal P}$. The projection of this disk on ${\mathcal P}$ gives a filling disk of ${\frak C}_0$ with area of order $\ell^2$ and which is contained in $\partial {\mathcal P}$. Let $k$ be the minimal number of codimension one faces of $\partial {\mathcal P}$ needed to cover this filling disk. We prove by induction on $k$ that the filling area of ${\frak C}$ in $\partial {\mathcal{N}}_R({\mathcal P})$ is at most $c_k \ell^2$. For $k=1$ the statement follows from the cases {\bf{(1)}} and {\bf{(2)}} discussed above in which one replaces ${\mathcal P}$ with the respective codimension one face. Suppose the statement is true for $k$. In the case when the filling disk of ${\frak C}_0$ is covered by $k+1$ faces of ${\mathcal P}$, one can push the loop $\frak C$ in $\partial {\mathcal{N}}_R({\mathcal P})$ so that the trace of ${\frak C}_0$ on one face is pushed into the boundary of this face, and in order to do this an area of order $\ell^2$ is needed. This and the induction hypothesis imply the statement for $k+1$.

Since ${\mathcal P}$ has a uniformly bounded number of faces, it follows that every loop of length $\ell $ satisfying the hypothesis of this case can be filled with an area of at most $c_m \ell^2$. \hspace*{\fill }$\diamondsuit $


\end{document}